\documentclass[12pt,leqno]{amsart}

\usepackage{amssymb}
\usepackage{amsthm}
\usepackage{amsmath}
\usepackage{amscd}
\usepackage[mathscr]{eucal}
\usepackage{verbatim}

\numberwithin{equation}{section}

\theoremstyle{plain}
\newtheorem{theorem}{Theorem}[section]
\newtheorem{corollary}[theorem]{Corollary}
\newtheorem{lemma}[theorem]{Lemma}
\newtheorem{proposition}[theorem]{Proposition}

\theoremstyle{definition}
\newtheorem{definition}[theorem]{Definition}
\newtheorem{remark}[theorem]{Remark}
\newtheorem{example}[theorem]{Example}

\theoremstyle{remark}

\newcommand{\R}{\mathbb{R}}
\newcommand{\Q}{\mathbb{Q}}
\newcommand{\Z}{\mathbb{Z}}

\newcommand{\C}{\mathbb{C}}
\newcommand{\h}{\mathbb{H}}

\newcommand{\A}{\mathbb{A}}

\newcommand{\G}{\Gamma}
\newcommand{\g}{\gamma}
\newcommand{\la}{\lambda}

\newcommand{\x}{\mathbf{x}}

\newcommand{\w}{\mathbf{w}}

\newcommand{\back}{\backslash}

\newcommand{\kzxz}[4]{\left(\begin{smallmatrix} #1 & #2 \\ #3 & #4\end{smallmatrix}\right) }

\newcommand{\im}{\operatorname{Im}}

\newcommand{\wwedge}[1]{\sideset{}{^{#1}}\bigwedge}

\newcommand{\calA}{\mathcal{A}}

\newcommand{\calC}{\mathcal{C}}
\newcommand{\calD}{\mathcal{D}}

\newcommand{\calH}{\mathcal{H}}

\newcommand{\calL}{\mathcal{L}}

\newcommand{\calS}{\mathcal{S}}
\newcommand{\calT}{\mathcal{T}}
\newcommand{\calV}{\mathcal{V}}
\newcommand{\calW}{\mathcal{W}}

\newcommand{\frakg}{\mathfrak g}
\newcommand{\frakgl}{\mathfrak{gl}}

\newcommand{\frakk}{\mathfrak k}
\newcommand{\frakp}{\mathfrak p}
\newcommand{\frakn}{\mathfrak n}
\newcommand{\frakm}{\mathfrak m}
\newcommand{\frako}{\mathfrak o}
\newcommand{\fraksl}{\mathfrak{sl}}

\newcommand{\ad}{\operatorname{ad}}
\newcommand{\eps}{\varepsilon}

\newcommand{\diag}{\operatorname{diag}}

\newcommand{\tr}{\operatorname{tr}}

\newcommand{\sgn}{\operatorname{sgn}}

\newcommand{\Ker}{\operatorname{Ker}}
\newcommand{\Span}{\operatorname{span}}

\newcommand{\SL}{\operatorname{SL}}
\newcommand{\GL}{\operatorname{GL}}

\newcommand{\Sp}{\operatorname{Sp}}
\newcommand{\Symp}{\operatorname{Sp}}

\newcommand{\Mp}{\operatorname{Mp}}
\newcommand{\Orth}{\operatorname{O}}
\newcommand{\Uni}{\operatorname{U}}
\newcommand{\Hom}{\operatorname{Hom}}
\newcommand{\Sym}{\operatorname{Sym}}

\newcommand{\SO}{\operatorname{SO}}

\newcommand{\End}{\operatorname{End}}

\newcommand{\Stab}{\operatorname{Stab}}

\begin{document}

\title[Boundary behavior of special cohomology classes]
{Boundary behavior of special cohomology classes arising from the
Weil representation}

\date{\today}

\author[Jens Funke and John Millson]{Jens Funke* and John Millson**}
\thanks{* Partially supported by NSF grants DMS-0305448 and DMS-0710228}
\thanks{** Partially supported by NSF grant DMS-0405606, NSF FRG grant DMS-0554254, and the Simons Foundation}
\address{Department of Mathematical Sciences, University of Durham, Science Laboratories,
South Rd, Durham DH1 3LE, United Kingdom}
\email{jens.funke@durham.ac.uk}
\address{Department of Mathematics, University of Maryland, College Park, MD
20742, USA} \email{jjm@math.umd.edu}

\begin{abstract}

In our previous paper \cite{FMII}, we established a correspondence between vector-valued holomorphic Siegel modular forms and cohomology with local coefficients for local symmetric
spaces $X$ attached to real orthogonal groups of type $(p,q)$. This correspondence is
realized using theta functions associated to explicitly constructed ``special'' Schwartz forms. Furthermore, the theta functions give rise to generating series of certain ``special cycles'' in $X$ with coefficients.

In this paper, we study the boundary behaviour of these theta functions in the non-compact case and show that the theta functions extend to the Borel-Sere compactification $\overline{X}$ of $X$. 
However, for the $\Q$-split case for signature $(p,p)$, we have to construct and consider a slightly larger compactification, the ``big'' Borel-Serre compactification. The restriction to each face of $\overline{X}$ is again a theta series as in \cite{FMII}, now for a smaller orthogonal group and a larger coefficient system.

As application we establish the cohomological nonvanishing of the special (co)cycles when passing to an appropriate finite cover of $X$. In particular, the (co)homology groups in question do not vanish.

\end{abstract}

\maketitle

\section{Introduction}

The cohomology of arithmetic quotients $X = \G \back D$ of a symmetric space $D$ associated to a reductive Lie group $G$ is of fundamental interest in number theory and for the field of automorphic forms.  For dual reductive pairs, one can apply the ``geometric theta correspondence'' (see below) obtained by the Weil representation to construct cohomology classes on locally symmetric spaces associated to these groups.  
One very attractive aspect of this method is that the classes obtained in this way often give rise to Poincar\'e dual forms for geometrically defined, ``special'' cycles arising via the embedding $H \hookrightarrow G$ of suitable subgroups $H$. 
\medskip

Let $V$ be a rational quadratic space of signature $(p,q)$ with for simplicity even dimension $m$. Let $\underline{G}= \SO(V)$ and let $G = \underline{G}(\R)_0= \SO_0(V_{\R})$. Let $D_V=D=G/K$ be the symmetric space of $G$ of dimension $pq$ with $K$ a maximal compact subgroup. We
let $\mathfrak{g} = \mathfrak{k} \oplus \mathfrak{p}$ be the associated
Cartan decomposition of the Lie algebra of $G$.

Every partition $\la$ of a non-negative integer $\ell'$ into at most
$n$ parts gives rise to a dominant weight $\la$ of $\GL(n)$. We write
$i(\la)$ for the number of nonzero entries of $\la$. We explicitly realize the corresponding irreducible representation of highest weight $\la$ as the image $\mathbb{S}_{\la}(\C^n)$ of the
Schur functor $\mathbb{S}_{\la}(\cdot)$ associated to $\la$ applied to the tensor
space $T^{\ell'}(\C^n)$. We can apply the same Schur functor to
$T^{\ell'}(V_{\C})$ to obtain the space $\mathbb{S}_{\la}(V_{\C})$, and the harmonic $\ell'$-tensors in $\mathbb{S}_{\la}(V_{\C})$ give the irreducible representation $\mathbb{S}_{[\la]}(V_{\C})$ for $G$ with highest weight $\widetilde{\la}$ (under some restrictions). If $i(\la) \leq [\tfrac{m}2]$, then $\tilde{\la}$ has the same nonzero entries as $\la$ (when $\tilde{\la}$ is expressed in coordinates relative to the standard basis $\{\epsilon_i\}$ of \cite{Bourbaki}, Planche II and IV). 

The Weil representation induces an action of $\Symp_n(\R) \times \Orth(V_\R)$ on $\mathcal{S}(V^n_\R)$, the Schwartz functions on $V^n_\R$. The main point of our previous paper \cite{FMII} is the construction of certain $(\mathfrak{g},K)$-cocycles
\begin{equation}
\varphi^V_{nq,[\la]} \in \left[ \wwedge{nq} (\mathfrak{p}^{\ast}_{\C})
 \otimes \mathcal{S}(V^n_\R)\otimes \mathbb{S}_{[\la]}(V_{\C})\right]^{K} \notag
\end{equation}
with values in $\mathcal{S}(V^n_\R) \otimes  \mathbb{S}_{[\la]}(V_{\C})$. These classes generalize the work of Kudla and Millson (e.g. \cite{KM90}) to the case of nontrivial coefficients systems $\mathbb{S}_{[\la]}(V_{\C})$. The cocycle $\varphi^V_{nq,[\la]}$ corresponds to a closed differential $nq$-form $\tilde{\varphi}^V_{nq,[\la]}$ on $D$ with values in $\mathcal{S}(V_\R^n) \otimes  \mathbb{S}_{[\la]}(V_{\C})$. For a coset of a lattice $\mathcal{L}$ in $V^n$, we define the theta distribution $\Theta_{\mathcal{L}}= \sum_{\ell \in \mathcal{L}} \delta_{\ell}$, where $ \delta_{\ell}$ is the delta measure concentrated at $\ell$. It is obvious that $\Theta_{\calL}$ is invariant under $\Stab(\mathcal{L}) \subset G$. 
Hence we can apply the theta distribution to $\tilde{\varphi}^V_{nq,[\la]}$ to obtain a closed $nq$-form $\theta_{\varphi^V_{nq,[\la]}}$ with values in (the local system associated to) $ \mathbb{S}_{[\la]}(V_{\C})$ on the finite volume quotient $X = \Gamma \backslash D$ given by 
\[
\theta_{\varphi^V_{nq,[\la]}}(\mathcal{L})= \langle \Theta_{\mathcal{L}}, \tilde{\varphi}^V_{nq,[\la]} \rangle.
\]
Here $\G \subseteq \Stab(\mathcal{L})$ is a congruence subgroup. Furthermore, it is shown in \cite{FMII} that $\theta_{\varphi^V_{nq,[\la]}}$ also gives rise to a non-holomorphic vector-valued Siegel modular form for the representation $S_{\la}(\C^n) \otimes \det^{m/2}$ on the Siegel space $\h_n$. We may then use $\theta_{\varphi^V_{nq,[\la]}}$ as the integral kernel of a pairing of Siegel modular forms $f$ with (closed) differential $(p-n)q$-forms $\eta$ or $nq$-chains (cycles) $C$ in $X$. The resulting pairing in $f$, $\eta$ (or $C$), and (possibly different) Schwartz cocycles $\varphi$, we call the {\it geometric theta correspondence}. 

Special cycles $Z_U$ arise from the embedding $G_U \hookrightarrow G$ of the stabilizer of a positive definite rational subspace $U \subset V$ of dimension $n$. Hence $G_U$ is an orthogonal group of signature $(p-n,q)$. The special cycles $Z_U$ for varying $U$ give rise to a family of composite cycles $Z_T$ parametrized by symmetric positive definite integral $n\times n $ matrices $T$. We obtain (by Poincar\'e duality) classes $[Z_T]$ in $H^{nq}(X,\Z)$, and in \cite{FMII} we explain how to attach $\mathbb{S}_{[\la]}(V_{\C})$-coefficients to the cycles to obtain classes 
\begin{equation}\label{FMcoeffcycles}
[Z_{T,[\la]}]  \in  \mathbb{S}_{\la} (\C^n)^{\ast} \otimes H^{nq}(X, \mathbb{S}_{[\la]}(V_{\C})).
\end{equation}
Then the main result in \cite{FMII} is that 
\begin{equation}\label{FMcoeffresult}
[\theta_{\varphi^V_{nq,[\la]}}]= \sum_{T \geq 0} [Z_{T,[\la]}] e^{2\pi i \tr(T\tau)}
\end{equation}
is a holomorphic vector-valued Siegel modular form with values in $ H^{nq}(X, \mathbb{S}_{[\la]}(V_{\C})$. Here $\tau \in \h_n$. (We omit the definition of $ [Z_{T,[\la]}]$ for $T$ semi-definite). This result gives further justification to the term geometric theta correspondence. 

\medskip

Recently, it has now been shown \cite{BMM} for all $\SO(p,q)$ with $p+q>6$ and $p \geq q$ in the cocompact  case that the geometric theta correspondence specialized to $\varphi^V_{q,[\la]}$ ($n=1$) induces on the adelic level an {\it isomorphism} from the appropriate space of classical modular forms to the complement of the space spanned by invariant forms in the direct limit of the cohomology groups $H^{q}(X, \mathbb{S}_{[\la]}(V_{\C}))$. In particular, for any congruence quotient, the cohomology groups $H^{q}(X, \mathbb{S}_{[\la]}(V_{\C}))$ are spanned by Poincar\'e duals of cycles and invariant forms. For $X$ finite volume one to change to cuspidal cohomology. (Their result for $n>1$ is more difficult to state). 
This result highlights the importance of the cohomology classes constructed via the Weil representation. 

\medskip

It is therefore a very natural question to study $\theta_{\varphi^V_{nq,[\la]}}$ for non-compact $X$, in particular to analyze its boundary behavior. This is what we do in this paper. 

\medskip

We let $P=\underline{P}(\R)_0$ be the connected component of the
identity of the real points of a rational parabolic subgroup
$\underline{P}$ in $\underline{G}$ stabilizing a flag ${\bf F}$ of
totally isotropic rational subspaces in $V$. Conversely, for signature different than $(p,p)$ all such flags give rise to a unique rational parabolic. Then the Borel-Serre compactification $\overline{X}$ compactifies $X$ by adding to each rational $P$ a face $e'(P)$, which is a nilmanifold bundle over a suitable quotient of the symmetric space associated to the semi-simple part of the Levi subgroup of $P$, see \cite{BS,BJ}. This makes $\overline{X}$ a manifold with corners. 

However, for the $\Q$-split case in signature $(p,p)$, the rational parabolics are {not} in 1-1 correspondence with the stabilizers of rational totally isotropic flags in $V$ (but rather of so-called oriflammes). This turns out to be a critical issue for us. 
To remedy this we consider instead the spherical building of proper rational parabolic subgroups for the full (non-connected) orthogonal group $\Orth(p,p)$ instead. The space $X$ does not change, but now isotropic flags do parameterize the parabolics. The resulting compactification we call the {\it big} Borel-Serre compactification of $X$ which turns out to be (slightly) bigger and denote by abuse of notation also by $\overline{X}$. For an alternative construction of the big $\overline{X}$, we embed $X=X_{p,p}$ into a locally symmetric space $X_{p+1,p}$ for signature $(p+1,p)$ and then consider the closure of $X_{p,p}$ in $\overline{X}_{p+1,p}$. 

To illustrate the big Borel-Serre compactification, we consider the split case for $\SO(2,2)$, when $X=X_1 \times X_2$ is the product of two modular curves. Then the Borel-Serre compactification of $X$ is the product of the two individual compactifications $\overline{X}_1 \times \overline{X}_2$ which adds to each cusp of the modular curves a circle $S_1$. Hence the corner at the cusp $(z_1,z_2)=(i\infty,i\infty)$ of $X$ is given by a $2$-torus $T^2$. Then the big Borel-Serre compactification of $X$ blows up the corner to $T^2 \times \R_+$ with the new coordinate $Im(z_1)/Im(z_2) \in \R_+$ measuring the ``slope'' by which one enters the corner from the interior. 
We explain the details of the big Borel-Serre compactification in section~\ref{bigBS}.

Let $E$ be the largest element in the rational isotropic flag ${\bf F}$ with dimension $\ell$ corresponding to $\underline{P}$. Set $W = E^{\perp} /E$, which is naturally a quadratic space of signature $(p-\ell,q-\ell)$. Then a suitable arithmetic quotient $X_W$ of the symmetric space $D_W$ associated to $W$ occurs as a factor in the base of the nilmanifold bundle $e'(P)$. 

The main result of this paper is

\begin{theorem}\label{globalrestrictionformula}  
\begin{enumerate}
\item The form $\theta_{\mathcal{L}_V}(\varphi^V_{nq,[\lambda]})$ extends to a smooth differential form on
the (big) Borel-Serre compactification $\overline{X}$ considered as a smooth manifold with corners.  In fact, the form
$\theta _{\mathcal{L}_V}(\varphi^V_{nq,[\lambda]}))$ is the sum of a rapidly decreasing differential form and  a special differential form in the sense of \cite{GHMP}, p.169.

\item  For a given face $e'(P)$, let $\tilde{r}_{P}$ be the restriction map from $\overline{X}$ to $e'({P})$. Then there exists a theta distribution $\widehat{\mathcal{L}}_W$ for $W$ such that
\[
[\tilde{r}_{{P}}(\theta _{\mathcal{L}_V}(\varphi^V_{nq,[\lambda]}))] = [\widetilde{\iota}_P (\theta_{\widehat{\mathcal{L}}_W}(\varphi^W_{n(q-\ell),[\ell\varpi_n+\la]}))].
\]
Here $\widetilde{\iota}_P$ is an embedding
\begin{align}
\widetilde{\iota}_P: H^{n(q-\ell)}(X_W,\mathbb{S}_{[\ell\varpi_n+\la]}(W_{\C})) &\hookrightarrow H^{nq}(e'(P),\mathbb{S}_{[\la]}(V_{\C})), \notag
\end{align} 
where $\varpi_n = (1,\dots,1)$ is the $n$-th fundamental weight for $\GL(n)$, so that the Young diagram associated to $\ell \varpi_n$ is an $n$ by $\ell$ rectangle.

\end{enumerate}

\smallskip

In particular, $[\tilde{r}_{{P}}(\theta _{\mathcal{L}_V}(\varphi^V_{nq,[\lambda]}))] =0$ for $n>\min\left(p,\left[\tfrac{m}{2}\right]\right)-\ell$ (if $\ell \geq 2$) and $n> p-1$ or $n>m-2-i(\la)$ (if $\ell=1$).

\end{theorem}

Loosely speaking Theorem~\ref{globalrestrictionformula}
can be summarized by saying that the restriction of our theta series for
$\SO(V)$ to a face of $\overline{X}$ is the
theta series for $\SO(W)$ of the same type corresponding to an enlarged coefficient
system corresponding to placing an $n$ by $\ell$ rectangle on the
left of the Young diagram corresponding to $\lambda$ to obtain a
bigger Young diagram corresponding to $\ell \varpi_n + \la$. The theta series $\theta_{\mathcal{L}_V}(\varphi^V_{nq,[\lambda]})$ is termwise moderately increasing, so the statement of the theorem is rather delicate. To capture the boundary behavior we switch to a mixed model of the Weil representation.

We can also interpret our result in terms of weighted cohomology \cite{GHMP}, see Remark~\ref{weightedremark}. 

As stated above, for the split $\SO(p,p)$-case, the differential form $\theta _{\mathcal{L}_V}(\varphi^V_{nq,[\lambda]}))$ does {\it not} extend to the usual Borel-Serre boundary. 

\medskip

{\it Non-vanishing at the boundary}. 
As an easy and direct application we obtain a non-vanishing result for the special (co)cycles.

\begin{theorem}\label{TH:nonvan}
Assume that the $\Q$-rank and the $\R$-rank of $\underline{G}$ coincide. Then for
\begin{equation*}\label{nonvan-crit}
 i(\la) \leq n \leq 
\begin{cases}  
\left[ \frac{p-q}{2} \right] & \; \text{if $q\geq2$} \\ 
  p-1 -i(\la) & \text{ if $q=1$}, 
\end{cases}
\end{equation*}
there exists a finite cover $X'$ of $X$ such that
\[
[\theta(\varphi^V_{nq,[\lambda]})]  \ne 0.
\]
Using \eqref{FMcoeffresult} this gives $[Z_{T,[\la]}]  \ne 0$ for infinitely many $T$. In particular,
\[
H^{nq}(X',\mathbb{S}_{[\la]}(V_{\C})) \ne 0.
\]
Finally, $H^{nq}(X',\mathbb{S}_{[\la]}(V_{\C}))$ is not spanned by classes given by invariant forms on $D$. 
\end{theorem}

The basic idea for the proof is to study the restriction to a face of  $\overline{X}$ associated to a minimal rational parabolic subgroup. At such a face, the space $W$ is positive definite, and hence the restriction becomes a positive definite theta series for which we establish non-vanishing.  

There are numerous non-vanishing results in the literature, and we mention a few related ones. In the case of nontrivial coefficients for compact hyperbolic manifolds, Millson \cite{M-Raghu} proved the nonvanishing of the special cycles with coefficients in codimension $n$ in the range $i(\la) \leq  n\leq p -i(\la)$. Bergeron \cite{Bergeron} in the compact case established non-vanishing of the classes introduced by Kudla and Millson (trivial coefficients) by considering the analogous classes in $\Uni(p,q)$. Li \cite{Li}  also used the theta correspondence to establish non-vanishing for the cohomology of orthogonal groups, again in the compact (or $L^2$)-case (without giving a geometric interpretation of the classes). Speh and Venkataramana \cite{SV} gave in general a criterion for the non-vanishing of certain modular symbols  in terms of the compact dual. In contrast to our result, their non-vanishing occurs from classes defined by invariant forms on $D$. 

\medskip

{\it Vanishing at the boundary}.
We first describe the general main motivation for our work. From \eqref{FMcoeffcycles} and \eqref{FMcoeffresult} we see that theta series $\theta_{\mathcal{L}}(\varphi_{nq,0})$ (for simplicity, we only consider trivial coefficients for the moment) gives rise to a map
\begin{equation}\label{KMlift}
\Lambda_{nq}: H_c^{(p-n)q}(X,\C) \to M_{m/2}^{(n)}(\G')
\end{equation}
from the cohomology with compact supports to the space of holomorphic Siegel modular forms of degree $n$ of weight $m/2$. We are interested in extending the lift \eqref{KMlift} (also for non-trivial coefficients) to other cohomology groups of the space $X$ which capture its boundary. This paper should be considered in this context, and is central to our efforts. This program is in particular motivated by the work of Hirzebruch-Zagier \cite{HZ}, which is the $\Q$-rank $1$ case for signature $(2,2)$ when $X$ is a Hilbert modular surface and the cycles in question are the famous Hirzebruch-Zagier curves ($n=1$). 

Whenever the restriction of $\theta_{\mathcal{L}}(\varphi_{nq,0})$ to $\partial X$ is cohomologically trivial, then such an extension exists. Namely, in this case, one can utilize a mapping cone construction to modify $\theta_{\mathcal{L}}(\varphi_{nq,0})$ to represent a class in the compactly supported cohomology of $X$ - in principal. The main problem is to explicitly construct suitable primitives for the restriction (again using the theta correspondence). Then one obtains an extension of $\Lambda_{nq}$ to the full cohomology of $X$. 

We have already carried this out in several instances. First and foremost, the restriction vanishes in the Hirzebruch-Zagier case, and based on this, we give in \cite{FM-HZ} a new treatment and extension of the results in \cite{HZ} using the theta correspondence. The $\Q$-rank 2 case when $X$ is the product of two modular curves is of course highly interesting as well. Now the boundary faces in the big Borel-Serre compactification are no longer isolated, and in addition some subtle analytic complications arise when constructing the primitives at the boundary. We consider this case in the near future. 

The case which resembles Hirzebruch-Zagier most closely is the one for Picard modular surfaces (quotients of $\Uni(2,1)$; the results of this paper generalize to unitary groups). Cogdell \cite{Cogdell} considered this case in the spirit of Hirzebruch-Zagier. We will consider this case from our point of view also in a subsequent paper.

Another case is $\SO(2,1)$ when $X$ is a modular curve, and the cycles are geodesics. For {\it non-trivial} coefficients, the restriction to the boundary vanishes. This case is particularly attractive since one can interpret our classes as (co)homology classes for even powers of the universal elliptic curve.  We discuss this case in detail in \cite{FMspectacle}.

\medskip

Finally, we mention that \cite{FKyoto} gives an introductory survey of the results obtained in this paper.

\medskip

We would like to thank G. Gotsbacher, L. Saper, and J. Schwermer for fruitful discussions and E. Freitag and R. Schulze-Pillot for answering a question on positive definite theta series. As always it is a pleasure to thank S. Kudla for his encouragement. The work on this paper has greatly benefitted from three visits of the first named author at the Max Planck Institute from 2005 to 2008. He gratefully acknowledges the excellent research environment in Bonn.

\tableofcontents

\section{Basic Notations}\label{basicstuff}

\subsection{Orthogonal Symmetric Spaces}

Let $V$ be a rational vector space of dimension $m= p+q$ and let
$(\;,\;)$ be a non-degenerate symmetric bilinear form on $V$ with
signature $(p,q)$. 
We fix a standard orthogonal basis $e_1,\dots,e_p, e_{p+1},\dots, e_m$ of
$V_{\R}$ such that $(e_{\alpha},e_{\alpha})=1$ for $1 \leq \alpha
\leq p$ and $(e_{\mu},e_{\mu})=-1$ for $p+1 \leq \mu \leq m$. (We
will use "early" Greek letters to denote indices between $1$ and
$p$, and "late" ones for indices between $p+1$ and $m$). With
respect to this basis the matrix of the bilinear form is given by
the matrix $I_{p,q}= \left(
\begin{smallmatrix}
1_p & \\
    &-1_q
\end{smallmatrix} \right)$.

We let $\underline{G} = \SO(V)$ viewed as an algebraic group over
$\Q$. We let $G:=G(\R)_0$ be the connected component of the identity
of $G(\R)$ so that $ G \simeq \SO_0(p,q)$. We let $K$ be the maximal
compact subgroup of $G$ stabilizing $\Span\{e_{\alpha}; 1 \leq \alpha \leq p\}$. Thus $K \simeq \SO(p) \times \SO(q)$. Let $D = G/K$ be the symmetric space
of dimension $pq$ associated to $G$. We realize $D$ as the space of
negative $q$-planes in $V_{\R}$:
\begin{equation}
D \simeq \{z \subset V_{\R}: \dim z =q; \, (\;,\;)|_z <0\}.
\end{equation}
Thus $z_0 =  \Span\{e_{\mu}; p+1 \leq \mu \leq m \}$ is the base point of $D$. Furthermore, we can also
interpret $D$ as the space of minimal majorants for $(\,,\,)$. That
is, $z \in D$ defines a majorant $(\,,\,)_z$ by $(x,x)_z = -(x,x)$
if $x \in z$ and $(x,x)_z = (x,x)$ if $x \in z^{\perp}$. We write
$(\,,\,)_0$ for the majorant associated to the base point $z_0$.

The Cartan involution $\theta_0$ of $G$ corresponding to the
basepoint $z_0$ is obtained by conjugation by the matrix $I_{p,q}$.
We will systematically abuse notation below and write $\theta_0(v)$
for the action of the linear transformation of $V$ with matrix
$I_{p,q}$ relative to the above basis acting on $v \in V$. Let
$\mathfrak{g}$ be the Lie algebra of $G$ and
$\mathfrak{k}$ be the one of $K$. We obtain the Cartan decomposition
\begin{equation}
\mathfrak{g} = \mathfrak{k} \oplus \mathfrak{p},
\end{equation}
where
\begin{equation}
\mathfrak{p} =\Span \{X_{\alpha\mu}:=e_{\alpha} \wedge {e_{\mu}}; 1
\leq \alpha \leq p, p+1 \leq \mu \leq m\}.
\end{equation}
Here $w \wedge w' \in \bigwedge^2 V_{\R}$ is identified with an
element of $\mathfrak{g}$ via
\begin{equation}
(w \wedge w')(v) = (w,v)w' - (w',v)w.
\end{equation}
We let $\{\omega_{\alpha\mu}\}$ be the dual basis of
$\mathfrak{p}^{\ast}$ corresponding to $\{D_{\alpha\mu}\}$. Finally
note that we can identify $\mathfrak{p}$ with the tangent space
$T_{z_0}(D)$ at the base point $z_0$ of $D$.

We let $r$ be the Witt rank of $V$, i.e., the dimension
of a maximal totally isotropic subspace of $V$ over $\Q$ and assume $r>0$. Let $F$ be an isotropic subspace of $V$ of dimension $\ell$. Then we
can describe the $\ell$-dimensional isotropic subspace $\theta_0(F)$
as follows. For $U$ a subspace of $V$, let $U^{\perp}$, resp.
$U^{\perp_0}$
 be the orthogonal complement of $U$
for the form $(\;,\;)$, resp. $(\;,\;)_0$. Then $\theta_0(F) = (F^{\perp})^{\perp_0}$. We fix a maximal totally isotropic subspace $E_r$ and choose a basis
$u_1,u_2,\dots,u_r$ of $E_r$. Let $E_r' = \theta_0(E_r)$.
 We pick a basis
$u'_{r},\cdots,u'_1$ of $E_r'$ such that $(u_i,u'_j) = \delta_{ij}$. More generally, we let
\begin{equation}\label{isospace}
E_{\ell} := \Span\{u_1,\dots,u_{\ell}\},
\end{equation}
and we call $E_{\ell}$ a \emph{standard} totally isotropic subspace.
Furthermore, we set $E'_{\ell}= \theta_0(E_{\ell}) =
\Span(u'_{\ell},\dots,u'_1)$. Note that $E_{\ell}'$ can be naturally
identified with the dual space of $E_{\ell}$. We can assume that
with respect to the standard basis of $V_{\R}$ we have $e_{\alpha} = \tfrac{1}{\sqrt{2}}(u_{\alpha} + u'_{\alpha})$
and $e_{m+1-\alpha} = \tfrac{1}{\sqrt{2}}( u_{\alpha} - u'_{\alpha})$ for $\alpha=1,\dots,\ell$. We let
\begin{equation}
W_{\ell} = E_{\ell}^{\perp}/E_{\ell},
\end{equation}
and note that $W_{\ell}$ is a non-degenerate space of signature
$(p-\ell,q-\ell)$. We can realize $W_{\ell}$ as a subspace of $V$ by
\begin{equation}
W_{\ell} = (E_{\ell}\oplus E_{\ell}')^{\perp},
\end{equation}
where the orthogonal complement is either with respect to $(\,,\,)$
or $(\,,\,)_0$. This gives 
\begin{equation}
V = E_{\ell} \oplus W_{\ell} \oplus E_{\ell}',
\end{equation}
a $\theta_0$-invariant Witt splitting for $V$. Note that with these choices $\theta_0$ restricts to a Cartan
involution for $\Orth(W_{\ell})$. We obtain a Witt basis
$u_1,\dots,u_{\ell},e_{\ell+1},...,e_{m-\ell},u'_{\ell},\dots,u'_{1}$
for $V_{\R}$. We will denote coordinates with respect to the Witt
basis with $y_{i}$ and coordinates with respect to the standard
basis with $x_{i}$. 

We often drop the subscript $\ell$ and just write $E$, $E'$, and $W$. 

\subsection{Parabolic Subgroups}\label{parabolicsubgroups}

We describe the rational parabolic subgroups of ${G}$. 

\subsubsection{Isotropic flags and parabolic subgroups}

 We let $\mathbf{F}$ be a flag of totally
isotropic subspaces $F_1 \subset F_2 \subset \dots
\subset F_k$ of $V$ over $\Q$. Then we let $\underline{P}=
\underline{P}_{\mathbf{F}}$ be the parabolic subgroup of
$\underline{G}$ stabilizing the flag $\mathbf{F}$:
\begin{equation}
\underline{P}_{\mathbf{F}} = \{ g \in \underline{G} ; g F_i = F_i\},
\end{equation}
and write $P= P_{\bf F} = (\underline{P}_{\mathbf{F}}(\R)) _0$ for the resulting rational parabolic in $G$. The first fundamental fact is 

\begin{lemma}
Assume that $V$ is not a rational $\Q$-split space of signature $(p,p)$. Then the assignment $
\mathbf{F} \mapsto \underline{P}_{\mathbf{F}}$ defines a bijection between the rational totally isotropic flags in $V$ and rational parabolic subgroups in ${G}$. Furthermore, under this map isotropic subspaces give rise to maximal parabolics. 
\end{lemma}

In this situation, we can assume by conjugation that the flag $\mathbf{F}$ consists of
standard totally isotropic subspaces $E_i$ \eqref{isospace} and call such parabolics
a \emph{standard} $\Q$-parabolic. 

However, if $V$ is a rational $\Q$-split space of signature $(p,p)$ then the map from totally isotropic flags to parabolics is surjective but not 1-1. We need a more involved incidence relation between totally isotropic subspaces than inclusion to describe parabolic subgroups which gives rise to a configuration called oriflammes, see eg \cite{Garrett}, chapter 11. 

\begin{definition}(Oriflammes)
We define the incidence relation $\sim$ on non-zero totally isotropic subspaces of $V$ of dimension different than $p-1$ by $F_1 \sim F_2$ if either 

(i) $F_1 \subset F_2$ or $F_2 \subset F_1$, or 

(ii) If $\dim F_1= \dim F_2 =p$, then $F_1 \cap F_2$ has dimension $p-1$. 

\noindent Then an oriflamme is a collection of such subspaces in which any two members are incident. 
\end{definition}

One then has (see eg \cite{AB,Garrett}, also Example~\ref{p1flag2})

\begin{lemma}\label{oriflamme}
Assume that $V$ is a rational $\Q$-split space of signature $(p,p)$. Then the rational parabolic subgroups in ${G}$ are in 1-1 correspondence with the rational oriflammes in $V$ by taking the stabilizer of the oriflamme. Concretely, 

(1) The maximal parabolics are attached to totally isotropic subspaces of dimension different than $p-1$. The totally isotropic subspaces of dimension $p-1$ do not give rise to a maximal parabolic. 

(2) All totally isotropic flags which do {\it not} include a constituent of dimension $p-1$ gives rise to different standard parabolic subgroups. 

(3) Let $F_{p-1}$ be a totally isotropic space of dimension $p-1$ and $\mathbf{F} = F_1 \subset F_2 \subset \dots F_k \subset F_{p-1}$ be a totally isotropic flag. Since $F_{p-1}^{\perp}/F_{p-1}$ is naturally a $\Q$-split space of signature $(1,1)$ there are exactly two totally isotropic spaces $F_{p,1}$, $F_{p,2}$ of (maximal) dimension $p$ which contain $F_{p-1}$. Then the three flags $\mathbf{F}$, $\mathbf{F} \subset F_{p,1}$, $\mathbf{F}  \subset F_{p,2}$, are fixed by the same parabolic in $G$. This parabolic fixes the oriflamme $(F_1, F_2 ,\dots, F_{k},F_{p,1},F_{p,2})$. 
\end{lemma}

Let $E_+ = E_p=\Span(u_1,\dots,u_{p-1}, u_p)$ and $E_- = \Span(u_1,\dots u_{p-1}, u_p')$. Then we define the {standard} $\Q$-parabolics to be the ones given by fixing a suboriflame of the maximal oriflamme $(E_1,E_2, \dots, E_{n-2}, E_+, E_-)$. We discuss the case when $V$ is a rational $\Q$-split space of signature $(p,p)$ in more detail in section~\ref{bigBS}.

\subsubsection{The Langlands decomposition}

We let $\underline{N}_{\underline{P}}$ be the unipotent radical of
$\underline{P}$. It acts trivially on all quotients of the flag.
We let $\underline{L}_{\underline{P}} =
\underline{N}_{\underline{P}} \back \underline{P}$ and let
$\underline{S}_{\underline{P}}$ be the split center of
$\underline{L}_{\underline{P}}$ over $\Q$. Note that
$\underline{S}_{\underline{P}}$ acts by scalars on each quotient.
Let $\underline{M}_{\underline{P}} = \cap_{\chi \in
X(\underline{L}_{\underline{P}})} \Ker(\chi^2)$. We let $N=N_P$ and
$L=L_{\underline{P}}$ be their respective real points in $G$, and as before we set
 $M = M_P= (\underline{M}_{\underline{P}}(\R))_0$, and
$A=A_{\underline{P}} = (\underline{S}_{\underline{P}}(\R))_0$. We can realize
$\underline{L}_{\underline{P}}$ (and also
$\underline{S}_{\underline{P}}, \underline{M}_{\underline{P}}$) as
$\theta_0$-stable subgroups of $\underline{P}$:
\begin{equation}
\underline{L}_{\underline{P}} = \underline{P} \cap
\theta_0(\underline{P}).
\end{equation}
Then $\underline{M}_{\underline{P}}$ is the semi-simple part of the
centralizer of $\underline{S}_{\underline{P}}$ in $\underline{P}$.
We will regularly drop the subscripts $\mathbf{F}$, $\underline{P}$,
and $P$.

We obtain the (rational) Langlands decomposition of $P$:
\begin{equation}
P= N A M \simeq N \times A \times M,
\end{equation}
and we write $\mathfrak{n}$, $\mathfrak{a}$, and $\mathfrak{m}$ for
their respective Lie algebras. The map $P \to N \times A \times M$
is equivariant with the $P$-action defined by
\begin{equation}
n'a'm'(n,a,m)  = \left(n'Ad(a'm')(n),a'a,m'm \right).
\end{equation}

\subsubsection{The Levi}

We let $\mathbf{F}$ be a standard rational totally isotropic flag $0 =E_{0} \subset 
E_{i_1} \subset \cdots \subset E_{i_k}= E_{\ell}=E$ and assume that
the last (biggest) totally isotropic space in the flag $\mathbf{F}$
is equal to $E_{\ell}$ for some $\ell$. The reader will make the necessary adjustments when considering an oriflamme in the $\Q$-split $\SO(p,p)$-case.

Let $U_{i_j}= \Span(u_{i_{j-1}+1},\dots, u_{i_{j}})$ be the orthogonal complement of $E_{i_{j-1}}$ in
$E_{i_{j}}$ with respect to $(\,,\,)_0$ and $U_{i_j}'$ be the orthogonal complement of
$E_{i_j}'$ in $E_{i_{j+1}}'$ and let $W =W_{\ell}= (E_{\ell} \oplus
E_{\ell}')^{\perp}$. We obtain a refinement of the Witt
decomposition of $V$ such that the subspaces $U_{i_j},U_{i_s}'$,
and $W$ are mutually orthogonal for $(\;,\;)_0$ and defined over
$\Q$:
\begin{equation}
V = \left( \bigoplus_{i_j=1}^k U_{i_j} \right) \oplus W \oplus
\left( \bigoplus_{i_j=1}^k U_{i_j}'\right).
\end{equation}
Then $\underline{L}_P$ is the subgroup of $\underline{P}$ that stabilizes each
of the subspaces in the above decomposition of $V$. In what follows
we will describe matrices in block form relative to the above direct
sum decomposition of $V$. We first note that we naturally have $\Orth(W) \times \GL(E) \subset
\Orth(V)$ via
\begin{equation}\label{O-GL-embed}
\left\{ \left( \begin{smallmatrix}
g &  &   \\
  &h &   \\
  &  & \tilde{g}
\end{smallmatrix} \right)
; h \in \Orth(W), g \in \GL(E) \right\},
\end{equation}
 where $\tilde{g} = J g^{\ast} J$, $g^{\ast} = {^t g
 ^{-1}}$, and  $J = \left( \begin{smallmatrix} & & 1 \\ & \cdots & \\
1 & &   \end{smallmatrix}\right)$. In particular, we can view the corresponding Lie algebras
 $\frako(W_{\R})$ and $\frakgl(E_{\R})$ as subalgebras of $\frakg$.
 Namely,
\begin{align}
\frako(W_{\R}) &\simeq \Span \{ e_i \wedge e_j; \, \ell < i < j \leq
m-\ell \}, \\
\frakgl(E_{\R}) &\simeq \Span \{ u'_i \wedge u_j; \, i,j \leq \ell
\}, 
\end{align}
via $\mathfrak{g} \simeq \bigwedge^2 V_{\R}$. We see
\begin{equation}\label{M-embedding}
\underline{L} \simeq \left\{ \left( \begin{smallmatrix}
g &  &   \\
  &h &   \\
  &  & \tilde{g}
\end{smallmatrix} \right)
; h \in \SO(W), g =\diag(g_1, \dots,g_k) \in \prod_{j=1}^{k}
\GL(U_{i_j}),  \right\}.
\end{equation}

We now consider the isotropic flag ${\bf F}$ in $V$ as a flag ${\bf
F}(E)$ of subspaces inside
$E$. We let $\underline{P}'$ be the parabolic subgroup of $\GL(E)$
stabilizing ${\bf F}(E)$. Then for the real points
$P'=(\underline{P'}(\R))_0$, we have
\begin{equation}
P' = N_{P'} A M_{P'},
\end{equation}
with unipotent radical $N_{P'}$ and Levi factor
\begin{equation}
M_{P'} = \prod_{j=1}^k \SL(U_{i_j}(\R)).
\end{equation}
Here $A$ is as above, viewed as a subgroup of $\GL_+(E_{\R})$. Furthermore, we can view $\underline{P}'$ and its subgroups naturally as subgroups of of $\underline{P}$ via the
embbeding of $\GL(E)$ into $\Orth(V)$ given by \eqref{O-GL-embed}. We obtain
 \begin{equation}\label{M-decomposition}
M \simeq \SO_0(W_{\R}) \times M_{P'}.
\end{equation}
We also define
\begin{equation}\label{pM-decomposition}
\mathfrak{p}_M = \mathfrak{p} \cap \mathfrak{m} = \mathfrak{p}_W \oplus \mathfrak{p}_E,
\end{equation}
where $\frakp_E = \fraksl(E) \cap \frakp$ and
\begin{equation}\label{W-forms}
\mathfrak{p}_W = \mathfrak{o}_W \cap \mathfrak{p}= \Span\{
X_{\alpha\mu} = e_{\alpha} \wedge e_{\mu}; \, \ell+1 \leq \alpha
\leq p, \, p+1 \leq \mu \leq m-\ell\}.
\end{equation}

\subsubsection{Roots}

We let $\underline{S}$ be the maximal $\Q$-split torus of
$\underline{G}$ given by
\begin{equation}
\underline{S} = \left\{ a(t_1,\dots,t_r):= \left(
\begin{smallmatrix}
\diag(t_1,\dots,t_r) & & \\
& 1 &  \\
& & \diag(t_r^{-1},\dots, t_1^{-1}) \end{smallmatrix}
\right)\right\}.
\end{equation}
Note $(\underline{S})(\R))_0=A_{P_0}$, where $P_0$ is the minimal parabolic contained in all standard parabolics. 
We write $\mathbf{t} = (t_1,\dots,t_r)$ and $\widetilde{\mathbf{t}}
= \mathbf{t} J = (t_r,\dots,t_1)$. Note $ a(0,\dots,0,1,0,\dots,0) = \exp (u'_i \wedge u_i) $. The set of simple rational roots for $\underline{G}$ with respect to
$\underline{S}$ is given by
$\Delta=\Delta(\underline{S}, \underline{G}) = \{
\alpha_1,\dots,\alpha_r\}$, where
\begin{align}
\alpha_i(a) &= t_i t_{i+1}^{-1}, \qquad \qquad \qquad \qquad (1 \leq i \leq r-1) \\
\alpha_r(a) &= \begin{cases}  t_r & \text{if $W_r \ne 0$}
\\ t_{r-1} t_r & \text{if $W_r =0$}.
\end{cases}
\end{align}
We write $\Phi({P},A_{\underline{P}})$ for the positive roots of $P$ with
respect to $A_{\underline{P}}$ and $\Delta(P,A_{\underline{P}})$ for
the simple roots of P with respect to $A_{\underline{P}}$, which are
those $\alpha \in \Delta$ which act nontrivially on
$\underline{S}_{\underline{P}}$. We let $\underline{Q}$ be the standard maximal parabolic stabilizing the totally isotropic rational subspace $E_{\ell}$ of
dimension $\ell \leq r$. We have
$A_{\underline{Q}}= \{ a(t,\dots,t,1,\dots1) \}$ and $\Delta(Q,A_{\underline{Q}}) = \{\alpha_{\ell}\}$ unless in the $\Q$-split case for $\SO(p,p)$ and $Q$ stabilizes $E_-$ in which case 
$A_{\underline{Q}}= \{ a(t,\dots,t,t^{-1}) \}$ and $\Delta(Q,A_{\underline{Q}}) = \{\alpha_{p-1}\}$.
For general $P$, we have
\begin{equation}
\Delta(P,A_{\underline{P}}) = \{ \alpha_{i_1},\dots,\alpha_{i_k} \};
\end{equation}
the reader will make the necessary adjustments in the $\Q$-split case for $\SO(p,p)$.

\subsubsection{The nilradical} 

With $P$ and $P'$ as before, we can naturally view $\underline{N}_{\underline{P'}} \subset \SL(E)$ as
a subgroup of $\underline{N}_{\underline{P}}$.
We then have a semidirect product decomposition
\begin{equation}
\underline{N}_{\underline{P}} = \underline{N}_{\underline{P'}}
\ltimes \underline{N}_{\underline{Q}},
\end{equation}
where $\underline{Q}$ is as above the maximal parabolic containing $\underline{P}$. Furthermore, we let $\underline{Z}_{\underline{Q}}$ be the center of
$\underline{N}_{\underline{Q}} \subseteq
\underline{N}_{\underline{P}}$. It is given by
\begin{equation}\label{N-center}
\underline{Z}_{\underline{Q}} =\left\{ z(b):= \left( \begin{smallmatrix}
1    &           &  b     \\
 &    1      &           \\
 &                 &            1
\end{smallmatrix} \right); \, J{^t}bJ = -b \right\}.
\end{equation}
Then for the coset space  $\underline{N}_{\underline{P}}/
(\underline{N}_{\underline{P'}} \ltimes
\underline{Z}_{\underline{Q}})$, we have
\begin{equation}\label{N-embed1}
\underline{N}_{\underline{P}}/ (\underline{N}_{\underline{P'}}
\ltimes \underline{Z}_{\underline{Q}}) \simeq
\underline{N}_{\underline{Q}}/
 \underline{Z}_{\underline{Q}} \simeq W \otimes E
\end{equation}
as vector spaces. Explicitly, the basis of $E$ gives rise to an isomorphism $W \otimes
E \simeq W^{\ell}$. Then for $(w_1,\dots,w_{\ell}) \in W^{\ell}$,
the corresponding coset is represented by
\begin{equation}\label{N-embed2}
n(w_1,\dots,w_{\ell})  :=
\begin{pmatrix}
I_{\ell}&(\cdot,w_1)     &           &       &-w_1^2 \\
 &    \vdots       &           & \cdot  &      \\
 &(\cdot,w_{\ell})&-w_{\ell}^2 &       &      \\
 &   I_W           & -w_{\ell}  & \dots &-w_1   \\
 &                 &           & I_{\ell} &    \\
\end{pmatrix}.
\end{equation}
Here we write $w_i^2=\tfrac12(w_i,w_i)$ for short. On the Lie algebra level, we let $\mathfrak{z}_Q$ be the center of
$\mathfrak{n}_Q \subseteq \frakn_P$, whence corresponding to
\eqref{N-center}
\begin{equation}
\mathfrak{z}_Q \simeq \wwedge{2} E_{\R}.
\end{equation}
We let $\frakn_{P'}$ be the Lie algebra of $N_{P'}$; thus
$\frakn_{P'} \subset E_{\R}' \wedge E_{\R} = \frakgl(E_{\R})$. Corresponding to
\eqref{N-embed2}, we can realize $W_{\R} \otimes E_{\R}$ as a subspace
of $\mathfrak{n}$. Namely, we obtain an embedding
\begin{equation}\label{WE-embedding}
W_{\R}\otimes E_{\R} \hookrightarrow \mathfrak{n},
\end{equation}
\begin{equation}\label{n-embed10}
w \otimes u \to  w \wedge u =:\mathfrak{n}_u(w),
\end{equation}
and we denote this subspace by $\mathfrak{n}_W$, which we frequently
identify with $W_{\R}\otimes E_{\R}$. Furthermore, this embedding is
$\frako(W_{\R}) \oplus \frakgl(E_{\R})$-equivariant, i.e.,
\begin{equation}
[X,\frakn_u(w)] = \frakn_u(Xw) \qquad \qquad [Y,\frakn_u(w)] =
\frakn_{Yu}(w)
\end{equation}
for $X \in \frako(W_{\R})$ and $Y \in \frakgl(E_{\R})$. We easily see
\begin{equation}
\exp( \mathfrak{n}_{u_i}(w)) = n(0,\dots,w,\dots,0).
\end{equation}
A standard basis of $\frakn_W$ is given by
\begin{equation}
X_{\alpha i} := n_{u_i}(e_{\alpha})=  e_{\alpha} \wedge u_i, \qquad
\qquad X_{\mu i} := n_{u_i}(e_{\mu}) = e_{\mu} \wedge u_i
\end{equation}
with $1 \leq i \leq \ell$, $ \ell+1 \leq  \alpha \leq p$, and $ p+1
\leq \mu \leq m-\ell$. The dual space $\frakn_W^{\ast}$ we can identify with $W_{\R} \otimes E'_{\R}$, and we denote the elements of the corresponding
dual basis by $\nu_{\alpha i} = e_{\alpha} \wedge u_i'$ and $\nu_{\mu i}= - e_{\mu} \wedge u_i'$.

Summarizing, we obtain

\begin{lemma}\label{nP-decomposition}
We have a direct sum decomposition (of vector spaces)
\[
\mathfrak{n}_P = \mathfrak{n}_{P'} \oplus \frakn_W \oplus
\mathfrak{z}_Q.
\]
Furthermore, the adjoint action of $\mathfrak{o}(W_{\R}) \oplus
\frakgl(E_{\R})$ on $\frakn_P$ induces an action on the space
$\mathfrak{n}_P/(\mathfrak{n}_{P'} \oplus \mathfrak{z}_Q) \simeq
\frakn_W$ such that
\[
\frakn_W \simeq W_{\R} \otimes E_{\R}
\]
as $\mathfrak{o}(W_{\R}) \oplus \frakgl(E_{\R})$-representations.
\end{lemma}

\subsection{The Maurer Cartan forms and horospherical coordinates}

The Langlands decomposition of $P$ gives rise to the (rational)
horospherical coordinates on $D$ associated to $P$ by
\begin{align}\label{horo-coord}
\sigma=&\sigma_P: N \times A \times D_{\underline{P}}\longrightarrow
D, \\
&\sigma(n,a,m) = n \, a \,m z_0. \notag
\end{align}
Here $D_{\underline{P}} = M_P / K_P$ is the boundary symmetric space associated to $\underline{P}$ with $K_P = M \cap K$. We note that $D_P$ factors into
a product of symmetric spaces for special linear groups and one
orthogonal factor, the symmetric space $D_W$ associated to
$\SO(W)$. We call $D_W$ the
\emph{orthogonal factor} in the boundary symmetric space $D_P$. We
have
\begin{equation}\label{boundary-factors}
D_P = D_W \times \prod_{j=1}^{k} D_{U_{i_j}},
\end{equation}
where $ D_{U_{i_j}}$ denotes the symmetric space associated to
$\SL(U_{i_j})$.

We now describe the basic cotangent vectors $\omega_{\alpha\mu}
= (e_{\alpha} \wedge e_{\mu})^{\ast} \in \mathfrak{p}^{\ast} \simeq
T_{z_0}^{\ast}(D)$ in $NAM$ coordinates.
We extend $\sigma$ to $ N \times A \times M \times K \longrightarrow
G$ by $\sigma (n,a,m,k) = namk$, and this induces an isomorphism
between the left-invariant forms on $NAM$ (which we identify with
$\mathfrak{n}^{\ast} \oplus \mathfrak{a}^{\ast} \oplus
\mathfrak{p}_M^{\ast}$) and the horizontal left-invariant forms on
$G$ (which we identify with $\mathfrak{p}^{\ast}$). Thus we have an
isomorphism
\begin{equation}\label{pnam-iso}
\sigma^{\ast}: \mathfrak{p}^{\ast} \longrightarrow
\mathfrak{n}^{\ast} \oplus \mathfrak{a}^{\ast} \oplus
\mathfrak{p}_M^{\ast}.
\end{equation}

\begin{lemma}\label{forms-formulas}
Let $1 \leq i \leq \ell $. For the preimage under $\sigma^{\ast}$ of
the elements in $\mathfrak{n}^{\ast}_W$ coming from $W_+ \otimes E$,
we have
\begin{equation}\label{invariantforms1}
\sigma^{\ast} \; \omega_{\alpha m+1-i} = - \tfrac1{\sqrt{2}}\nu_{\alpha i},
\end{equation}
where $\ell+1 \leq \alpha \leq p$. Furthermore, for the ones coming
from $W_- \otimes E$, we have
\begin{equation}
\sigma^{\ast} \omega_{ i \mu}  =  \tfrac1{\sqrt{2}} \nu_{\mu i},
\end{equation}
where $p+1 \leq \mu \leq m+1-\ell$. On $\mathfrak{p}_M^{\ast}$, the
map $\sigma^{\ast}$ is the identity. In particular, for $ \ell+1 \leq \alpha \leq p$ and any $\mu\geq p+1$,
we have
\begin{equation}\label{throw-away-key}
\sigma^{\ast} \omega_{\alpha \mu} \in \mathfrak{p}^{\ast}_W \oplus
\mathfrak{n}^{\ast}_W.
\end{equation}
The remaining elements of $\mathfrak{p}^{\ast}$ are of the form
$\omega_{i \mu}$ with $ p+1 \leq \mu \leq m+1-\ell$. These elements
are mapped under $\sigma^{\ast}$ to $ \frakn_{P'}^{\ast} \oplus
\mathfrak{a}^{\ast} \oplus \mathfrak{p}^{\ast}_E \subset
\mathfrak{gl}(E_{\R})^{\ast}$.
\end{lemma}

\subsection{Borel-Serre Compactification}\label{BS-C1}

We now briefly describe the Borel-Serre compactification of $D$ and of $X=\G
\back D$. For a more detailed discussion see also the last section where we discuss the $\Q$-split case for $\SO(p,p)$ in detail. In that situation the Borel-Serre compactification is {\it not} the right compactification for our purposes, and we need to work with a slightly larger compactification.

We follow \cite{BJ}, III.9. We first partially compactify
the symmetric space $D$. For any rational parabolic $\underline{P}$,
we define the boundary component
\begin{equation}
e(\underline{P}) = N_P \times D_{\underline{P}} \simeq P /
A_{\underline{P}} K_{\underline{P}}.
\end{equation}
Then as a set the (rational) Borel-Serre enlargement
$\overline{D}^{BS} = \overline{D}$ is given by
\begin{equation}
\overline{D} = D \cup \coprod_{\underline{P}} e(\underline{P}),
\end{equation}
where $\underline{P}$ runs over all rational parabolic subgroups of
$\underline{G}$. As for the topology of $\overline{D}$, we first
note that $D$ and $e(\underline{P})$ have the natural topology.
Furthermore, a sequence of $y_j =\sigma_P(n_j,a_j,z_j) \in D$ in
horospherical coordinates of $D$ converges to a point $(n,z) \in
e(\underline{P})$ if and only if $n_j \to n$, $z_j \to z$ and
$\alpha(a_j) \to \infty$ for all roots $\alpha \in
\Phi(\underline{P},A_{\underline{P}})$. For convergence within
boundary components, see \cite{BJ}, III.9.

With this, $\overline{D}$ has a canonical structure of a real
analytic manifold with corners. Moreover, the action of
$\underline{G}(\Q)$ extends smoothly to $\overline{D}$. The action
of $g = kp=kman \in KMAN = G$ on $e(\underline{P})$ is given by
\begin{equation}
g \cdot(n',z') = k \cdot( Ad(am)(nn'),mz') \in e(Ad(k)\underline{P})
=e(Ad(g)\underline{P})
\end{equation}
with $k\cdot(n',z') = (Ad(k)n, Ad(k)m K_{Ad(k)\underline{P}}) \in
e(Ad(k)\underline{P})$. Finally,
\begin{equation}
\overline{X}:= \G \back \overline{D}
\end{equation}
is the Borel-Serre compactification of $X = \G \back D$ to a
manifold with corners. If $\underline{P}_1, \dots,\underline{P}_k$
is a set of representatives of $\G$-conjugacy classes of rational
parabolic subgroups of $\underline{G}$, then
\begin{equation}
\G \back \overline{D} = \G \back D \cup \coprod_{i=1}^k
\G_{{P}_i} \back e(\underline{P})_i,
\end{equation}
with $\G_{{P}_i} = \G \cap P_i$. We will write
$e'(\underline{P}) = \G_{\underline{P}} \back e(\underline{P})$. We write $\G_M$ for the image of $\G_P$ under the quotient map $P \to P/N$. Furthermore, $\G_P$ acts on $E^{\perp}_{\R}/E_{\R}$, and we denote this transformation group by $\G_W$. Note that $\G_M$ and $\G_W$ when viewed as subgroups of $P$ contain $\G \cap M$ and $\G \cap \SO_0(W_{\R})$ respectively as subgroups of finite index.

\bigskip

We now describe Siegel sets. For $t \in \R_+$, let
\begin{equation}
A_{P,t} = \{ a \in A_P; \, \alpha(a) >t \text{ for all } \alpha \in
\Delta(P,A_{\underline{P}})\},
\end{equation}
and for bounded sets $U \subset N_P$ and $V \subset D_P$, we define
the Siegel set
\begin{equation}
\mathfrak{S}_{P,U,t,V} = U \times A_{P,t} \times V \subset N_P
\times A_P \times D_P.
\end{equation}
Note that for $t$ sufficiently large, two Siegel sets for different
parabolic subgroups are disjoint. Furthermore, if $P_1, \dots, P_k$
are representatives of the $\underline{G}(\Q)$-conjugacy classes of
rational parabolic subgroups of $G$, then there are Siegel sets
$\calS_i$ associated to $P_i$ such that the union $\bigcup \pi
(\calS_i)$ is a fundamental set for $\G$. Here $\pi$ denotes the
projection $\pi: D \to \G \back D$.

\section{Review of representation theory for general linear and
orthogonal groups}\label{rep-review}

In this section, we will briefly review the construction of the
irreducible finite dimensional (polynomial) representations of
$\GL(\C^n)$ and $\Orth(V)$. Here, in this section, we assume that
$V$ is an orthogonal complex space of dimension $m$. Basic references are
\cite{FultonHarris}, \S 4.2 and \S 6.1, \cite{GoodmanWallach}, \S
9.3.1-9.3.4 and \cite{Boerner}, Ch.~V, \S 5 to which we refer for
details.

\subsection{Representations of $\GL_n(\C)$}

Let $\la = (b_1,\dots, b_n)$ be a partition of $\ell'$ with the $b_i$'s arranged in decreasing order. We will
use $D(\la)$ to denote the associated Young diagram. We identify the partition $\la$ with the dominant weight $\la$ for
$GL(n)$ in the usual way. A standard filling $\la$ of the Young
diagram $D(\la)$ by the elements of the set $[\ell']=\{1,2, \cdots,
\ell' \}$ is an assignment of each of the numbers in $[\ell']$ to a
box of $D(\la)$ so that the entries in each row strictly increase
when read from left to right and the entries in each column strictly
increase when read from top to bottom. A Young diagram equipped with
a standard filling will be also called a standard tableau.

We let $s_{t(\la)}$ be the idempotent in the group algebra of the
symmetric group $S_{\ell'}$ associated to a standard tableau
with $\ell'$ boxes corresponding to a standard filling $t(\lambda)$
of a Young diagram $D(\lambda)$. Note that $S_{\ell'}$ acts on the
space of $\ell'$-tensors $T^{\ell'}(\C^n)$ in the natural fashion on
the tensor factors. Therefore $s(t(\la))$ gives rise
to a projection operator in $\End(T^{\ell'}(\C^n))$, which we also denote by $s_{t(\la)}$. We write
\begin{equation}
\mathbb{S}_{t(\la)}(\C^n) = s_{t(\la)}(T^{\ell'}(\C^n)).
\end{equation}
We have a direct sum decomposition
\begin{equation}
T^{\ell'} (\C^n) = \bigoplus_{\la} \bigoplus_{t(\la)}
\mathbb{S}_{t(\la)}(\C^n),
\end{equation}
where $\la$ runs over all partitions of $\ell'$ and $t(\la)$ over
all standard fillings of $D(\la)$. This gives the decomposition of
$T^{\ell'}(\C^n)$ into irreducible constituents, i.e, for every
standard filling $t(\lambda)$, the $\GL(\C^n)$-module
$\mathbb{S}_{t(\la)}(\C^n)$ is irreducible with highest weight
$\lambda$. In particular, $\mathbb{S}_{t(\la)}(\C^n)$ and
$\mathbb{S}_{t'(\la)}(\C^n)$ are isomorphic for two different
standard fillings $t(\la)$ and $t'(\la)$. We denote this isomorphism
class by $\mathbb{S}_{\la}(\C^n)$ (or if we do not want to specify
the standard filling).

Explicitly, we let $A$ be the standard filling of a Young diagram
$D(A)$ corresponding to the partition $\lambda$ with less than or
equal to $n$ rows and $\ell'$ boxes by $1,2,\cdots,\ell'$ obtained
by filling the rows in order beginning at the top with
$1,2,\cdots,\ell'$.  We let $R(A)$ be the subgroup of $S_{\ell'}$
which preserves the rows of $A$ and $C(A)$ be the subgroup that
preserves the columns of $A$. We define elements $r_A$ and $c_A$
by
\begin{equation}
r_A = \sum_{s \in R(A)}s \quad \text{and} \quad c_A = \sum_{s \in
C(A)} \sgn(s) s.
\end{equation}
Let $h(A)$ be the product of the hook lengths of the boxes in
$D(A)$, see \cite{FultonHarris}, page 50. Then the idempotent $s_A$
is given
\begin{equation}
s_A =\frac{1}{h(A)}c_Ar_A.
\end{equation}
We will also need the "dual" idempotent $s_A^{\ast}$ given by $s_A^{\ast} = \tfrac{1}{h(A)}r_Ac_A$.
We let $\eps_1,\dots,\eps_n$ denote the standard basis of $\C^n$ and
$\theta_1,\dots,\theta_n \in (\C^n)^{\ast}$ be its dual basis. We
set
\begin{equation}
\eps_A = \eps_1^{b_1} \otimes \cdots \otimes \eps_n^{b_n}
\end{equation}
and let $\theta_A$ be the corresponding element in
$T^{\ell'}(\C^n)^{\ast}$. Then $s_A(\eps_{A})$ is a highest weight
vector in $\mathbb{S}_{A}(\C^n)$, see \cite{GoodmanWallach},
\S9.3.1. We have

\begin{lemma}
Let $|R(A)|$ be the order of $R(A)$. Then 
\[
s_A^*\theta_{A}(s_A\eps_{A}) = \frac{|R(A)|}{h(A)}.
\]
\end{lemma}
\begin{proof}
We compute
\begin{align*}
s_A^*\theta_{A}(s_A\eps_{A}) = \theta_{A}(s_A^2 \eps_A)
=\theta_A(s_A \eps_A)  = \frac{|R(A)|}{h(A)}\theta_A(c_A \eps_A) =
\frac{|R(A)|}{h(A)} \theta_A( \eps_A).
\end{align*}
The last equation holds because $\theta_A(q\eps_A)=0$ for any
nontrivial $q$ in the column group of $A$ as the reader will easily
verify. We have used $r_A \eps_A = |R(A)|\eps_A$ (since all elements of $R(A)$ fix $\eps_A)$ and $s_A= \frac{1}{h(A)}c_Ar_A$.
\end{proof}

\subsection{Enlarging the Young diagram}\label{box-section}
We let $B=B_{n,\ell}$ be the standard tableau with underlying shape
$D(B)$ an $n$ by $\ell$ rectangle with the standard filling obtained
by putting $1$  through $\ell$ in the first row, $\ell +1$ through
$2\ell$ in the second row etc. Then $D(B)$ is the Young diagram
corresponding to the dominant weight $\ell \varpi_n$. Here $\varpi_n
= (1,1,\cdots,1)$ is the $n$-th fundamental weight for $GL(n)$. We
note that we have $\eps_B = \epsilon_1^{\ell} \otimes \cdots \otimes
\epsilon_n^{\ell}$ and $\theta_{B}= \theta_1^{\ell} \otimes \cdots
\otimes \theta_n^{\ell}$. 

\begin{lemma}\label{basisforonedimensionalspace}
The space $s_BT^{n\ell}(\C^n)$ is $1$-dimensional, and 
\[
s_BT^{n\ell}(\C^n) = \C s_B\eps_{B}
\]
as $\GL(n,\C)$-modules. Correspondingly, $s_B^*T^{n\ell}(\C^n)^*$ is $1$-dimensional and
\[
s_B^*T^{n\ell}(\C^n)^* = \C s_B^*\theta_{B}.
\]
In particular, 
\[
s_B^*T^{n\ell}(\C^n)^* \cong \left(
\wwedge{n}(\C^n)^{\ast}\right)^{\otimes \ell}.
\]
\end{lemma}

We let $A$ be the standard filling of the Young diagram $D(\la)$ as
above. Then $B|A$ denotes the standard tableau with underlying shape
$D(B|A)$ given by making the shape of $A$ abut $B$ (on the right),
using the above filling for $B$ and filling $A$ in the standard way
(as above) with $n\ell + 1$ through $n\ell + \ell'$. For example, if
\begin{center}
$B=$
\begin{tabular}{| c | c | c |}
\hline 1  & 2  & 3  \\ \hline 4  & 5  & 6  \\ \hline 7  &  8 & 9\\
\cline{1-3}
\end{tabular}
\quad and \quad $A=$
\begin{tabular}{| c | c | c |}
\hline
    1  &     2  &     3  \\ \hline
    4  &     5 \\ \cline{1-2}
\end{tabular},
\quad then \quad $B|A = $
\begin{tabular}{| c | c | c |c|c|c|}
\hline 1  & 2  & 3 & 10 & 11 & 12  \\ \hline 4  & 5  & 6 & 13 & 14
\\ \cline{1-5} 7  &  8 & 9\\ \cline{1-3}
\end{tabular}
\end{center}
We have an idempotent $s_{B|A}$ in the group ring of $S_{n\ell
+\ell'}$ and $\eps_{B|A} \in
T^{n\ell+\ell'}(\C^n)$, which give rise to a highest weight vector
$s_{B|A}\eps_{B|A}$ in $s_{B|A}(T^{n\ell+\ell'}(\C^n))$. Note
\begin{equation}\label{P1}
\eps_{B|A} = \eps_B \otimes \eps_A.
\end{equation}

\begin{lemma}\label{John}
There is a positive number $c(A,B)$ such that
$$s_B\eps_B \otimes s_A \eps_A = c(A,B) s_{B|A}\eps_{B|A}.$$
\end{lemma}
\begin{proof}
Since the Young diagrams $D(B)$ and $D(A)$ are abutted along their
vertical borders, we see
\begin{equation}\label{P2} c_{B|A} =(c_B\otimes
1_{\ell'})\circ (1_{n\ell} \otimes c_A) = (1_{n\ell} \otimes
c_A)\circ (c_B\otimes 1_{\ell'}).
\end{equation}
Also $r(C) \eps_C = |R(C)|\eps_C$. Then we easily compute (using \eqref{P1} and \eqref{P2})
\begin{align*}
&s_B\eps_B \otimes s_A \eps_A  
=
\frac{h(B|A)}{h(B)h(A)}\frac{|R(B)||R(A)|}{|R(B|A)|}s_{B|A}\eps_{B|A}. \qedhere
\end{align*}
\end{proof}

\begin{corollary}\label{repsA-Bcoincide}
Under the identification of $T^{n\ell}(\C^n) \otimes T^{\ell'}(\C^n)
\to T^{n\ell+\ell'}(\C^n)$ given by tensor multiplication, we have
the equality of maps
\[
s_B \otimes s_A = s_{B|A}.
\]
That is,
\[
\mathbb{S}_B(\C^n) \otimes \mathbb{S}_A(\C^n) =
\mathbb{S}_{B|A}(\C^n)
\]
as (physical) subspaces of $T^{n\ell+\ell'}(\C^n)$. The same
statements hold for the dual space
$\mathbb{S}_{B|A}^{\ast}(\C^{n\ell+\ell'})^{\ast}$ etc.
\end{corollary}

\begin{proof}
Since $\mathbb{S}_B(\C^n)$ is one-dimensional, the tensor product
$\mathbb{S}_B(\C^n) \otimes \mathbb{S}_A(\C^n)$ defines an
irreducible representation for $\GL_n(\C^n)$ (under the tensor
multiplication map $T^{n\ell}(\C^n) \otimes T^{\ell'}(\C^n)$ inside
$T^{n\ell+\ell'}(\C^n)$). But by Lemma~\ref{John} it has nonzero
intersection with the irreducible $\GL_n(\C)$-representation
$\mathbb{S}_{B|A}(\C^n)$ inside $T^{n\ell+\ell'}(\C^n)$. Hence the
two subspaces coincide.
\end{proof}

\subsection{Representations of $\Orth(V)$}

We extend the bilinear form $(\, ,\, )$ on $V$ to $T^{\ell'}(V)$ as
the $\ell'$-fold tensor product and note that the action of
$S_{\ell'}$ on $T^{\ell'}(V)$ is by isometries. We let $V^{[\ell']}$ be the space of harmonic $\ell'$-tensors (which
are those $\ell'$-tensors which are annihilated by all contractions
with the form $(\,,\,)$). We let $\mathcal{H}$ be the orthogonal
projection $\mathcal{H}:T^{\ell'}(V) \to V^{[\ell']}$ onto the
harmonic $\ell'$-tensors of $V$. Note that $V^{[\ell']}$ is invariant under the action of $S_{\ell'}$. We then
define for $\la$ as above the harmonic Schur functor $\mathbb{S}_{[t(\lambda)]}(V)$ by
\begin{equation}\label{harmonic}
\mathbb{S}_{[t(\lambda)]}(V) = \mathcal{H}
\mathbb{S}_{t(\lambda)}(V).
\end{equation}
If the sum of the lengths of the first two columns of $D(\la)$ is at most $m$, then $\mathbb{S}_{[t(\la)]}(V_{\C})$ is a nonzero irreducible representation for $\Orth(V_{\C})$, see \cite{FultonHarris} section~19.5. Otherwise, it vanishes. Of course, for different fillings $t(\la)$ of $D(\la)$, these representations are all isomorphic and we write $\mathbb{S}_{[\lambda]}(V)$ for the isomorphism class. Furthermore, it is also irreducible when restricted to $G$ unless $m$ is even and $i(\la)=\tfrac{m}2$, in which case it splits into two irreducible representations. If $i(\la) \leq [\tfrac{m}2]$, then the corresponding highest weight $\tilde{\la}$ for the representation $\mathbb{S}_{[\lambda]}(V)$ of $G$ has the same nonzero entries as $\la$.

\section{The Weil representation}

We review different models of the Weil representation. In this
section, $V$ denotes a real quadratic space of signature $(p,q)$ and
dimension $m$.

We let $V'$ be a real symplectic space of dimension $2n$. We denote
by $G'= \Mp(n,\R)$ the metaplectic cover of the symplectic group
$\Sp(V') = \Sp(n,\R)$ and let $\mathfrak{g}'$ be its Lie algebra. We
let $K'$ be the inverse image of the standard maximal compact
$\Uni(n) \subset \Sp(n,\R)$ under the covering map $\Mp (n,\R)
\rightarrow \Sp(n,\R)$. Note that $K'$ admits a character
$\det^{1/2}$, i.e., its square descends to the determinant character
of $\Uni(n)$. The embedding of $\Uni(n)$ into $\Sp(n,\R)$ is given
by $A + iB \mapsto \left(
\begin{smallmatrix}
A & B \\
-B &A
\end{smallmatrix}\right)$.
We write $\mathcal{W}_{n,V}$ for (an abstract model of) the
$K'$-finite vectors of the restriction of the Weil representation of
$\Mp(V' \otimes V)$ to $\Mp(n,\R) \times \Orth(V)$ associated to the
additive character $t \mapsto e^{2\pi i t}$.

\subsection{The Schr\"odinger model}

We let $V_1'$ be a Langrangian subspace of $V'$. Then $V \otimes
V_1'$ is a Langrangian subspace of $V \otimes V'$ (which is
naturally a symplectic space of dimension $2nm$). The Schr\"odinger
model of the Weil representation consists of the space of
(complex-valued) Schwartz functions on the Lagrangian subspace
$V_1' \otimes V\simeq V^n$. We write $\calS(V^n)$ for the space of
Schwartz functions on $V^n$ and write $\omega=\omega_{n,V}$ for the
action.

The Siegel parabolic $P'=M'N'$ has Levi factor
\begin{equation}
M'=\left\{ m'(a) =
\begin{pmatrix} a&0\\0&^ta^{-1}
    \end{pmatrix}; a \in \GL(n,\R)
\right\}
\end{equation}
and unipotent radical
\begin{equation}
N'=\left\{ n'(b) = \begin{pmatrix} 1&b\\0&1
    \end{pmatrix}; b \in \Sym_n(\R) \right\}.
\end{equation}
It is well known that we can embed $P'$ into $\Mp(n,\R)$, and the action of $P'$ on $S(V^n)$ is given by
\begin{align}
\omega \left( m'(a) \right) \varphi(\x) &=
 (\det a)^{m/2}\varphi(\x a) \qquad \qquad (\det a >0),\\
\omega \left( n'(b) \right) \varphi(\x) &= e^{\pi i tr(b(\x,\x))}
\varphi(\x)
\end{align}
with $\x =(x_1, \dots, x_n) \in V^n$. 
The orthogonal group $G$ acts on $\calS(V^n)$
via
\begin{equation}
\omega(g) \varphi(\x) = \varphi(g^{-1}\x),
\end{equation}
which commutes with the action $G'$. The standard Gaussian is given
by
\begin{equation}\label{Gauss}
\varphi_0(\x) = e^{-\pi tr(\x,\x)_{z_0}} \in \mathcal{S}(V^n)^K.
\end{equation}
 Here $(\x,\x)$ is the inner product matrix $(x_i,x_j)_{ij}$. 

We let $S(V^n)$ be the space of $K'$-finite vectors
inside the space of Schwartz functions on $V^n$. It consists of
those Schwartz functions of the form $p(\x)\varphi_0(\x)$, where $p$
is a polynomial function on $V^n$.

\subsection{The mixed model and local restriction for the Weil representation}\label{mixed}

We let $P$ be a standard parabolic of $G$ stabilizing a totally isotropic flag in $V$ with $E=E_{\ell}$ be the largest constituent of the flag and associated Witt decomposition $V = E \oplus W \oplus E'$. 

We describe a different model for the Weil representation, the
so-called mixed model. Furthermore, we will define a "local"
restriction $r_P^{\mathcal{W}}$ from $\mathcal{S}(V^n)$ to the space
of Schwartz functions $\mathcal{S}(W^n)$ for $W$, a subspace
of signature $(p-\ell,q-\ell)$.

\subsubsection{The mixed model}

We let $E=E_{\ell}$ be one of the standard totally isotropic
subspaces of $V $, see \eqref{isospace}. As before, we identify the
dual space of $E$ with $E'$. Accordingly, we write $\x = \left( \begin{smallmatrix} u \\\x_W \\u' \end{smallmatrix} \right)$ for $\x \in V^n$, where  $u \in E^n$, $u' \in (E')^n$, and $\x_W \in W^n$. We then have
an isomorphism of two models of the Weil representation given by
\begin{align}
\calS(V^n) &\longrightarrow \calS((E')^n) {\otimes} S(W^n) {\otimes}
\calS((E')^n) \\
\varphi &\longmapsto \widehat{\varphi} \notag
\end{align}
given by the partial Fourier transform operator
\begin{equation}
\widehat{\varphi} \left(\begin{smallmatrix} \xi \\\x_W \\u' \end{smallmatrix} \right) =
\int_{E^n} \varphi \left(\begin{smallmatrix} u \\\x_W \\u' \end{smallmatrix} \right)
e^{-2\pi i tr(u,\xi)} du
\end{equation}
with $\xi,u' \in (E')^n$ and $\x_W \in W^n$. We need some formulae relating the action of $\omega$ in the two models. 

\begin{lemma}\label{mixedformulas}
Let $\left( \begin{smallmatrix} \xi \\\x_W
\\ u'   \end{smallmatrix} \right) \in (E' \oplus W\oplus E')^n$.
\begin{itemize}
\item[(i)] Let $n \in N_Q$ and write $n(u')_W$ for the image of $n(u')$ under the orthogonal projection onto $W$. Then
      \begin{equation*}
      \widehat{n\varphi} ({^t(\xi,\x_W,u')}) =
              e \left(\tr(n(\x_W+u'),\xi)\right)
             \widehat{\varphi} ({^t(\xi,\x_W + n(u')_W,u')}).
      \end{equation*}

\item[(ii)] For $g \in \SL(E)\subset G$ (in particular, $g \in N_{P'}$ or $g \in M_{P'}$) we have
\[
\widehat{g\varphi}({^t(\xi,\x_W,u')}) =
\widehat{\varphi}({^t(\widetilde{g}\xi,\x_W,\widetilde{g}^{-1}u')})
\]
with $\widetilde{g} = Jg^{\ast}J$ and $g^{\ast} = {^t g^{-1}}$.

\item[(iii)]For $\mathbf{t} =(t_1,\dots,t_{\ell})$, set
$\widetilde{\mathbf{t}} = \mathbf{t} J = (t_{\ell},\dots,t_1)$ and
$|\mathbf{t}| = t_1 \cdot t_2 \cdots t_{\ell}$. Then
      \begin{equation*}
\widehat{a(\mathbf{t})\varphi} ({^t(\xi,\x_W,u')}) = |\mathbf{t}|^n
\widehat{\varphi}({^t(\widetilde{\mathbf{t}}\xi,\x_W,\widetilde{\mathbf{t}}u')}).
      \end{equation*}

\item[(iv)]
For $h \in \SO_0(W) \subset M$, we have
\[
\widehat{h\varphi}({^t(\xi,\x_W,u')})  =
\widehat{\varphi}({^t(\xi,h^{-1}\x_W,u')}).
\]

\item[(v)] For $m'(a) = \left( \begin{smallmatrix} a&0\\0&^ta^{-1} \end{smallmatrix}  \right)
        \in M' \subset Sp(n,\R)$ with $a \in GL^+_n(\R)$,
      \[
       \widehat{(m'(a)\varphi)}({^t(\xi,\x_W,u')})
       = (\det a)^{\frac{m}2 -\ell} \widehat{\varphi}({^t(\xi a^{\ast},\x_Wa,u'a)})
      \]

\item[(vi)] For $n'(b) = \left( \begin{smallmatrix} 1&b\\0&1 \end{smallmatrix} \right)
 \in N' \subset Sp(n,\R)$ with $b \in \Sym_n(\R)$,
       \[
       \widehat{(n'(b)\varphi)}({^t(\xi,\x_W,u')}) =
       e \left(tr(b\tfrac{(\x_W,\x_W)}{2}) \right) \widehat{\varphi}({^t(\xi-u'b,\x_W,u')}).
       \]
\end{itemize}
\end{lemma}

\begin{proof}
This is straightforward.
\end{proof}

We obtain

\begin{proposition}\label{intertwiner}
Let $\varphi \in \calS(V^n)$. Then the restriction of
$\widehat{\varphi}$ to $W^n$,
\[
\varphi \mapsto \widehat{\varphi}|_{W^n},
\]
defines a $G' \times MN$ intertwiner from $\calS(V^n)$
to $\calS(W^n)$. Here, we identify $W$ with $E^{\perp}/E$ to define
the action of $MN$ on $W$. In particular, $N$ and $M_{P'}$ (see
\ref{M-decomposition}) act trivially on $\mathcal{S}(W^n)$.

\end{proposition}

\subsubsection{Weil representation restriction}

\begin{definition}\label{def:localresWeil}
Let $\varphi \in \calS(V^n)$ and let $P$ be the parabolic as before. We {define} the "local" restriction
$r^{\mathcal{W}}_P (\varphi) \in \calS(W^n)$ with respect to $P$ for
the Schr\"odinger model of the Weil representation $\mathcal{W}$ by
\[
r^{\mathcal{W}}_P (\varphi) = \widehat{\varphi}|_{W^n}.
\]
\end{definition}

We now describe this restriction on a certain class of Schwartz functions on $V^n$. For $\x =(x_1, \dots, x_n) \in V^n$, we write $\left(
\begin{smallmatrix} x_{1j}  \\ \vdots \\ x_{mj} \end{smallmatrix}
\right)$ for the standard coordinates of $x_j$. We define a family
of commuting differential operators on $\calS(V^n)$ by
\begin{equation}
\calH_{r j}  = \left(x_{rj} - \frac{1}{2\pi}
\frac{\partial}{\partial x_{r j}} \right),
\end{equation}
where $1 \leq r \leq m$ and $1 \leq j \leq n$. Define a polynomial $\widetilde{H}_k$ by 
\begin{equation}
\widetilde{H}_k(x)=(2\pi)^{-k/2} H_k\left(\sqrt{2\pi}x\right),
\end{equation}
where $H_k(x) =(-1)^k e^{x^2} \tfrac{d^k}{dx^k} e^{-x^2}$ is the $k$-th Hermite polynomial. Then \begin{equation}
\calH_{r j}^{k} \varphi_0(\x) = \widetilde{H}_k(x_{rj})
\varphi_0(\x),
\end{equation}
where $\varphi_0(\x)$ is the standard Gaussian, see \eqref{Gauss}. We let $\Delta \in M_{m \times n}(\Z) =(\delta_{rj})$ be an integral
matrix with non-negative entries and split $\Delta$ into $\Delta_+
\in M_{p\times n}(\Z)$ and $\Delta_- \in M_{q \times n}(\Z)$ into
its "positive" and "negative" part, where $\Delta_+$ consists of the
first $p$ rows of $\Delta$ and $\Delta_-$ of the last $q$. (Recall
$m = p+q$). We define operators
\begin{align*}
\mathcal{H}_{\Delta} = \prod_{\substack{ 1 \leq r \leq m \\ 1 \leq j
\leq n}} \mathcal{H}^{\delta_{r j}}_{rj}, \qquad
 \mathcal{H}_{\Delta_+} = \prod_{\substack{ 1 \leq \alpha \leq p \\ 1
\leq j \leq n}} \mathcal{H}^{\delta_{\alpha j}}_{\alpha j},  \qquad
\mathcal{H}_{\Delta_-} = \prod_{\substack{ p+1 \leq \mu \leq m \\ 1
\leq j \leq n}} \mathcal{H}^{\delta_{\mu j}}_{\mu j},
\end{align*}
so that $\mathcal{H}_{\Delta} = \mathcal{H}_{\Delta_+}
\mathcal{H}_{\Delta_-}$. Here again we make use of our convention to
use early Greek letters for the "positive" indices of $V$ and late
ones for the "negative" indices.
\begin{definition}\label{varphidelta}
For $\Delta$ as above, we define the Schwartz function
$\varphi_{\Delta}$ by
\[
\varphi_{\Delta}(\x) = \mathcal{H}_{\Delta} \varphi_0(\x) =
\prod_{\substack{ 1 \leq \alpha \leq p \\  p+1 \leq \mu \leq m  \\1
\leq j \leq n}} \widetilde{H}_{\delta_{\alpha j}}(x_{\alpha j})
\widetilde{H}_{\delta_{\mu j}}(x_{\mu j}) \varphi_0(\x).
\]
\end{definition}

We now describe $\varphi^V_{\Delta}$ in the mixed model. The
superscript $V$ emphasizes that the Schwartz function is associated
to the space $V$. We begin with some auxiliary considerations. The
following little fact will be crucial for us.

\begin{lemma}\label{keylemma}
For a Schwartz function $f \in \calS(\R)$, let $\widehat{f}(\xi) =
\int_{\R} f(y) e^{-2\pi i y\xi} dy$ be its Fourier transform. Let
$g_k(y) =\widetilde{H}_k(-\tfrac{y}{\sqrt{2}}) e^{- \pi y^2}$. Then
\[
\widehat{g_k}(\xi) = (-\sqrt{2}i \xi)^k e^{-\pi \xi^2}.
\]
\end{lemma}

\begin{proof}
We use induction and differentiate the equation $
\widehat{(\widehat{g_k})}(-y) =\widetilde{H}_k(\tfrac{y}{\sqrt{2}}) e^{- \pi
y^2}$. The assertion follows from the recursion
$\widetilde{H}_{k+1}(y) = 2y\widetilde{H}_k(y)- \tfrac{1}{2\pi}
\widetilde{H}'_k(y)$, which is immediate from the definition of
$\widetilde{H}_k$. The claim also follows easily from \cite{Leb},
(4.11.4).
\end{proof}

\begin{remark}
Recall that on the other hand $\widetilde{H}_k(y) e^{- \pi
y^2}$ is an eigenfunction under the Fourier transform with
eigenvalue $(-i)^k$, see \cite{Leb}, (4.12.3). This fact is
underlying the automorphic properties of the theta series associated
to the special forms $\varphi_{nq,[\la]}$.
\end{remark}

The Gaussian is given in standard coordinates by $ \varphi^V_0(\x) =
\exp ( -\pi \sum_{j=1}^n \sum_{i=1}^m x^2_{ij} )$. In Witt coordinates, we have $x_{rj}= \tfrac1{\sqrt{2}}(y_{rj} -
y_{(m-r)j})$ and $x_{(m-r)j} = \tfrac1{\sqrt{2}}(y_{rj} + y_{(m-r)j})$; thus
$x_{rj}^2 + x^2_{(m-r),j} = y_{rj}^2 + y_{(m-r)j}^2$ for
$r \leq \ell$. Thus
\begin{equation}
\varphi_0^V \left(\begin{smallmatrix}  u \\ \x_W \\ u'
\end{smallmatrix} \right)  = \exp\left(-\pi \sum_{j=1}^n
\sum_{r=1}^{\ell}(y_{rj}^2 + y_{(m-r)j}^2) \right)
\varphi_0^W(\x_W).
\end{equation}
We write slightly abusing
\begin{equation}
\varphi_0^E(u,u') := \varphi_0^V\left(\begin{smallmatrix}  u \\ 0\\
u' \end{smallmatrix} \right) =\exp\left(-\pi \sum_{j=1}^n
\sum_{r=1}^{\ell}(y_{rj}^2 + y_{(m-r)j}^2) \right).
\end{equation}

We let $\Delta'$ be the truncated matrix of size $(m-2\ell)\times n$
given by eliminating the first and the last $\ell$ rows from
$\Delta$. We let $\Delta''$ be the matrix of these eliminated rows.
Note that $\mathcal{H}_{\Delta'}$ now defines an operator on
$\mathcal{S}(W^n)$ and $\mathcal{H}_{\Delta''}$ on
$\mathcal{S}((E\oplus E')^n)$. We also obtain matrices $\Delta'_+$
of size $(p-\ell)\times n$ and $\Delta'_-$ of size $(q-\ell) \times
n$ by eliminating the first $\ell$ and the last $\ell$ rows from
$\Delta_+$ and $\Delta_-$ respectively. Similarly we obtain
$\Delta''_+$ and $\Delta''_-$. With these notations we obtain

\begin{lemma}\label{mixedformula10}

\begin{equation}
\widehat{\varphi^V_{\Delta}}
\left(\begin{smallmatrix}  \xi \\ \x_W \\ u'
\end{smallmatrix} \right)
= \varphi^W_{\Delta'}(\x_W)
 \widehat{\varphi^E_{\Delta^{\prime \prime}}}
(\xi,u'). \tag{i}
\end{equation}

\begin{equation}
r_P^{\mathcal{W}}\left( \varphi^V_{\Delta}\right)(\x_W) =
\varphi^W_{\Delta'}(\x_W) \widehat{\varphi^E_{\Delta^{\prime
\prime}}} (0,0). \tag{ii}
\end{equation}

\end{lemma}

In our applications all entries of $\Delta_-$ will be zero, so
$\Delta =\Delta_+$ (by abuse of notation). 

\begin{lemma}\label{keylemma2}
\[
\widehat{\varphi^E_{\Delta_+^{\prime\prime}}}(\xi,0) =
\left( \prod_{j=1}^n \prod_{\alpha=1}^{\ell} (-\sqrt{2}i\xi_{\alpha
j})^{\delta_{\alpha j}} \right) \varphi_0^E(\xi,0).
\]
In particular, if in addition \emph{all} entries of $\Delta_+''$
vanish, then
\[
\widehat{\varphi^V_{\Delta_+}}
\left(\begin{smallmatrix}  \xi \\ \x_W \\ 0
\end{smallmatrix} \right)=
\varphi^W_{\Delta_+'}(\x_W) \varphi_0^E(\xi,0).
\]
\end{lemma}

\begin{proof}
This follows from applying Lemma~\ref{keylemma}.
\end{proof}

We conclude

\begin{proposition}\label{thetavanishing}

\begin{itemize}
\item[(i)]
Assume that \emph{one} of the entries of $\Delta_+''$ is nonzero,
then
\[
r_P^{\mathcal{W}}(\varphi^V_{\Delta_+}) =0.
\]

\item[(ii)]
If all of the entries of $\Delta_+''$ vanish, then
\[
r_P^{\mathcal{W}}\left(\varphi^V_{\Delta_+}\right) = 
\varphi^W_{\Delta_+'}.
\]

\end{itemize}

\end{proposition}

\begin{remark}
Analogous results hold for
$r_P^{\mathcal{W}}(\varphi^V_{\Delta_-})$. However, a general
formula for the restriction of
$r_P^{\mathcal{W}}(\varphi^V_{\Delta})$  
is more complicated (and is not needed in this paper).
\end{remark}

\subsection{The Fock model}

It will be convenient to also consider the Fock model $\mathcal{F}=
\mathcal{F}_{n,V}$ of the Weil representation. For more details on $
\mathcal{F}_{n,V}$, see the appendix of \cite{FMII}.

 There is an intertwining map $\iota: S(V^n) \to
\mathcal{P}(\C^{n(p+q)})$ from the $K'$-finite Schwartz functions to
the infinitesimal Fock model of the Weil representation acting on
the space of complex polynomials $\mathcal{P}(\C^{n(p+q)})$ in
$n(p+q)$ variables such that $\iota(\varphi_0) = 1$. We denote the
variables in $\mathcal{P}(\C^{n(p+q)})$ by $z_{\alpha i}$ ($1 \leq
\alpha \leq p$) and $z_{\mu i}$ ($p+1 \leq \mu \leq p+q$) with
$i=1,\dots ,n$. Moreover, the intertwining map $\iota$ satisfies
\begin{align}
\iota \left(  x_{\alpha i} - \frac{1}{2\pi} \frac{\partial}{\partial
x_{\alpha i}}  \right)  \iota^{-1} =  \frac{1}{2\pi i} z_{\alpha
i}, \qquad \qquad 
 \iota \left(  x_{\mu j} - \frac{1}{2\pi} \frac{\partial}{\partial
x_{\mu j}} \right)  \iota^{-1} = - \frac{1}{2\pi i} z_{\mu j}.
\notag
\end{align}
By slight abuse of notation, we use the same symbol for
corresponding objects in the Schr\"odinger and Fock model. In the Fock model, $\varphi_{\Delta}^V$ looks as follows.

\begin{lemma}
\[
\varphi_{\Delta}^V= \prod_{\substack{ 1 \leq \alpha \leq p \\  p+1
\leq \mu \leq m  \\1 \leq j \leq n}}\left( \frac{1}{2\pi i}
z_{\alpha j}\right)^{\delta_{\alpha j}}\left(- \frac{1}{2\pi i}
z_{\mu j}\right)^{\delta_{\mu j}}.
\]
\end{lemma}
Proposition~\ref{thetavanishing} translates to

\begin{proposition}\label{Fock-thetavanishing}
If \emph{one} of the entries of $\Delta_+''$ is nonzero, then
\[
r_P^{\mathcal{W}}\left(\varphi^V_{\Delta_+}\right) =0.
\]
If all of the entries of $\Delta_+''$ vanish, then
\[
r_P^{\mathcal{W}}\left(\varphi^V_{\Delta_+}\right) = 
\prod_{\substack{ \ell+1 \leq \alpha \leq p  \\1 \leq j \leq
n}}\left( \frac{1}{2\pi i} z_{\alpha j}\right)^{\delta_{\alpha j}}.
\]
\end{proposition}

\section{Differential graded algebras associated to the Weil
representation}

In this section, we construct certain differential graded algebras
$C^{\bullet}_V$ and $A_P^{\bullet}$ and define a local restriction map $r_P$ from $C^{\bullet}_V$ to $A_P^{\bullet}$. Again $V$ will denote a
non-degenerate real quadratic space of dimension $m$ and signature
$(p,q)$.

\subsection{Relative Lie algebra complexes}

For convenience of the reader, we briefly review some basic
facts about relative Lie algebra complexes, see e.g., \cite{BW}. For
this subsection, we deviate from the notation of the paper and let
$\frakg$ be any real Lie algebra $\frakg$ and let $\frakk$ be any
subalgebra. We let $(U,\pi)$ be a representation of $\frakg$. We set
\begin{equation}
C^{q}(\frakg,\frakk;U) = \left[ \Hom\left(\sideset{}{^q}\bigwedge
(\mathfrak{g}/\frakk),U\right) \right]^{\frakk}  \simeq \left[
\wwedge{q} (\mathfrak{g}/\frakk)^{\ast} \otimes U \right]^{\frakk},
\end{equation}
where the action of $\frakk$ on $\sideset{}{^q}\bigwedge
(\mathfrak{g}/\frakk)$ is induced by the adjoint representation. Then in the setting of
$\left[\wwedge{q} (\mathfrak{g}/\frakk)^{\ast} \otimes U
\right]^{\frakk}$, the differential $d$ is given by
\begin{equation}\label{Lie-d}
d= \sum_i A(\omega_i) \otimes \pi(X_i) + \frac12 \sum_{i}
A(\omega_i) \ad^{\ast}(X_i) \otimes
 1.
\end{equation}
Here $A(\omega_i)$ denotes the left multiplication with $\omega_i$
in $\wwedge{\bullet} (\frakg/\frakk)^{\ast}$, and $\ad^{\ast}(X)$
is the dual of the adjoint action on $\wwedge{\bullet}
(\frakg/\frakk)^{\ast}$, that is, $
(\ad^{\ast}(X)(\alpha))(Y_1,\cdots,Y_q) = \sum_{i=1}^q
\alpha(Y_1,\dots,[Y_i,X],\dots,Y_q)$. We easily see
\begin{lemma}\label{magic}
Consider two relative Lie algebra complexes
$C^{\bullet}(\frakg,\frakk;U)$ and
$C^{\bullet}(\frakg',\frakk';U')$. Then the following datum,
\begin{itemize}

\item[(i)]

$\rho: \frakg \to \frakg'$, a Lie algebra homomorphism such that
$\rho(\frakk) \subseteq \frakk'$,

\item[(ii)]

$T:U' \to U$, an intertwining map with respect to $\rho$ (i.e., $T(
\rho(X)\cdot u') = X \cdot T(u')$ for $X \in \frakg$),

\end{itemize}
induces a natural map of complexes
\[
C^{\bullet}(\frakg',\frakk';U') \to C^{\bullet}(\frakg,\frakk;U)
\]
given by
\[
\varphi \mapsto T \circ \varphi \circ \rho.
\]
When realizing $\varphi$ as an element  $\left[
\sideset{}{^q}{\bigwedge} (\mathfrak{g'}/\frakk')^{\ast} \otimes U'
\right]^{\frakk'}$, then the map is given by
\[
\rho^{\ast} \otimes T,
\]
where $\rho^{\ast}: (\mathfrak{g'}/\frakk')^{\ast} \to
(\mathfrak{g}/\frakk)^{\ast}$ is the dual map.

\end{lemma}

Now let $G$ be any real Lie group with Lie algebra $\frakg$ and
let $K$ be a closed connected subgroup of $G$ (not necessarily
compact) with Lie algebra $\frakk$. For $U$ a smooth $G$-module, we
let $\calA^{q}(G/K;U)$ be the $U$-valued differential $q$-forms on
$G/K$ (with the usual exterior differentiation). The $G$-action on
$\calA^{q}(G/K;U)$ is given by
\begin{equation}
(g \circ w)_x(X) = g( \omega_{g^{-1} \cdot x}(g^{-1} \cdot X)),
\end{equation}
for $\omega \in \calA^{q}(G/K;U)$, $x \in G/K$, and $X \in
T^q_x(G/K)$. Then evaluation at the base point of $G/K$ gives rise
to an isomorphism of complexes
\begin{equation}
\calA^{\bullet}(G/K;U)^G \simeq C^{\bullet}(\frakg,\frakk;U)
\end{equation}
of the $G$-invariant forms on $G/K$ with
$C^{\bullet}(\frakg,\frakk;U)$.

\subsection{The differential graded algebra
$C^{\bullet}_V$} We begin this section by defining a differential
graded (but not graded-commutative) algebra $C^{\bullet}_V$. The complex  $C^{\bullet}_V$ is essentially the relative Lie algebra
complex for $\Orth(V)$ with values in $\mathcal{W}_{n,V}$ tensored
with the tensor algebra of $V_{\C}$ and twisted by some factors
associated to $\C^n$. Precisely, it is the complex  given by
\begin{align*}
C^{j,r,k}_V &= \left[T^j(U)[-\tfrac{p-q}{2}] \otimes T^{k}(\C^n)^*
\otimes \mathcal{W}_{n,V} \otimes \wwedge{r}
\mathfrak{p}^{\ast}_{\C}
\otimes T^{k}(V_{\C})\right]^{K' \times K \times S_k}\\
& \simeq \left[T^j(U)[-\tfrac{p-q}{2}] \otimes T^k(\C^n)^* \otimes
\mathcal{W}_{n,V} \otimes  \mathcal{A}^{r}(D) \otimes
T^k(V_{\C})\right]^{K' \times G \times S_k }.
\end{align*}
Here $j,r,k$ are nonnegative integers and $\mathcal{A}^r(D)$ denotes
the space of complex-valued  differential $r$-forms on $D$. We let
$U = \bigwedge^n(\C^n)^{\ast}$, and we define the action of $K'$ on
$T^{j}(U)[\frac{p-q}{2}]$ by requiring $K'$ to act by the
character $\det^{-j -\frac{p-q}{2}}$ on  $T^j(U)$. Thus $K'$ acts by
algebra homomorphisms shifted by the character
$\det^{-\frac{p-q}{2}}$. We will usually drop the $[\frac{p-q}{2}]$ in what follows.
The differential is the usual relative Lie algebra differential for the action of
$\Orth(V)$. The group $K'$ acts on the first three factors, while
the maximal compact subgroup $K_V=K$ of $\SO_0(V)$ fixing the
basepoint $z_0$ acts on the last three factors. Finally, the
symmetric group $S_k$ acts on the second and the last factor.

We now give the complex $C^{\bullet}_V$ an associative
multiplication. In order to give the complex the structure of a
graded algebra we choose as a  model for the Weil representation
that has an algebra structure, the Fock model $\mathcal{F}_{n,V}$,
the multiplication law is multiplication of polynomials.  However,
it is important to observe that $K'$ does not act on
$\mathcal{F}_{n,V}$ by algebra homomorphisms (but rather by
homomorphisms twisted by the character $\det^{\frac{p-q}{2}}$). Now
the vector space underlying $C^{\bullet}_V$  is a subspace (of
invariants under a group action) of  a tensor product of graded
algebras. Thus it remains to prove that the group acts by
homomorphisms of the product multiplication.

\begin{lemma}
The group $K' \times K \times S_k$ acts by algebra homomorphisms on
the tensor product of algebras $T^{\bullet}(U) \otimes
T^{\bullet}(\C^n)^* \otimes \mathcal{W}_{n,V} \otimes \wwedge
{\bullet} \mathfrak{p}^{\ast}_{\C} \otimes T^{\bullet}(V_{\C})$.
\end{lemma}
\begin{proof}
The statement is obvious except possibly for the action of the group
$K'$. The group $K'$ acts on the algebra $\mathcal{F}_{n,V}$ by
algebra homomorphisms twisted by the character
$\det^{\frac{p-q}{2}}$. It acts on the tensor product
$T^{\bullet}(U)$ by algebra homomorphisms twisted by the inverse
character $\det^{-\frac{p-q}{2}}$, see e.g. \cite{FMII} Lemma A.1. The two twists cancel on the
tensor product and we find that $K'$ acts by algebra homomorphisms.
\end{proof}

Sometimes it is more convenient to view an element $\varphi \in
C^{j,r,k}_V$ as an element in
\begin{equation}\label{CHom}
\left[\Hom\left(T^{k}(\C^n); T^j(U) \otimes \mathcal{W}_{n,V}
\otimes \wwedge{r} \mathfrak{p}^{\ast}_{\C} \otimes
T^{k}(V_{\C})\right)\right]^{K' \times K \times S_k}.
\end{equation}
For $w \in T^k(\C^n)$, we write $\varphi(w)$ for its value in
$T^j(U) \otimes \mathcal{W}_{n,V} \otimes \wwedge{r}
\mathfrak{p}^{\ast}_{\C} \otimes T^{k}(V_{\C})$.

By Schur-Weyl theory, see \cite{FultonHarris}, Lecture 6, we have
the decomposition
\begin{equation}\label{SWforCn}
T^k(\C^n)^* \simeq  \bigoplus_{\lambda} s_{t(\lambda)}(T^k(\C^n)^*)
\otimes V_{\lambda}^{\ast}.
\end{equation}
Here the sum is over the Young diagrams $\lambda$ with $k$ boxes and
no more than $n$ rows,  $t(\lambda)$ is a chosen standard filling of
$\lambda$ for each  $\lambda$ and $V_{\lambda}$ is the irreducible
representation of $S_k$ corresponding to $\lambda$. We also have the
corresponding decomposition
\begin{equation}
T^k(V_{\C}) \simeq \bigoplus_{\mu} s_{t'(\mu)}( T^k(V_{\C}))\otimes
V_{\mu}.
\end{equation}
Combining the two decompositions we obtain
\begin{align}
C^{j,r,k}_V \simeq  \bigoplus_{\lambda,\mu} \left[T^j(U) \otimes
\mathbb{S}_{t(\lambda)}(\C^n)^* \otimes V_{\lambda}^{\ast} \otimes
\mathcal{W}_{n,V} \otimes \wwedge{r} \mathfrak{p}^{\ast}_{\C}
\otimes \mathbb{S}_{t'(\mu)}(V_{\C}) \otimes V_{\mu}\right]^{K'
\times K \times S_k}.
\end{align}
Noting that \begin{equation}
(V_{\lambda}^{\ast} \otimes V_{\mu})^{S_k} \simeq \begin{cases} 0 & \text{if}  \ \lambda \neq \mu \\
\C & \text{if} \ \lambda = \mu,
\end{cases}
\end{equation}
we obtain
\begin{lemma}\label{SchurWeylforC}
\begin{align*}
 C^{j,r,k}_V \simeq \bigoplus_{\lambda} \left[T^j(U) \otimes
\mathbb{S}_{t(\lambda)}(\C^n)^{\ast}\otimes \mathcal{W}_{n,V}
\otimes \wwedge{r} \mathfrak{p}^{\ast}_{\C} \otimes
\mathbb{S}_{t(\lambda)}(V_{\C})\right]^{K'\times K}.
\end{align*}
\end{lemma}
We have assumed (as we may do) that the fillings $t(\lambda)$ and
$t'(\lambda)$ are the same. For the summands in the lemma we write
$C^{j,r,t(\la)}_V$ (or just $C^{j,r,\la}_V$ if we do not want to
specify the filling). The application of the Schur functor
$\mathbb{S}_{t(\la)}^{\ast} (\cdot)$ on $T^k(\C^n)^{\ast}$ or
equivalently applying $\mathbb{S}_{t(\la)}(\cdot)$ on $T^k(V_{\C})$,
gives rise to a projection map
\begin{equation}\label{SWC}
\pi_{t(\la)}: C^{j,r,k}_V \longrightarrow C^{j,r,t(\la)}_V.
\end{equation}
That is,
\begin{align}
\pi_{t(\la)} & = 1_U \otimes {s_{t(\la)}}_{(\C^n)^{\ast}} \otimes
1_{\calW_{n,V}} \otimes 1_{\mathfrak{p}^{\ast}} \otimes 1_V \\
& = 1_U \otimes 1_{\C^n} \otimes 1_{\calW_{n,V}} \otimes
1_{\mathfrak{p}^{\ast}} \otimes {s_{t(\la)}}_V. \notag
\end{align}
Here we use subscripts to indicate which spaces the respective
identity transformations $1$ operate on. 
We apply the harmonic projection $\mathcal{H}_V$, see
\eqref{harmonic}, on the last factor to obtain
$\mathbb{S}_{[t(\la)]}(V_{\C})$, and we obtain a complex
$C^{\bullet,[t(\la)]}_V$ (or $C^{\bullet,[\la]}_V$) and a projection
map
\begin{equation}\label{SWC1}
\pi_{[t(\la)]}: C^{j,r,k}_V \longrightarrow C^{j,r,[t(\la)]}.
\end{equation}
That is,
\begin{align}
\pi_{[t(\la)]} &= 1_U \otimes 1_{\C^n} \otimes 1_{\calW_{n,V}}
\otimes 1_{\mathfrak{p}^{\ast}} \otimes {s_{[t(\la)]}}_V.
\end{align}

\begin{remark}\label{SchurWeylforC1}
We can interpret an element $\varphi \in C^{j,r,k}_V$ as a $K'
\times K \times S_k$-invariant homomorphism from $T^{k}(\C^n)$ to
$T^j(U)\otimes  \mathcal{W}_{n,V} \otimes {\bigwedge} {^r}
\mathfrak{p}^{\ast}_{\C} \otimes T^{k}(V_{\C})$, see \eqref{CHom}.
In this setting, we can interpret $\pi_{t(\la)}\varphi$ as the
restriction of the homomorphism $\varphi$ to the
$\mathbb{S}_{t(\la)}(\C^n)$. From this point of view,
Lemma~\ref{SchurWeylforC} states that the homomorphism
$\pi_{t(\la)}\varphi$ on $\mathbb{S}_{t(\la)}(\C^n)$ automatically
takes values in $\mathcal{W}_{n,V} \otimes {\bigwedge} {^r}
\mathfrak{p}^{\ast}_{\C} \otimes \mathbb{S}_{t(\la)}(V_{\C})$.
\end{remark}

\subsection{The face differential graded algebra $A^{\bullet}_P$ and the map $r_P$}
In this section we assume $P$ is the stabilizer of a standard flag
$E_{i_1} \subset E_{i_2} \subset \cdots \subset E_{i_k}=E_{\ell}=E$
and $N_P$ is the unipotent radical of $P$. We let $Q$ be the stabilizer of $E$. We will now construct a differential graded
algebra $A^{\bullet}_P$, which is the relative Lie algebra version
of a differential graded subalgebra of the de Rham complex of the
face $e(P)$ of the Borel-Serre enlargement of $D$.
We will continue with the notation of section 2.

We define the differential graded algebra $A^{\bullet}_P$ associated
to the face $e(P)$ of the Borel-Serre boundary corresponding to $P$
by
\begin{align}
A^{j,r,k}_P   &= \left[T^j(U) \otimes T^{k}(\C^n)^* \otimes
\mathcal{W}_{n,W} \otimes  \wwedge{r} (\mathfrak{n} \oplus
\mathfrak{p}_M)^*_{\C} \otimes T^{k}(V_{\C})\right]^ {K' \times K_P
\times S_k} \\
& \simeq \left[ T^j(U) \otimes T^{k}(\C^n)^* \otimes
\mathcal{W}_{n,W} \otimes \mathcal {A}^{r}(e(P)) \otimes
T^{k}(V_{\C})\right]^{K'\times NM \times S_k}. \notag
\end{align}
Furthermore, we define $A^{\bullet,\la}_P$ and $A^{\bullet,[\la]}_P$
as for $C^{\bullet}_V$.

\begin{definition}\label{def:localres}
The "local" restriction map of de Rham algebras with coefficients
\[
r_P: C^{\bullet}_V  \to A^{\bullet}_P
\]
of de Rham algebras with coefficients is given by
\[
1\otimes 1 \otimes r_P^{\mathcal{W}}  \otimes \iota^{\ast} \otimes
1.
\]
Here $ \iota:\mathfrak{n} \oplus \mathfrak{m} \hookrightarrow \frakg$
is the underlying Lie algebra homomorphism, and the map from the
coefficients of $C^{\bullet}_{V}$ to the coefficients of
$A^{\bullet}_P$ is given by the tensor product
\[
1\otimes 1 \otimes r_P^{\mathcal{W}} \otimes 1,
\]
where $r_P^{\mathcal{W}}: \mathcal{W}_{n,V} \to \mathcal{W}_{n,W}$
is the restriction map of the Weil representation (see
Definition~\ref{def:localresWeil}). By Lemma~\ref{magic} we
therefore see that $r_P$ is a map of complexes. We note that
$r_P^{\mathcal{W}} $ is not a ring homomorphism so $r_P$ is not a
map of algebras. Since $r_P$ commutes with the action of the
symmetric group $S_k$, we obtain maps $C^{\bullet,\la}_V  \to
A^{\bullet,\la}_P$ and $C^{\bullet,[\la]}_V  \to
A^{\bullet,[\la]}_P$ as well, which we also denote by $r_P$.
\end{definition}

Note that the induced map $\iota^{\ast}: (\frakg/\frakk)^{\ast}
\simeq \frakp^{\ast} \to ((\frakn \oplus \frakm)/\frakk_M)^{\ast}
\simeq (\mathfrak{n} \oplus \mathfrak{p}_M)^*$ is the composition of
the isomorphism $\sigma^{\ast}: \mathfrak{p}^{\ast} \longrightarrow
(\mathfrak{n} \oplus \mathfrak{a} \oplus \mathfrak{p}_M)^{\ast}$,
see \eqref{pnam-iso}, with the restriction $ (\mathfrak{n} \oplus
\mathfrak{a}\oplus \mathfrak{p}_M)^{\ast} \to (\mathfrak{n} \oplus
\mathfrak{p}_M)^{\ast}$.

Finally observe that on the level of homogeneous spaces, the map
$r_P$ arises by realizing $e(P)$ as the orbit of the basepoint $z_0$
under the group $NM$. So in this setting, we are no longer thinking
of $e(P)$ as being at the boundary of $D$; we have pushed $e(P)$ far
inside $D$.

\section{Aspects of nilpotent Lie algebra cohomology and the map $\iota_P$}\label{nilcoh}

\subsection{An explicit constituent in the Lie algebra cohomology of $\frak{n}_P$}

In what follows, we discuss some aspects of the Lie algebra cohomology of the nilpotent Lie algabra $\frak{n}_P$ which we need later. Some parts we develop here could have been deduced from the general work of Kostant \cite{Kostant}. However, our concern here is proving that certain explicit cocycles are (non)zero rather than computing the cohomology itself.

As before, we let $P$ be a standard parabolic subgroup of $G$. Recall that we have the decomposition of vector spaces $\frak{n}_P = \frak{n}_{P'} \oplus \frak{n}_Q$, where $Q$ is the maximal parabolic containing $P$. For the two-step nilpotent algebra $\frak{n}_Q$,  we have the central extension $\mathfrak{z}_Q \to \mathfrak{n}_Q\to \frakn_W$ with $\mathfrak{z}_Q \simeq \wwedge{2}E$ and  $\frakn_W \simeq W\otimes E$. On the other hand, $\frakn_{P'}$ is a nilpotent subgroup of $\frak{sl}(E) \subset E' \otimes E$.

We assume for the next subsections that $V,W,\frak{n}_P$ etc. are defined over $\C$. We let 
\[
\calC^{\bullet,\ell'} = \wwedge{\bullet}( \frakn_P^{\ast}) \otimes T^{\ell'}(V)
\]
be the complex for the nilpotent cohomology with coefficients in $T^{\ell'}(V)$ and define analogously $\calC^{r,\la}=\calC^{r,A}$ and $\calC^{r,[\la]} = \calC^{r,[A]}$ for $\mathbb{S}_{\la}(V)$ and $\mathbb{S}_{[\la]}(V)$ respectively.

We are interested in certain cohomology classes arising from $\wwedge{r} \frakn_W^{\ast}$. By Lemma~\ref{nP-decomposition} $\frakn_W^{\ast}
\simeq W \otimes E'$ as $\Orth(W) \times \GL(E)$-modules.
Furthermore (e.g. \cite{FultonHarris}, p. 80),
\begin{equation}\label{FH-decomp}
\wwedge{r} (\frakn_W^{\ast}) \simeq  \wwedge{r} \left(W
\otimes E' \right) \simeq \bigoplus_{\mu}
\mathbb{S}_{\mu}(W) \otimes \mathbb{S}_{\mu'}(E'),
\end{equation}
as $\Orth(W) \times \GL(E)$-modules. Here the sum extends
over all partitions $\mu$ of $r$ with at most $\dim W = m- 2\ell$
rows and at most $\dim E = \ell$ columns, and $\mu'$ denotes the
conjugate partition of $\mu$.

We will be mainly interested in the case $r=n\ell$. Then we can take
$\mu=\ell\varpi_n=(\ell,\ell,\dots,\ell)$, so that $\mu' = n
\varpi_{\ell}= (n,n,\dots,n)$ and $\mathbb{S}_{\mu'}(E') =
\left( \wwedge{\ell}E'\right)^{\otimes n} \simeq \C$ is the 
trivial (one-dimensional) $\SL(E)$-module. We obtain
\begin{equation}\label{SB-iso}
\mathbb{S}_B(W) \otimes \mathbb{S}_{B'} (E') \simeq \left[ \wwedge{n \ell}(W \otimes E') \right]^{SL(E)}
\simeq \left[ \wwedge{n\ell}
(\frakn_W^{\ast})_{\C} \right]^{\SL(E)}
\end{equation}
as $\Orth(W) \times \SL(E)$-modules. Here $B=B_{n,\ell}$ is the filling of the
Young diagram associated to $\mu$ described in
section~\ref{box-section}.

To realize this isomorphism, we define a $\GL(W) \times \GL(E)$ intertwining map 
\begin{equation}
\tau_{r,\ell'} : T^{r}(W) \otimes T^{\ell'}(W) \otimes T^r(E') \to \wwedge{r}\left(W \otimes E'\right) \otimes T^{\ell'}(V)  \subset \calC^{r,\ell'}
\end{equation}
given by
\begin{align}
\tau_{r,\ell'}( (w_1 \otimes \cdots \otimes w_r) \otimes  {\bf \bar{w}} &\otimes (v'_1 \otimes \cdots \otimes v'_r)) = (w_1 \otimes v'_1) \wedge \cdots \wedge (w_r \otimes v'_r) \otimes {\bf \bar{w}}, \notag 
\end{align}
where ${\bf \bar{w}} \in T^{\ell'}(W)$. We also write $\tau_r$ for $\tau_{r,0}$. We immediately see

\begin{lemma}\label{tau-comm}
The map $\tau_{r,\ell'}$ is $\Orth(W) \times \SL(E) \times S_{r+\ell'} \times S_{r}$-equivariant. Here the
action of the symmetric group $S_{r+\ell'}$(respectively $S_r$) is on the tensor factors involving $W$(respectively $E'$). 
\end{lemma}

For $r=n\ell$, the map $\tau_{n\ell}$ realizes the isomorphism \eqref{SB-iso}. Furthermore,

\begin{lemma}\label{SLE-invariance}
Let $\w \in T^{n\ell+\ell'}(W)$ and ${\bf v}' \in T^{n\ell}(E')$. Then
\[
\tau_{n\ell,\ell'}(s_{B|A} (\w) \otimes {\bf v}') \in \left(\mathcal{C}^{n\ell,A}\right)^{\SL(E)}.
\]
\end{lemma}

We view from now on $\tau_{n\ell,\ell'}$ as a map of $T^{n\ell+\ell'}(W)$ by setting 
\[
\tau_{n\ell,\ell'}({\bf w}) := \tau_{n\ell,\ell'}({\bf w} \otimes  (u_1' \otimes \cdots \otimes u_{\ell}')^{n}).
\]

We let $V^{[k]}$ ($W^{[k]}$) be the space of harmonic $k$-tensors in $V$ ($W$); i.e., the tensors which are annihilated by all the contractions $C_{ij}$. We let $\mathcal{E}^k(V) \subset T^{k}(V)$ be the orthogonal complement of the harmonic tensors. Thus $\mathcal{E}^k(V)$ is the sum of the images of the insertion maps $E_{ij}(g_V^{\ast}):T^{k-2}(V) \to T^k(V), 1 \leq i < j \leq k$ with the metric $g_V^{\ast}$ of $V$. Similarly, we define $\mathcal{E}^k(W) \subset T^{k}(W)$. Note $\mathbb{S}_{[\la]}(W) \subset \mathbb{S}_{[\la]}(V)$. However note, that if ${\bf \bar{w}} \in T^{\ell'}(W)$ is a nonzero tensor in the orthogonal complement of $T^{[\ell']}(W)$ (i.e., spanned by tensors in the image of the inclusion with the metric for $\mathcal{E}^{\ell'}(W)$), then ${\bf \bar{w}}$ does \emph{not} necessarily lie in the orthogonal complement in $T^{[\ell']}(V)$ (since the metric of $V$ is different).

\begin{proposition}\label{nonvanish}
Let $B$ again be the given filling of the Young diagram associated to $\ell \omega_n$ and $A$ be a filling for $\la$. 

\begin{itemize}
\item[(i)]
Let ${\bf w}  \in \mathbb{S}_{B|A}(W)$. Then $\tau_{n\ell,\ell'}({\bf w})$ defines a cocycle in $\calC^{n\ell,\ell'}$. More precisely, we obtain a map
\[
\mathbb{S}_{B|A}(W) \to H^{n\ell}(\frakn_P,\mathbb{S}_A(V))^{\SL(E)}.
\]

\item[(ii)]
Let $n \leq \left[ \frac{\dim W}{2} \right]$ and let ${\bf w} \in \mathbb{S}_{[B|A]}(W) $. Then the cohomology class
\[
[\tau_{n\ell,\ell'}\left({\bf w} \right)] \in H^{n\ell}(\frakn_P,\mathbb{S}_{[A]}(V))^{\SL(E)}
\]
does not vanish.Thus we obtain an embedding
\[
\mathbb{S}_{[B|A]}(W) \hookrightarrow H^{n\ell}(\frakn_P,\mathbb{S}_{[A]}(V))^{\SL(E)}
.
\]
\item[(iii)]
Let ${\bf w} \in \mathbb{S}_{B|A}(W) \cap \mathcal{E}^{n\ell+\ell'}(W)$ be in the orthogonal complement of $\mathbb{S}_{[B|A]}(W)$ inside $\mathbb{S}_{B|A}(W)$. Then
\[
[\pi_{[A]} \circ \tau_{n\ell,\ell'}({\bf w})] =0 
\]
in $ H^{n\ell}(\frakn_P,\mathbb{S}_{[A]}(V))$. Here $\pi_{[A]}$ is the natural projection from $H^{\bullet}(\frakn_P,\mathbb{S}_A(V))$ to $H^{\bullet}(\frakn_P,\mathbb{S}_{[A]}(V))$ induced by the orthogonal projection $\mathbb{S}_{A}(V) \to \mathbb{S}_{[A]}(V)$. In particular, for ${\bf w}  \in \mathbb{S}_{B|A}(W)$, we have
\[
[\pi_{[A]} \circ \tau_{n\ell,\ell'}({\bf w})] = [\tau_{n\ell,\ell'} (\pi_{[B|A]}(\bf w))].
\]

\end{itemize}

\end{proposition}

The next section will be concerned with proving this proposition.

\subsection{Proof of Proposition~\ref{nonvanish}}

We give $V$ the positive definite Hermitean metric coming from the majorant $(\ ,\ )_0$. This induces positive definite metrics on $\bigwedge^2E$, $W\otimes E$, and $E' \otimes E$ and hence an admissible metric on the entire Lie algebra complex $\calC^{\bullet,\ell'}$, which we also denote by $(\ ,\ )_0$. Using $(\ ,\ )_0$ we obtain an adjoint $d^*$ to the differential $d$  on $\calC^{\bullet}$. We then have the finite-dimensional analogues of Hodge theory. Namely,  we define the Laplacian $\Delta = d d^* + d^* d$ and say a form
in $\calC^{\bullet}$ is harmonic if it is in the kernel of Delta. It is immediate that
\[
\ker \Delta = \ker d \cap \ker d^*.
\]
We let $\mathcal{H}^{r,\ell'}$ be the harmonic forms of degree $r,$ be the intersection $\ker \Delta \cap \mathcal{C}^{r,\ell'}$. In particular, we have the Hodge decomposition

\begin{lemma}\label{Hodge}
The space $\calC^{r,\bullet}$ is the orthogonal direct sum
of the exact forms $\im d$, the coexact forms $\im d^*$ and the harmonic
forms. Furthermore, the map $\mathcal{H}^{r,\bullet} \to H^{r,}(\mathfrak{n}_P,\bullet)$ is an isomorphism.
\end{lemma}

The Lie algebra complex $\calC^{\bullet,\ell'}$ is in fact triple-graded via 
\[
\calC^{r,s,t,\ell'} := \wwedge{r}(W\otimes E') \otimes \wwedge{s}\left(\wwedge{2}E'\right) \otimes 
 \wwedge{t} \frakn^{\ast}_{P'} \otimes T^{\ell'}(V)
\]
and define analogously $\calC^{r,s,t,\la}$ and $\calC^{r,s,t,[\la]}$ for $U =\mathbb{S}_{\la}(V)$ and $U=\mathbb{S}_{[\la]}(V)$ respectively. Here again we have used the form $(\,,\,)$ to identify $W^{\ast} \simeq W$ and $E^{\ast} \simeq E'$. 

We now give explicit formulas for the Lie algebra differential $d$ and its adjoint $d^{\ast}$ on $\calC$. We omit the proofs. We write $d=d_{\frakn} + d_V$ with a "Lie algebra part" $d_{\frakn}$ and a "coefficient" part  $d_V$. That is,
\begin{equation}\label{dV}
d_{\frakn} = d_{\frakn_Q} + d_{\frakn_{P'}} \qquad \text{and} \qquad d_V= d_V^W + d_V^E +d_V^{\frakn_{P'}}
\end{equation}
with
\begin{align}\label{dn}
d_{\frakn_Q} &= \frac12 \sum_{\alpha, i} A(e_{\alpha} \otimes u_i') ad^{\ast}(e_{\alpha} \wedge u_i') 
 + \frac12 \sum_{1 \leq i <j \leq \ell} A(u'_i \wedge u'_j) ad^{\ast}(u_i \wedge u_j)
 \end{align}
and
\begin{align*}
d^W_{V}  = \sum_{\alpha, i} A(e_{\alpha} \otimes u_i') \otimes \rho(e_{\alpha} \wedge u_i) \qquad \text{and} \qquad 
d_V^E = \sum_{1 \leq i <j \leq \ell} A(u'_i \wedge u'_j) \otimes \rho (u_i \wedge u_j).
\end{align*}
Here $\rho$ denotes the action of $\frakn_P$ on the coefficient system $T^{\ell'}(V)$. Finally, $d_{\frakn_{P'}}+d_V^{\frakn_{P'}}$ is the part of the differential arising from $\frakn_{P'}$. (We don't need it more precisely).

Since $[\frakn_W,\frakn_W] \subseteq \mathfrak{z} _Q$, we first note that $d_{\frakn_Q}$ has triple-degree (2,-1,0). In particular, all elements of degree $(r,0,t)$ are $d_{\frakn_Q}$-closed. Accordingly, $d_{\frakn_Q}$ is determined by its values on $\mathcal{C}^{0,s,0,\bullet}$. In fact, it suffices to consider $s=1$.

\begin{lemma}\label{exteriorderivative}
Let $v_1',v_2' \in E'$ and ${\bf v} \in T^{\ell'}(V)$. Then
\begin{align*}
d_{\frakn_Q} \left( ( v_1'\wedge v_2')  \otimes {\bf v}\right)= 
-\tau_{2}\left( E_{1,2}(g_W^{\ast}) \otimes (v'_1 \otimes v'_2) \right) \otimes {\bf v}.
\end{align*}
\end{lemma}

It suffices to compute the dual $d^*_{\frakn_Q}$ on basic forms. 
\begin{proposition}\label{adjointn}
\begin{align*}
&d^*_{\frakn_Q} ( (w_1\otimes v_1')\wedge  \cdots \wedge (w_k \otimes v_k')  \otimes {\bf v}
) \\
&= \sum_{ i < j } \biggl\{ (-1)^{i+j}(w_i,w_j)   
(w_1\otimes v_1')\wedge  \cdots \wedge (\widehat{w_i \otimes v_i'}) \wedge \cdots \wedge (\widehat{w_j \otimes v_j'}) \wedge \cdots \wedge (w_k \otimes v_k')\\
& \hskip 12 cm \otimes (v_i' \wedge v_j')\biggr\} \otimes {\bf v}. 
\end{align*}
\end{proposition}
   
 For the $\frakn_{P'}$-contribution, we have the following:
 
\begin{lemma}\label{nP'}
The differential $d_{\frakn_{P'}}+d_V^{\frakn_{P'}}$ has triple-degree $(0,0,1)$. The adjoint action of $\frakn_{P'} \subset \mathfrak{sl}(E)$ on $\frakn_Q =  (W \otimes E) \oplus  \wwedge{2}E$ arises from the natural action of $ \mathfrak{sl}(E)$ on $E$. Hence $d_{\frakn_{P'}}+d_V^{\frakn_{P'}}$ vanishes on $\left(\mathcal{C}^{r,s,0,\ell'}\right)^{\SL(E)}$. In particular, 
\[
(d_{\frakn_{P'}}+d_V^{\frakn_{P'}}) \tau_{n\ell,\ell'}(s_{B|A}(\w)) =0
\]
for $\w \in T^{n\ell+\ell'}(W)$. Finally, the dual $d_{\frakn_{P'}}^*+\left(d_V^{\frakn_{P'}}\right)^*$ vanishes on $\mathcal{C}^{r,s,0,\ell'}$. 
\end{lemma}

We now turn our attention to $d_V$ and $d_V^{\ast}$. It suffices to consider the case $\ell'=1$.

\begin{lemma}\label{dV-formula}

\begin{itemize}

\item[(i)] Let ${\bf w} \in T^{k}(W)$, ${w} \in W$, and ${\bf v'} \in T^k(E')$. Then
\begin{align*}
& d_V^W\left(\tau_{k,1}( {\bf w} \otimes w \otimes {\bf v'})\right) =  \sum_{i=1}^{\ell} \tau_{k+1}
 \left( (w \otimes {\bf w})  \otimes (u_i \otimes {\bf v'}) \right) \otimes u_i.
\end{align*}

\item[(ii)] Let  ${\bf w} \in T^{k}(W)$, ${\bf v'} \in T^k(E')$, and $u' \in E'$. Then
\begin{align*}
&d_V^W \left( \tau_k ({\bf w} \otimes {\bf v'}) \otimes u'\right)  = - \tau_{k+1,1} \left(    E_{1,k+1}(g_W^{\ast}) ({\bf w} ) \otimes (u' \otimes {\bf v'}) \right).
\end{align*}

\end{itemize}

\end{lemma}

\begin{lemma}\label{adjointV}
Let ${\bf w} \in T^{k}(W)$, ${w} \in W$, and ${\bf v'} \in T^k(E')$. Then
\begin{align*}
&d^{\ast}_V \left(\tau_{k,1}( {\bf w} \otimes w  \otimes {\bf v'}) \right) \\
&= \sum_{i=1,\dots,k} (-1)^{i-1} (w_i,w) \left( (w_1\otimes v_1')\wedge  \cdots \wedge (\widehat{w_i \otimes v_i'}) \wedge \cdots \wedge (w_k \otimes v_k')  \right)\otimes v_k'.
\end{align*}

\end{lemma}

As a consequence of Lemma~\ref{adjointn}, Lemma~\ref{adjointV}, and Lemma~\ref{nP'} we  obtain
   
\begin{proposition}\label{d-ast-key}
Let ${\bf w} \in W^{[k+\ell']}$ be a harmonic $(k+\ell')$-tensor. Then for any ${\bf v'} \in T^k(E')$, we have
\[
d^{\ast} \tau_{k,\ell'}({\bf w} \otimes {\bf v'}) = 0.
\]
\end{proposition}

We are now ready to prove Proposition~\ref{nonvanish}. For (i), first note that the action of $\sigma \in S_{\ell'}$ on the coefficients $T^{\ell'}(V)$ commutes with the differentiation $d$: $d\circ (1 \otimes \sigma \otimes 1) = (1 \otimes \sigma \otimes 1) \circ d$.  Furthermore, in the first factor $T^{n\ell}(W)$, $\tau_{n\ell,\ell'}$ factors through $c_B$, the column anti-symmetrizer for Young tableau $B$, that is, $\tau_{n\ell,\ell'} \circ (c_B \otimes 1) = \tau_{n\ell,\ell'}$. Combining this with Lemma~\ref{tau-comm} gives $
\tau_{n\ell,\ell'} \circ (c_{B|A}) = (1 \otimes c_A) \circ \tau_{n\ell,\ell'}$
on $T^{n\ell+\ell'}(W)$. Therefore it suffices to show that $\tau_{n\ell,\ell'}(r_{B|A}({\bf w}))$ is closed. Indeed, we have
\[
d \left(\tau_{n\ell,\ell'}(s_{B|A}(\w)) \right) = d\left( (1 \otimes c_A) \circ \tau_{n\ell,\ell'}(r_{B|A}(\w))\right) = (1 \otimes c_A) \circ d(\tau_{n\ell,\ell'}(r_{B|A}(\w))).
\]

Furthermore, it suffices to establish closedness for $n=1$. Indeed, if the Young diagram $A$ arises from the partition $(\ell_1',\ell_2',\dots,\ell_n')$ of $\ell'$, we write $\w = \w_1 \otimes \cdots \otimes \w_n \in T^{n\ell}(W)$ with $\w_i \in T^{\ell}(W)$ and ${\bf \bar{w}} = {\bf \bar{w}}_1 \otimes \cdots \otimes {\bf \bar{w}}_n$ with ${\bf \bar{w}}_i \in T^{\ell_i'}(
W)$. We then have a natural product decomposition
\begin{equation}\label{p-dec}
\tau_{n\ell,\ell'}(\w \otimes {\bf \bar{w}}) = \tau_{\ell,\ell_1'}(\w_1 \otimes  {\bf {\bar{w}}}_1) \wedge \cdots \wedge 
 \tau_{\ell,\ell_n'}(\w_n\otimes  {\bf {\bar{w}}}_2),
\end{equation}
for which $d$ acts as a derivation. Also note that $d_{\mathfrak{n}_Q}$ vanishes on the image of $\tau_{n\ell,\ell'}$ and by Lemma~\ref{nP'} so does the $\frakn_{P'}$-contribution. Now for $n=1$, using Lemma~\ref{dV-formula} (i), we see that applying $d_V$ to $\tau_{\ell,\ell'}(\w)$ with $\w \in \Sym^{\ell +\ell'}(W)$ gives rise to a map
\begin{align}\label{LR}
\Sym^{\ell+\ell'}(W) & \to \bigoplus_{i=1}^{\ell} \wwedge{\ell+1} (W \otimes E') \otimes (E'_{i} \otimes T^{\ell'-1}(W))\\& \quad = \bigoplus_{i=1}^{\ell} \bigoplus_{C} \mathbb{S}_{C}(W) \otimes \mathbb{S}_{C'}(E') \otimes (E'_i \otimes T^{\ell'-1}(W)). \notag 
\end{align}
Here $E'_i = {\C} u_i'$, and the sum extends over all Young diagrams $C$ of size $\ell+1$, which have at least $2$ rows (otherwise the dual diagram $C'$ would have at least $\ell+1$ rows, which is impossible as $\dim E' = \ell$). By the Littlewood-Richardson rule we now see that in the decomposition of $\mathbb{S}_{C}(W) \otimes T^{\ell'-1}(W)$ into irreducibles only Young diagrams with at least $2$ rows can occur. Hence $\Sym^{\ell+\ell'}(W)$ does not occur on the right hand side of \eqref{LR}, and the map vanishes identically. This proves Proposition~\ref{nonvanish}(i). 

Proposition~\ref{nonvanish}(ii) now follows immediately from Proposition~\ref{d-ast-key} and Lemma~\ref{Hodge}.

For (iii), it suffices to show that for any $\w \in T^{n\ell+\ell-2}$, the form $\pi_{[A]} \circ \tau_{n\ell,\ell'}(s_{B|A}(E_{i,j}(g_W^{\ast})(\w)))$ is exact. For this, it suffices to show that $\tau_{n\ell,\ell'}(r_{B|A}E_{i,j}(g_W^{\ast})(\w)))$ is exact up to terms involving the inclusion of the metric $g_V^{\ast}$ into the coefficient system. The product decomposition \eqref{p-dec} reduces the claim to the cases of $n=1$ (in case the metric $g_W^{\ast}$ occurs in one factor for \eqref{p-dec}) or $n=2$ (if $g_W^{\ast}$ occurs in two factors). It is not too hard but a bit tedious to explicitly construct primitives for these cases. We omit this.

\subsection{The map $\iota_P$}

We now assume again that all objects are defined over $\R$. We construct a map $\iota_P: C_W^{\bullet} \hookrightarrow A_P^{\bullet}$ of complexes.

 We let $U, U'$ be two representations of $G$ and $T:U' \to U$ be  $G$-intertwiner. We let $\mathcal{C}^{\bullet}(\frakn_P,U)=  (\wwedge{\bullet} \frakn_P^*) \otimes U)$ be the complex computing the nilpotent cohomology $\mathcal{H}^s(\frakn_P,U)$, and we let $\mathcal{C}_{\text{closed}}^{\bullet}(\frakn_P,U)$ be the subspace of cocycles in $\mathcal{C}^{\bullet}(\frakn_P,U)$. 

\begin{lemma}\label{Harderlemma}

Define a map 
\[
\eta^{r,s}:\left[ \wwedge{r}(\frakp_M^*)\otimes \left( (\wwedge{s} \frakn_P^*) \otimes U' \right) \right]^{K_P}  \to 
\left[ \wwedge{r+s}(\frakp_M^* \oplus \frakn_P^*) \otimes U \right]^{K_P} 
\]
by 
\[
\eta^{r,s} (\omega^{(r)} \otimes (\omega^{(s)} \otimes u') = (\omega^{(r)} \wedge \omega^{(s)}) \otimes T(u').
\]
Then $\eta^{r,s}$ induces a map of relative Lie algebra complexes
\[
\eta: \mathcal{C}^{\bullet} \left( \mathfrak{m}, \mathfrak{k}_P; \mathcal{C}_{\text{closed}}^{s}(\frakn_P,U') \right) \longrightarrow \mathcal{C}^{\bullet+s} \left( \mathfrak{p}, \mathfrak{k}_P; U \right)
\]
and the induced map in cohomology factors through $H^{\bullet} \left( \mathfrak{m}, \mathfrak{k}_P;H^{s}(\frakn_P,U') \right)$.

\end{lemma}

\begin{proof}
This is essentially in \cite{Harder}, Lemma~2.6, see also \cite{Schwermer}, section 2, together with the standard spectral sequences argument in this context. Note that Harder actually considers instead of cocycles in $\mathcal{C}(\frakn_P,U')$ the nilpotent cohomology group $\mathcal{H}^s(\frakn_P,U)$ realized as subspace in $\mathcal{C}(\frakn_P,U')$ by harmonic forms as discussed in section~\ref{nilcoh}.
\end{proof}

\begin{definition}

We define the map $\iota_P$ on $C_W^{j,r,k}$ as follows. In fact, it
is defined on the underlying tensor spaces without taking the group
invariants. First we set $\iota_P$ to be zero if $k < n \ell$. If $k
\geq n \ell$ we split the two tensor factors
\[
T^{k}(\C^n)^* = T^{n\ell}(\C^n)^* \otimes T^{k-n\ell}(\C^n)^* \quad
\text{and} \quad T^{k}(W_{\C}) = T^{n\ell}(W_{\C}) \otimes
T^{k-n\ell}(W_{\C}).
\]
We define $\iota_P$ on tensors which are decomposable relative to
these two splittings. We let $u_1=\theta_1 \wedge \cdots \wedge
\theta_n$ be the standard generator of $U = \bigwedge ^n
(\C^n)^{\ast}$ (with the twisted $K'$-action). Let $u_1^{j} \otimes
x \otimes f \otimes \omega \otimes w$ be a single tensor
component of an element in $C^{j,r,k}_W$ and assume that $k \geq
n\ell$. Assume that $x$ and $w$ are decomposable, that is
\[
x = x_1 \otimes x_2 \in T^{n\ell}(\C^n)^* \otimes
T^{k-n\ell}(\C^n)^* \quad  \text{and} \quad  w = w_1 \otimes w_2 \in
T^{n\ell}(W_{\C}) \otimes T^{k-n\ell}(W_{\C}).
\]
Then we define
\begin{align*}
 \iota_P(u_1^{j} & \otimes x \otimes f \otimes \omega
\otimes w) \\
& = (-1)^{n\ell (\frac{(q-\ell)(n-1)}{2}+1)} 
\eta^{r,n\ell} \biggl((u_1^{j} \otimes s_B^*(x_1))
\otimes x_2 \otimes f \otimes
\omega \otimes \tau_{n\ell}(w_1) \otimes w_2 \biggr) \notag \\
& \quad \in T^{j+\ell}(U) \otimes T^{k-n\ell}(\C^n)^{\ast} \otimes
\mathcal{W}_{n,W} \otimes \wwedge{r} (\mathfrak{p}_W^*)_{\C} \otimes
\wwedge{{n\ell}} (\mathfrak{n}_W^{\ast})_{\C} \otimes
T^{k-n\ell}(W_{\C}). \notag
\end{align*}
Note here that by Lemma~\ref{basisforonedimensionalspace}, we see
that $ \mathbb{S}_B(\C^{n})^{\ast} =
s_B^{\ast}T^{n\ell}(\C^{n})^{\ast} \simeq T^{\ell}(U)[0]$ and
therefore $u_1^{j} \otimes s_B^*(x_1)$ lies in
$T^{j+\ell}(U)[-\tfrac{p-q}{2}]$ and is zero if and only if
$s_B^*(x_1)=0$. 

\end{definition}

\begin{proposition}
$\iota_P$ is a map of complexes
\[
\iota_P:C_W^{j,r,k} \to A_P^{j+\ell,r+n\ell,k-n\ell}.
\]
\end{proposition}

\begin{proof}

In view of Lemma~\ref{Harderlemma}, it suffices to show that the map on $C_W^{j,r,k}$ to 
\begin{equation}
 \mathcal{C}^{r} \left( \mathfrak{m}, \mathfrak{k}_P; \mathcal{C}^{n\ell}(\frakn_P, T^{k-n\ell}(W_{\C})) \otimes T^{j+\ell}(U)[-\tfrac{p-q}{2}] \otimes \mathcal{W}_{n,W}\right)
\end{equation}
induced  by 
\begin{equation}\label{abc}
u_1^{j}  \otimes x \otimes f \otimes \omega \otimes w \mapsto
(u_1^{j} \otimes s_B^*(x_1))
\otimes x_2 \otimes f \otimes \omega \otimes \tau_{n\ell}(w_1) \otimes w_2
\end{equation}
gives a cocycle for the nilpotent $\mathfrak{n}_P$-complex. Going through the proof of Proposition~\ref{nonvanish}(i), we see that the composition of the $\mathfrak{n}_P$-differential with \eqref{abc} factors when viewed as a map on $T^{k}(W_{\C})$ through representations $\mathbb{S}_C(W_{\C})$ with $C$ having at least $n+1$ rows. But now by Lemma~\ref{SchurWeylforC} such representations do not occur in $C_W^{j,r,k}$. 
\end{proof}

The reader easily checks from the definition that $\iota_P$
satisfies the following properties.

\begin{lemma}\label{propofiota}

\begin{enumerate}
\item $\iota_P$ is a $[T(U)\otimes \mathcal{W}_{n,W}\otimes \bigwedge \mathfrak{p}_{W}^*]^{
K' \times K_W}$-module homomorphism. That is,
\[
\iota_P(\varphi^W_{j',r',0}\cdot \varphi^W_{j,r,k})=
\varphi^W_{j',r',0}\cdot \iota_P(\varphi^W_{j,r,k})
\]
for $\varphi^W_{j',r',0} \in C_W^{j',r',0}$ and $\varphi^W_{j,r,k}
\in C_W^{j,r,k}$.

\item $\iota_P(\varphi^W_{j,r,k})$ is zero if $k < n \ell$.

\item Suppose $\varphi^W_{j,r,k} \in C^{j,r,k}_W$ with $k \geq n \ell$ and
$\varphi_{j',r',\ell'}^W \in C_W^{j',r',\ell'}$. Then
\[
\iota_P( \varphi^W_{j,r,k} \cdot \varphi^W_{j',r',\ell'})=
\iota_P(\varphi^W_{j,r,k}) \cdot \varphi^W_{j',r',\ell'}.
\]

\item Let $x \in T^{n\ell}(\C^n)^*$ and $w
\in T^{n\ell}(W_{\C})$. Then
\[
\iota_P(1_U \otimes x \otimes 1_{\mathcal{F}}  \otimes
1_{\mathfrak{p}_W^{\ast}}\otimes w) = x(\eps_B)( u_{\ell} \otimes
1_{\C^n} \otimes
 1_{\mathcal{F}} \otimes
1_{ \mathfrak{p}_W^*} \otimes \tau_{n\ell}(w)\otimes 1_{T(V_{\C})}).
\]

\end{enumerate}
\end{lemma}

\begin{proposition}\label{iotaP-Prop1}

Let $k = n\ell +\ell'$ as above. Let $\la$ be a dominant weight of
$\GL_n(\C)$, and we let $A$ be a standard filling of the associated
Young diagram $D(\la)$. We let $B|A$ be the associated filling for
the weight $\ell\varpi_n + \la$, see section~\ref{rep-review}. 

\begin{itemize}

\item[(i)]
Then the preimage of $A_P^{j+\ell,r+n\ell,A}$ under $\iota_P$ lies in 
$C_W^{j,r,B|A}$;i.e.,
\[
\iota_P^{-1} \left(A_P^{j+\ell,r+n\ell,A}\right) = C_W^{j,r,B|A}.
\]
Moreover, if $\iota_P(\varphi') = \varphi$ for $\varphi' \in C_W^{j,r,n\ell+\ell'}$ and $\varphi \in 
A_P^{j+\ell,r+n\ell,\ell'}$, then
\[
\pi_{A} (\varphi) = \iota_P(\pi_{B|A}(\varphi')). 
\]
Here $\pi_{B|A}$ is  the projection from $C_W^{j,r,n\ell+\ell'}$ to
$C_W^{j,r,B|A}$, see \eqref{SWC}, and $\pi_{A}$ is the one from $A_P^{j+\ell,r+n\ell,\ell'}$ to $A_P^{j+\ell,r+n\ell,A}$.

\item[(ii)]
Let $\varphi \in A_P^{j+\ell,r+n\ell,[A]}$ be a closed form such that $\iota_P(\varphi') = \varphi$ for some $\varphi' \in C_W^{j,r,B|A}$. Let $\pi_{[B|A]}$ be the projection from $C_W^{j,r,B|A}$ to
$C_W^{j,r,[B|A]}$. Then the cohomology class $[\varphi]$ satisfies
\[
[\varphi] = [\iota_P(\pi_{[B|A]}(\varphi'))].
\]
\end{itemize}

\end{proposition}

\begin{proof}
(i) We first observe that $\iota_P$ is invariant under $s_B$ in the
$T^{n\ell}(W)$-factor and also $s(B^{\ast})$-invariant in the $T^{n\ell}(\C^n)^{\ast}$-factor, that is,
\begin{align*}
\iota_P &= \iota_P \circ (1_U \otimes 1_{T^{n\ell}(\C^{n})^*}
\otimes 1_{T^{\ell'}(\C^n)^*} \otimes 1_{\mathcal{W}}
 \otimes 1_{\mathfrak{p}_W^{\ast}} \otimes s_B \otimes 1_{T^{\ell'}(W)}) \\
& =\iota_P \circ (1_U \otimes s(B^{\ast}) \otimes
1_{T^{\ell'}(\C^n)^{\ast}}\otimes 1_{\mathcal{W}}
 \otimes 1_{\mathfrak{p}_W^{\ast}} \otimes 1_{T^{n\ell}(W)} \otimes 1_{T^{\ell'}(W)}).
\end{align*}
Taking the $S_{\ell'}$-invariance into account, we see that
$\iota_P$ maps
\begin{equation}
\left[T^j(U) \otimes \mathbb{S}_{B}(\C^n)^{\ast} \otimes
\mathbb{S}_A(\C^n)^{\ast} \otimes \mathcal{W}_{n,W} \otimes
{\bigwedge} {^r} (\mathfrak{p}_W^*)_{\C} \otimes
\mathbb{S}_{B}(W_{\C}) \otimes \mathbb{S}_A(W_{\C})\right]^{K'
\times K_W}
\end{equation}
to $A^{j+2\ell,r+n\ell,A}$. But now 
\begin{lemma}\label{nicefact}

\begin{align}\label{firstappr}
\left[T^j(U) \otimes \mathbb{S}_{B}(\C^n)^{\ast} \otimes
\mathbb{S}_A(\C^n)^{\ast} \otimes \mathcal{W}_{n,W} \otimes
{\bigwedge} {^r} (\mathfrak{p}_W^*) \otimes \mathbb{S}_{B}(W_{\C})
\otimes \mathbb{S}_A(W_{\C})\right]& ^{K' \times K_W} \\ & =
C_W^{j,r,B|A}. \notag
\end{align}

\end{lemma}

\begin{proof}
In \eqref{firstappr}, we first observe $\mathbb{S}_{B}(\C^n)^{\ast}
\otimes \mathbb{S}_A(\C^n)^{\ast} = \mathbb{S}_{B|A}(\C^n)^{\ast}$
as subspaces of $T^{n\ell+\ell'}(\C^n)$, see
Corollary~\ref{repsA-Bcoincide}. But then by Schur-Weyl theory, see
Lemma~\ref{SchurWeylforC} or Remark~\ref{SchurWeylforC1}, we can now
replace $\mathbb{S}_{B}(W_{\C}) \otimes \mathbb{S}_A(W_{\C})$ with
its subspace $\mathbb{S}_{B|A}(W_{\C})$ in \eqref{firstappr}, that
is, the left hand side in \eqref{firstappr} is equal to $C_W^{j,r,B|A}$.
\end{proof}

From this we obtain Proposition~\ref{iotaP-Prop1}(i). Proposition~\ref{iotaP-Prop1}(ii) follows from Proposition~\ref{nonvanish}(iii) and Lemma~\ref{Harderlemma}.
\end{proof}

\section{Special Schwartz forms}

Again $V$ will denote a real quadratic space of
dimension $m$ and signature $(p,q)$.

\subsection{Construction of the special Schwartz forms}

 We recall the construction in \cite{FMII} of the
special Schwartz forms $\varphi_{nq,\ell'}$, $\varphi_{nq,\la}$, and
$\varphi_{nq,[\la]}$, which define cocycles in $C_V^{\bullet,\ell'}$,
$C_V^{\bullet,\la}$, and $C_V^{\bullet,[\la]}$ respectively. It will
be more convenient to use the model $C_V^{\bullet}$ consisting of
homomorphisms on $T^{\ell'}(\C^n)$ (and its subspaces
$\mathbb{S}_{t(\la)}(\C^n)$), see \eqref{CHom} and
Remark~\ref{SchurWeylforC1}.  We will initially use the
Schr\"odinger model $\mathcal{S}(V^n)$.

In \cite{FMII}, we construct for $n \leq p$ a family of Schwartz
forms $\varphi_{nq,\ell'}$ on $V^n$ such that $\varphi_{nq,\ell'}
\in C_V^{q, nq,\ell'}$. So
\begin{align}\label{varphi}
\varphi_{nq,\ell'} & \in  \left[\Hom\left(T^{\ell'}(\C^n), T^q(U)
\otimes \calS(V^n) \otimes \mathcal
{A}^{nq}(D) \otimes  T^{\ell'}(V_{\C})\right)\right]^{K' \times G \times S_{\ell'}} \\
& \notag \simeq \left[  \Hom\left(T^{\ell'}(\C^n), T^q(U) \otimes
\calS(V^n) \otimes \wwedge{nq} (\mathfrak{p}^{\ast}_{\C}) \otimes
T^{\ell'}(V_{\C})\right) \right]^{K' \times K \times S_{\ell'}}.
\end{align}
These Schwartz forms are generalizations of the Schwartz forms
considered by Kudla and Millson \cite{KMI,KMII,KM90}. Under the
isomorphism in \eqref{varphi}, the standard Gaussian $\varphi_0(\x)
= 1 \otimes  e^{-\pi tr(\x,\x)_{z_0}} \in \left[ T^0(U) \otimes
\mathcal{S}(V^n)\right]^{K' \times K}$ corresponds to
\[
\varphi_0(\x,z) = 1 \otimes  e^{-\pi tr(\x,\x)_z} \in \left[T^0(U)
\otimes \mathcal{S}(V^n) \otimes C^{\infty}(D)\right]^{K' \times G}.
\]

\begin{definition}\label{calD-def}
Let $n \leq p$. The form $\varphi_{nq,0}$ with trivial coefficients is given
by applying the operator
\begin{gather*}
\calD = \frac{1}{2^{nq/2}} \, A\left(u_1^q\right) \otimes
\prod_{i=1}^n \prod_{\mu = p+1}^{p+q} \left[ \sum_{\alpha =1}^{p}
\left(  x_{\alpha i} - \frac{1}{2\pi} \frac{\partial}{\partial
x_{\alpha i }}  \right)   \otimes A(\omega_{\alpha\mu}) \right]
\end{gather*}
to  $\varphi_0$:
\begin{equation*}
 \varphi_{nq,0} = \calD ( \varphi_0) \in C_V^{q,nq,0}=
 \left[ T^q(U) \otimes \calS (V^n) \otimes \wwedge{nq}(\mathfrak{p}^{\ast}_{\C}) \right ]^{K' \times K}.
\end{equation*}
Here as before $A(\cdot)$ denotes left multiplication and $u_1$ is
the generator of $U = \wwedge{n} (\C^n)^{\ast}$. Furthermore,
Theorem 3.1 of \cite{KMI} implies that $\varphi_{nq,0}$ is indeed
$K'$-invariant.
\end{definition}

For $T(V_{\C})$, we define another
$K$-invariant differential operator $\calD'_i$ which acts on
\begin{equation}
\calS(V^n) \otimes \wwedge{\bullet}(\mathfrak{p}^{\ast}_{\C})
\otimes T(V_{\C})
\end{equation}
by
\begin{equation}
 \calD'_i = \frac12 \sum_{\alpha=1}^p \left( x_{\alpha i} -
\frac{1}{2\pi} \frac{\partial}{\partial x_{\alpha i}} \right)
\otimes 1 \otimes A(e_{\alpha}).
\end{equation}
Let $I=(i_1,\dots,i_{\ell'}) \in \{1,\dots,n\}^{\ell'}$ be a
multi-index of length $\ell'$ and write
\begin{equation}
\eps_I = \eps_{i_1} \otimes \cdots \otimes \eps_{i_{\ell'}}
\end{equation}
for the corresponding standard basis element of  $T^{\ell'} (\C^n)$.
Then for $\eps_I \in T^{\ell'} (\C^n)$, we define an operator by
\begin{equation}\label{varphi0ell}
\calT_{\ell'}(\eps_I) = \calD'_{i_1} \circ \cdots \circ
\calD'_{i_{\ell'}}
\end{equation}
extending $\calT_{\ell'}$ linearly to $T^{\ell'}(\C^n)$.

\begin{definition}\label{calT-def}
Define
\begin{align*}
 \varphi_{nq,\ell'} \in  C_V^{q, nq,\ell'} =
\Hom_{\C}\left( T^{\ell'}(\C^n), T^q(U) \otimes \mathcal{S}(V^n)
\otimes \wwedge{nq}(\mathfrak{p}^{\ast}_{\C}) \otimes
T^{\ell'}(V_{\C}) \right)^{K' \otimes K \otimes S_{\ell'}}
\end{align*}
by
\[
\varphi_{nq,\ell'}(w) = \calT_{\ell'}(w) \varphi_{nq,0}
\]
for $w \in T^{\ell'} (\C^n)$. We put $ \varphi_{nq,\ell'} =0$ for
$\ell' <0$. Here the $S_{\ell'}$-invariance of $\varphi_{nq,\ell'}$
is shown in Proposition~5.2 in \cite{FMII}, while the
$K'$-invariance is Theorem~5.6 in \cite{FMII}.

\end{definition}

Using the projections $\pi_{t(\la)}$ and $\pi_{[t(\la)]}$, see
\eqref{SWC} and \eqref{SWC1}, we can therefore make the following
definitions.

\begin{definition}
For any standard filling $t(\la)$ of $D(\la)$, we define
\begin{align*}
\varphi_{nq,t(\la)} &  = \pi_{t(\la)} \varphi_{nq,\ell'} \in C_V^{q, nq,\la}, \\
\varphi_{nq,[t(\la)]} & = \pi_{[t(\la)]} \varphi_{nq,\ell'}\in
C_V^{q, nq,[\la]}.
\end{align*}
We write $\varphi_{nq,\la}$ and $\varphi_{nq,[\la]}$, if we do not
want to specify the standard filling.
\end{definition}

\begin{proposition}[Theorem~5.7 \cite{FMII}]
The form $\varphi_{nq,\ell'}$ is closed. That is, for $w \in
T^{\ell'}(\C^n)$ and $\x \in V^n$, the differential form
\[
\varphi_{nq,\ell'}(w)(\x) \in
\left[{A}^{nq}\left(D;T^{\ell'}(V_{\C}) \right)\right]^G
\]
is closed.
\end{proposition}

\subsection{Explicit formulas}

We give more explicit formulas for $\varphi_{nq,\ell'}$ in the
various models of the Weil representation.

\subsubsection{Schr\"odinger model}

We introduce multi-indices $\underline{\alpha_i} =
(\alpha_{i1},\cdots,\alpha_{iq})$ of length $q$ (typically) with $1
\leq i \leq n$ and $\underline{\beta} =
(\beta_1,\cdots,\beta_{\ell'})$ of length $\ell'$ (typically) with
values in $\{1,\dots,p\}$ (typically). Note that we suppressed their
length from the notation. We also write $\mathbf{\alpha} =
(\alpha_{ij})$ for the $n \times q$ matrix of indices. With $I$ as
above, we then define
\begin{gather}\label{indices}
\omega_{\underline{\alpha_i}} = \omega_{\alpha_{i1} p+1} \wedge
\cdots \wedge \omega_{\alpha_{iq} p+q} \\ \notag \omega_{\alpha} =
\omega_{\underline{\alpha_1}} \wedge \cdots \wedge
\omega_{\underline{\alpha_n}} \\ \notag \calH_{\underline{\alpha_i}}
= \calH_{\alpha_{i1}i} \circ \cdots \circ\calH_{\alpha_{iq}i}, \\
\notag \calH_{\alpha} = \calH_{\underline{\alpha_1}} \circ \cdots
\circ \calH_{\underline{\alpha_n}} \\
\notag \calH_{\underline{\beta},I} =\calH_{\beta_1 i_1} \circ \cdots
\circ \calH_{\beta_{\ell'} i_{\ell'}} \\ \notag
 e_{\underline{\beta}} = e_{\beta_1} \otimes \cdots \otimes
e_{\beta_{\ell'}}
\end{gather}
 Let $ 1 \leq \gamma \leq p$ and $1 \leq j \leq n$. For $I$, $\alpha$, and
 $\underline{\beta}$ fixed, let
\begin{equation}
\delta_{\gamma j} = \# \{k; \,  \alpha_{kj} = \gamma \} + \#\{k; \,
(\beta_k,i_k) = (\gamma,j)\}.
\end{equation}
This defines a $p \times n$ matrix
$\Delta_{\alpha,\underline{\beta},I}
=\Delta_{\alpha,\underline{\beta},I;+}$ and Schwartz functions
$\varphi_{\Delta_{\alpha,\underline{\beta},I}}$ as in
Definition~\ref{varphidelta}.

\begin{lemma}\label{varphi-formula11}
The Schwartz form $\varphi_{nq,\ell'}(\eps_I)$ is given by
\[
\varphi_{nq,\ell'}(\eps_I) = \frac{1}{2^{nq/2+\ell'}} \sum_{\alpha,
\underline{\beta}} u_1^q \otimes
\varphi_{\Delta_{\alpha,\underline{\beta},I}} \otimes
\omega_{\alpha} \otimes e_{\underline{\beta}}.
\]
\end{lemma}

\begin{proof}
With the above notation we have
\begin{align}
\varphi_{nq,\ell'}(\eps_I) &= \frac{1}{2^{nq/2+\ell'}}  \sum_
{\substack{\underline{\alpha_1},\dots, \underline{\alpha_n} \\
\underline{\beta}}} u_1^q \otimes
 ( (\calH_{\underline{\alpha_1}} \circ \cdots \circ
\calH_{\underline{\alpha_n}} \circ
\calH_{\underline{\beta},I})\varphi_0)
 \otimes ( \omega_{\underline{\alpha_1}} \wedge \cdots
\wedge \omega_{\underline{\alpha_q}} ) \otimes
e_{\underline{\beta}} \\
&= \notag \frac{1}{2^{nq/2+\ell'}} \sum_{\alpha, \underline{\beta}}
u_1^q \otimes ( \calH_{\alpha} \circ
\calH_{\underline{\beta},I})\varphi_0 \otimes \omega_{\alpha}
\otimes e_{\underline{\beta}}.
\end{align}
But now we easily see
\begin{equation}
\left(\calH_{\alpha} \circ
\calH_{\underline{\beta},I}\right)\varphi_0(\x) = \prod_{\gamma
=1}^p \prod_{j=1}^n \widetilde{H}_{\delta_{\gamma,j}}(x_{\gamma j})
\varphi_0(\x),
\end{equation}
which gives the assertion.
\end{proof}

\subsubsection{Mixed model}

We now describe the Schwartz form $\varphi_{nq,\ell'}$ in the mixed
model. We describe this in terms of the individual components
$\varphi_{\Delta_{\alpha,\underline{\beta},I}}$ described in the
Schr\"odinger model. From Lemma~\ref{mixedformula10},
Lemma~\ref{keylemma2}, and Proposition~\ref{thetavanishing} we see

\begin{lemma}\label{mixedvarphi}
\[
\widehat{\varphi^V_{ \Delta_{\alpha,\underline{\beta},I}}}
\left(\begin{smallmatrix}  \xi \\ \x_W \\ u'
\end{smallmatrix} \right)
= \varphi^W_{\Delta'_{\alpha,\underline{\beta},I}}(\x_W)
 \widehat{\varphi^E_{\Delta_{\alpha,\underline{\beta},I}^{\prime \prime}}}
(\xi,u').
\]
Note that  $\varphi^W_{\Delta'_{\alpha,\underline{\beta},I}}$ only
depends on the indices $\alpha_{ij}, \beta_j$ such that $
\alpha_{ij}, \beta_j\geq \ell+1$, while
$\widehat{\varphi^E_{\Delta_{\alpha,\underline{\beta},I}^{\prime
\prime}}}$ only depends on the indices $\alpha_{ij}, \beta_j$ such
that $\alpha_{ij}, \beta_j\leq \ell$. In particular, if \emph{all}
$\alpha_{ij}, \beta_j \geq \ell+1$, then
\[
\widehat{\varphi^V_{\Delta_{\alpha,\underline{\beta},I}}}
\left(\begin{smallmatrix}  \xi \\ \x_W \\ 0
\end{smallmatrix} \right)=
\varphi^W_{\Delta'_{\alpha,\underline{\beta},I}}(\x_W)
\varphi_0^E(\xi,0).
\]
On the other hand, if \emph{one} of the $\alpha_{ij}, \beta_j$ is
less or equal to $\ell$, then
\[
\widehat{\varphi}^E_{\Delta_{\alpha,\underline{\beta},I}^{\prime\prime}}(0,0)
= \widehat{\varphi^V_{\Delta_{\alpha,\underline{\beta},I}}}
\left(\begin{smallmatrix}  0 \\ \x_W \\ 0
\end{smallmatrix} \right) = 0.
\]
\end{lemma}

\subsubsection{Fock model}

In the Fock model, $\varphi_{nq,\ell'}$ looks particularly
simple. We have

\begin{lemma}\label{Fock-varphi}
\[
\varphi_{nq,\ell'}(\eps_I) = \frac{1}{2^{nq/2+\ell'}} \left(
\frac{1}{2\pi i}\right)^{nq+\ell'} \sum_
{\substack{\underline{\alpha_1},\dots, \underline{\alpha_n} \\
\underline{\beta}}} u_1^q \otimes z_{\underline{\alpha_1},1} \cdots
z_{\underline{\alpha_n},n} \cdot z_{\underline{\beta},I} \otimes
\left( \omega_{\underline{\alpha_1}} \wedge \cdots \wedge
\omega_{\underline{\alpha_q}} \right) \otimes e_{\underline{\beta}}.
\]
Here we use the notational conventions in \eqref{indices} and in
addition
\begin{align}\label{indices2}
z_{\underline{\alpha_j},j}  =  z_{\alpha_{j1}j} \cdots
z_{\alpha_{jq}j},  \qquad \qquad z_{\underline{\beta},I} =
z_{\beta_1 i_1} \cdots z_{\beta_{\ell'} i_{\ell'}}.
\end{align}
\end{lemma}

\subsection{The forms $\varphi_{0,k}$}

We now define another class of special forms. We will only do this
in the Fock model.

\begin{definition}

We define $\varphi_{0,k} \in \Hom\left(T^{k}(\C^n); T^0(U) \otimes
\mathcal{F}_{n,V} \otimes T^{k}(V_{\C})\right)$ by
\begin{equation}
\varphi_{0,k}(\eps_I) =  \frac{1}{2^{k}} \left( \frac{1}{2\pi
i}\right)^{k} \sum_ { \underline{\beta}} 1 \otimes
z_{\underline{\beta},I} \otimes e_{\underline{\beta}}.
\end{equation}

\end{definition}

\begin{remark} \label{varphikrem}
The element $\varphi_{0,k}$ is the image of the operator
$\mathcal{T}_{k}$ (see \eqref{varphi0ell}) applied to the Gaussian
$\varphi_0$ under the intertwiner from the Schr\"odinger to the Fock
model. Also note that $\varphi_{0,k}$ is \emph{not} closed, hence
they do not define cocycles.
\end{remark}

We also leave the proof of the following lemma to the reader. It
follows (in large part) from Remark~\ref{varphikrem} and the
corresponding properties of $\varphi_{nq,\ell'}$.

\begin{lemma}
\[
\varphi_{0,k} \in [T^{0}(U) \otimes T^{k}(\C^n)^*\otimes
\mathcal{F}_{n,V}  \otimes T^{k}(V_{\C})]^{K' \times K \times S_k},
\]
i.e.,
\[
\varphi_{0,k} \in C_V^{0,0,k}.
\]
\end{lemma}

From Lemma~\ref{Fock-varphi}, we immediately see

\begin{lemma}\label{firstproductrule}
\[
\varphi_{nq,\ell'} = \varphi_{nq,0} \cdot \varphi_{0,\ell'}
\]
and
\[
\varphi_{0,k_1} \cdot \varphi_{0,k_2} = \varphi_{0,k_1 + k_2},
\]
where the multiplication is the one in $C_V^{\bullet}$.
\end{lemma}

\begin{remark}
This kind of product decomposition for $\varphi_{nq,\ell'}$ and
$\varphi_{0,k}$ in Lemma~\ref{firstproductrule} only holds in the
Fock model. In the Schr\"odinger model this only makes sense in
terms of the operators $\calD$ and $\calT_{\ell'}$ of
Definition~\ref{calD-def} and Definition~\ref{calT-def}
respectively.
\end{remark}

We apply the projection $\pi_{t(\la)}$, see \eqref{SWC}, to define
$\varphi_{0,t(\la)}$:

\begin{definition}
\[
\varphi_{0,t(\la)} := \pi_{t(\la)} \varphi_{0,k} \in
C_V^{0,0,t(\la)}.
\]
\end{definition}

The following product formula will be important later.

\begin{proposition}\label{productrule}
Let $A=t(\la)$ be a filling of the Young diagram associated to $\la$
and let $B=B_{n,\ell}$ be the filling of the $n\times \ell$
rectangular Young diagram introduced in section 3. Then
\[
\varphi^W_{0,B}\cdot \varphi^W_{0,A} =\varphi^W_{0,B|A}.
\]
\end{proposition}
The proposition will follow from the next two lemmas.

\begin{lemma}\label{repscoincide}
Both $\varphi^W_{0,B} \cdot \varphi^W_{0,A}$ and $\varphi^W_{0,B|A}$
are elements of
\[
C_W^{0,B|A,0} = \left[T^{0}(U)\otimes \mathbb{S}_{B|A}(\C^n)^*
\otimes \mathcal{F}_{n,W} \otimes
\mathbb{S}_{B|A}(W_{\C})\right]^{K' \times K_W}.
\]

\end{lemma}

\begin{proof}
Since $\mathbb{S}_{B}(\C^n)^{\ast} \otimes \mathbb{S}_A(\C^n)^{\ast}
= \mathbb{S}_{B|A}(\C^n)^{\ast}$ as subspaces of
$T^{n\ell+\ell'}(\C^n)$, see Corollary~\ref{repsA-Bcoincide}, the
claim follows in the same way as Lemma~\ref{nicefact}.
\end{proof}

\begin{lemma}\label{thephiscoincide}
\[
(\varphi^W_{0,B}\cdot \varphi^W_{0,A})(s_B\eps_B \otimes
s_A\eps_A) = \varphi^W_{0,B|A}(s_B\eps_B \otimes s_A\eps_A).
\]
\end{lemma}
\begin{proof}
This is a little calculation using Lemma~\ref{John} and
Lemma~\ref{firstproductrule}. Indeed, we have
\begin{align*}
\left(\varphi^W_{0,B}\cdot \varphi^W_{0,A}\right) & \left(s_B\eps_B
 \otimes s_A\eps_A\right) = \left(\varphi^W_{0,n
\ell}\cdot\varphi^W_{0,\ell'}\right )\left(s_B\eps_B \otimes s_A
\eps_A\right)
\\
& = \varphi^W_{0,n\ell + \ell'}\left(s_B\eps_B \otimes s_A
\eps_A\right)=c(A,B)
\varphi^W_{0,n\ell + \ell'}\left(s_{B|A}\eps_{B|A}\right) \\
&= c(A,B)\varphi^W_{0,B|A}\left(s_{B|A}\eps_{B|A}\right) =
\varphi^W_{0,B|A}\left(s_B\eps_{B}\otimes s_A \eps_A\right). \qedhere
\end{align*} 
\end{proof}
Now we can prove Proposition \ref{productrule}. By
Lemma~\ref{repscoincide} we see that $\varphi^W_{0,B} \cdot
\varphi^W_{0,A}$ and $\varphi^W_{0,B|A}$ are $U(n)$-equivariant
homomorphisms from $\mathbb{S}_{B|A}(\C^n)^{\ast}$ to $ T^{0}(U)
\otimes \mathcal{F}_{n,W} \otimes \mathbb{S}_{B|A}(W_{\C}).$ By
Lemma~\ref{thephiscoincide} they agree on the highest weight vector
(see Lemma~\ref{John}), hence coincide.

\section{Local Restriction}

We retain the notation from the previous sections. In this section,
we will give formulas for the restrictions $r_P^{\mathcal{W}}$ and
$r_P$ of $\varphi_{nq,\ell'}$. The main result will be then the
local restriction formula for $\varphi_{nq,[\la]}$.

\begin{proposition}\label{throwaway1}
We have
\[
(r^{\mathcal{W}}_{\underline{P}} \varphi^V_{nq,\ell'})(\eps_I) =
\frac1{2^{nq/2+\ell'}} \sum_{\alpha',\underline{\beta}'}
u_1^q \otimes \varphi^W_{\Delta'_{\alpha',\underline{\beta'},I}}
\otimes \omega_{\alpha'} \otimes e_{\underline{\beta}'}.
\]
Here $\eps_I = \eps_{i_1} \otimes \cdots \otimes \eps_{i_{\ell}} \in
T^{\ell}(\C^n)$, $\alpha'$ and $\underline{\beta}'$ are the same
indices as before with
\[
\ell+1 \leq \alpha'_{ij}, \beta'_j \leq p.
\]
Loosely speaking,
$r_P^{\mathcal{W}}\left(\varphi^V_{nq,\ell'}\right)$ is obtained
from $\varphi^V_{nq,\ell'}$ by "throwing away" all the indices less
or equal to $\ell$. In particular, if $n>p-\ell$, we have
\[
r^{\mathcal{W}}_{\underline{P}} \varphi^V_{nq,\ell'} =0.
\]

\end{proposition}

\begin{proof}
This follows from Lemma~\ref{varphi-formula11}, the formula for
$\varphi_{nq,\ell}$ in the Schr\"odinger model, and from
Lemma~\ref{mixedvarphi}. For the last statement, we observe that
$\omega_{\alpha'}$ is in the $nq$-exterior power of a
$(p-\ell)q$-dimensional space.
\end{proof}

The local restriction looks particularly simple in the Fock model.

\begin{proposition}\label{throw-away-lemma}
Let $\underline{\alpha'_j}$ and $\underline{\beta'}$ be as before
in Proposition~\ref{throwaway1}. Then 
\begin{align*}
& r_P^{\mathcal{W}}\left(\varphi^V_{nq,\ell'}(\eps_I)\right) \\ &
\quad = \frac{1}{2^{nq/2+\ell'}} \left( \frac{1}{2\pi
i}\right)^{nq+\ell'} \sum_
{\substack{\underline{\alpha'_1},\dots, \underline{\alpha'_n} \\
\underline{\beta'}}} u_1^q \otimes z_{\underline{\alpha'_1}} \cdots
z_{\underline{\alpha'_n}} \cdot z_{\underline{\beta'},I} \otimes
\left( \omega_{\underline{\alpha'_1}} \wedge \cdots \wedge
\omega_{\underline{\alpha'_q}} \right) \otimes
e_{\underline{\beta'}}.
\end{align*}
\end{proposition}

\begin{proof}
This follows immediately either from Proposition~\ref{throwaway1}
and applying the intertwiner to the Fock model or also from
Proposition~\ref{Fock-thetavanishing} and Lemma~\ref{Fock-varphi}.
\end{proof}

\begin{proposition}
For the restriction of $\varphi^V_{nq,\ell'}$, we have
\[
r_{P} \varphi^V_{nq,\ell'} = \left( 1_U \otimes 1_{\C^n} \otimes
r^{\mathcal{W}}_{P} \otimes \sigma^{\ast} \otimes 1_V \right)
\varphi^V_{nq,\ell'}.
\]
Analogous statements hold for $\varphi^V_{nq,\la}$
and $\varphi^V_{nq,[\la]}$. 
\end{proposition}

\begin{proof}
By Definition~\ref{def:localres}, the restriction $r_P:C_V^{\bullet}
\to A_P^{\bullet}$ is given by $1_{U} \otimes 1_{\C^n} \otimes
r_P^{\mathcal{W}} \otimes ( \iota^{\ast} \circ \sigma^{\ast} )
\otimes 1_V$. Then the theorem follows from
Proposition~\ref{throw-away-lemma} and Lemma~\ref{forms-formulas},
in particular \eqref{throw-away-key}: The components of
$\sigma^{\ast} \varphi^V_{nq,\ell'}$ involving $\mathfrak{a}^{\ast}
$ already become annihilated under $r^{\mathcal{W}}_P$, so that
$\iota^{\ast}$ acts trivially on $\sigma^{\ast}  r^{\mathcal{W}}_{P}
\varphi^V_{nq,\ell'}$.
\end{proof}

We define
\begin{equation}\label{varphiP-def}
\varphi_{P,n\ell} =  \frac{1}{2^{n\ell  }}  \left( \frac{1}{2\pi i}\right)^{n\ell}
\sum_{\underline{\g_1},\dots, \underline{\g_n}} u_1^{\ell} \otimes
z_{\underline{\g_1},1} \cdots z_{\underline{\g_n},n} \otimes (
\nu_{\underline{\g_1}} \wedge \cdots \wedge \nu_{\underline{\g_n}}).
\end{equation}
Here $\underline{\g_j} = (\g_{j m-\ell+1},\dots,\g_{j m }) $ is a
multi-index of length $\ell$ such that $\ell+1 \leq \g_{ji} \leq p$,
and $z_{\underline{\g_j},j}$ as in \eqref{indices2}. Furthermore, we
have set
\begin{equation}\label{nu-convention}
\nu_{\underline{\g_j}}=  \nu_{\g_{j m-\ell+1} \ell} \wedge \cdots
\wedge \nu_{\g_{j m} 1} \in \wwedge{\ell} (\mathfrak{n}^{\ast}_W).
\end{equation}
We have

\begin{lemma}\label{imageunderiota}
\[
\iota_P (\varphi^W_{0,B}) = \iota_P (\varphi^W_{0,n\ell})=(-1)^{n\ell (\frac{(q-\ell)(n-1)}{2}+1)} 
\varphi_{P,n\ell}.
\]

\end{lemma}

\begin{proof}
First note that by Proposition~\ref{iotaP-Prop1} we have $\iota_P
(\varphi^W_{0,B}) =\iota_P (\varphi^W_{0,n\ell})$. We let
$\underline{\beta_1}, \dots, \underline{\beta_n}$ be $n$ indices of
length $\ell$ with $\ell+1 \leq \beta_{ji} \leq p$. For the
corresponding elements $e_{\underline{\beta_j}} \in T^{\ell}(W)$, we
easily see
\begin{align}
\sum_{\underline{\beta_1},\dots,\underline{\beta_n}}
(z_{\underline{\beta_1}1} \cdots z_{\underline{\beta_n}n}) \otimes
\tau_{n\ell}(e_{\underline{\beta_1}} \otimes \cdots \otimes
e_{\underline{\beta_n}}) = 
\sum_{\underline{\beta_1},\dots,\underline{\beta_n}}
(z_{\underline{\beta_1}1} \cdots z_{\underline{\beta_n}n}) \otimes
(\nu_{\underline{\beta_1}} \wedge \cdots \wedge
\nu_{\underline{\beta_n}})
\end{align}
with $\nu_{\underline{\beta_j}}$ as in \eqref{nu-convention}. With
that, we conclude
\begin{align}
\iota_P (\varphi^W_{0,B}) & = (-1)^{n\ell (\frac{(q-\ell)(n-1)}{2}+1)} 
\frac{1}{2^{n\ell}} \left( \frac{1}{2\pi i}\right)^{n\ell} \sum_
{\underline{\beta_1},\dots,\underline{\beta_n}} u_1^{\ell} \otimes
(z_{\underline{\beta_1}1}\cdots z_{\underline{\beta_n}n})\otimes
(\nu_{\underline{\beta_1}}
\wedge \cdots \wedge \nu_{\underline{\beta_n}})  \\
&= (-1)^{n\ell (\frac{(q-\ell)(n-1)}{2}+1)} \varphi_{P,n\ell} \notag
\end{align}
by \eqref{varphiP-def}. 
\end{proof}

We are now ready for the main result of this section, the
local restriction formula for $\varphi_{nq,[\la]}$.

\begin{theorem}\label{localrestrictiontheorem}

Let $A$ be a standard filling of Young diagram with $\ell'$ boxes
and let $B_{n,\ell}$ be the standard tableau associated to the $n$
by $\ell$ rectangle as in section~\ref{rep-review}. Then
\begin{align*}
r_P(\varphi^V_{nq,\ell'}) &= \iota_P(\varphi^{W}_{n(q-\ell),n\ell+\ell'}), \\
r_P(\varphi^V_{nq,A}) &= \iota_P(\varphi^{W}_{n(q-\ell),B|A}).
\end{align*}
Furthermore, for the form $\varphi^V_{nq,[A]}$ with harmonic
coefficients, we have
\[
[r_P(\varphi^V_{nq,[A]})] = 
[\iota_P(\varphi^{W}_{n(q-\ell),[B|A]})].
\]
\end{theorem}

\begin{proof}
We first note  
\[
r_{P} \varphi^V_{nq,\ell'} = (-1)^{n\ell (\frac{(q-\ell)(n-1)}{2}+1)}
\varphi^W_{n(q-\ell),0} \cdot \varphi_{P,n\ell} \cdot
\varphi^W_{0,\ell'}.
\]
Here we view $\varphi^W_{n(q-\ell),0} \in
A_P^{q-\ell,n(q-\ell),0}$ and $\varphi^W_{0,\ell'} \in
A_P^{0,0,\ell'}$ in the natural fashion. The analogous statements
hold for $\varphi^V_{nq,A}$ and $\varphi^V_{nq,[A]}$. Indeed, this follows immediately from Proposition~\ref{throw-away-lemma} and
\begin{equation}\label{sigma-formula}
\sigma^{\ast} \omega_{\underline{\alpha'_j}} = (-1)^{\ell} \frac1{2^{\ell/2}} \omega_{\alpha'_{j1}
p+1} \wedge \cdots \wedge \omega_{\alpha'_{j q-\ell} m-\ell} \wedge
\nu_{\alpha'_{j q-\ell+1} \ell } \wedge \cdots \wedge
\nu_{\alpha'_{j q} 1 },
\end{equation}
which follows from Lemma~\ref{forms-formulas}. The sign arises from
'sorting' $\sigma^{\ast} \left( \omega_{\underline{\alpha'_1}}
\wedge \cdots \wedge \omega_{\underline{\alpha'_q}} \right)$
according to \eqref{sigma-formula} into
$\omega_{\alpha'_{\bullet}}$'s (which lie in $\mathfrak{p}_W$) and
$\nu_{\alpha'_{\bullet}}$'s (which lie in $\mathfrak{n}^{\ast}_W$). From this and Lemma~\ref{imageunderiota} we conclude 
\[
r_P(\varphi^{V}_{nq,\ell'})= \iota_P(\varphi^{W}_{n(q-\ell),0} \cdot
\varphi^W_{0,n\ell} \cdot \varphi^W_{0,\ell'}) = \iota_P(\varphi^{W}_{n(q-\ell),n\ell+\ell'}).
\]
By $S_{\ell'}$-equivariance of $\iota_P$ we also obtain
\begin{align*}
r_P(\varphi^{V}_{nq,A}) = \iota_P(\varphi^{W}_{n(q-\ell),0} \cdot
\varphi^W_{0,B} \cdot \varphi^W_{0,A}) = \varphi^W_{n(q-\ell),B|A}
\end{align*}
since $\varphi^W_{0,B}\cdot \varphi^W_{0,A} =\varphi^W_{0,B|A}$ (see
Proposition~\ref{productrule}) and $\varphi^W_{n(q-\ell),B|A} =
\varphi_{n(q-\ell),0} \cdot \varphi^W_{0,B|A}$ (see Lemma~\ref{firstproductrule}).The cohomology statement now follows from Proposition~\ref{iotaP-Prop1}(ii).
 \end{proof}

\begin{corollary}
We have $[r_P(\varphi^V_{nq,[\lambda]})]=0$ for $n>\min\left(p,\left[\tfrac{m}{2}\right]\right)-\ell$ (if $\ell \geq 2$) and $n> p-1$ or $n>m-2-i(\la)$ (if $\ell=1$). 
\end{corollary}

\begin{proof}
The Schur functor $\mathbb{S}_{[B|A]}(W_{\C})$ vanishes in this range.
\end{proof}

On the other hand, we have

\begin{corollary}
Let $P$ be a (real) parabolic subgroup as above such that the associated space $W$ is positive definite. Assume 
\begin{equation*}
 i(\la) \leq n \leq 
\begin{cases}  
\left[ \frac{p-q}{2} \right] & \; \text{if $q\geq2$} \\ 
  p-1 -i(\la) & \text{ if $q=1$}. 
\end{cases}
\end{equation*}
Then
\[
[r_P(\varphi^V_{nq,[\lambda]})] \ne 0.
\]
\end{corollary}

\section{Global complexes, theta series, and the global restriction}\label{globalressection}

In this section, we return to the global situation and assume that $V,W,E$ etc. are $\Q$-vector spaces. Furthermore, $\underline{P}$ is a standard $\Q$-parabolic subgroup and $P=\underline{P}_0(\R)$ for its real points etc. All the 'local' notions (over $\R$) of the previous
sections carry over naturally to this situation, and we make use of
the already established notation.

Let $L \subset V$ be an even $\Z$-lattice of full rank, i.e., $(x,x)
\in 2 \Z$ for $x \in L$. In particular, $L \subset L^{\#}$, the dual
lattice. We fix $h \in (L^{\#})^n$ once and for all and pick a congruence subgroup $\G
\subset \underline{G}(\Z)$ of finite index which stabilizes $\calL := \calL_V= L^n +h$. 
The associated locally symmetric space $X= \G \back D$ is non-compact (since the Witt index of $V$ is positive) 
but has finite volume.

\subsection{Global complexes and theta series}

\subsubsection{Global complexes}

We first define "global" versions of the "local" complexes
$C^{\bullet}$ of forms on $X = \G \back D$, $A_P^{\bullet}$ of forms
on $e'(\underline{P}) = \G_P \back e(\underline{P})$. We set 
\begin{equation}
C^{\infty}(\G', j,\la) :=
C^{\infty}\left( \G' \back G'; T^{j}(U) \otimes
\mathbb{S}_{\la}(\C^n)^{\ast} \right)^{K'} 
\end{equation}
for $\G'$ an (appropriate) arithmetic subgroup of $\Sp(n,\Z)$. Note that we can identify $C^{\infty}(\G',j,\la)$ in the usual way with the space of vector-valued $C^{\infty}$-functions on the Siegel upper half space of genus $n$, transforming like a Siegel modular
form of type $\det^{j/2} \otimes \mathbb{S}_{\la}(\C^n)$. Furthermore, we denote by $Mod(\G',j,\la)$ the space of holomorphic Siegel modular forms of this type. 
We let
\begin{align}
\widetilde{C}^{j,r,\la}_V &=
C^{\infty}(\G', j,\la) {\otimes}
\left[\mathcal{A}^{r}(D) \otimes \mathbb{S}_{\la} (V_{\C})\right]^{\G},  \\
& \simeq 
C^{\infty}(\G', j,\la) {\otimes}
\left[\wwedge{r}(\mathfrak{p}_{\C}^*) \otimes \mathbb{S}_{\la} (V_{\C}) \otimes C^{\infty}(\G \back G)\right]^K \notag
\end{align}
and
\begin{align}
\widetilde{A}_P^{j,r,\la} & = C^{\infty}(\G', j,\la)
{\otimes}
\left[\mathcal{A}^{r}(e'(P) \otimes \mathbb{S}_{\la} (V_{\C})\right]^{\G_P} \\
& \simeq
C^{\infty}(\G', j,\la){\otimes}
\left[\wwedge{r}(\mathfrak{n} \oplus \mathfrak{p}_M)_{\C}^* \otimes \mathbb{S}_{\la} (V_{\C}) \otimes C^{\infty}(\G_P \back P) \right]^{K_P} \notag
\end{align}
We then define $\widetilde{C}^{j,r,[\la]}_V$ and $\widetilde{A}_P^{j,r,[\la]}$ as in the local case by harmonic
projection onto $\mathbb{S}_{[\la]} (V_{\C})$. 
The local map $\iota_P$ induces a global intertwining map of complexes
\begin{equation}
\tilde{\iota}_P:
\widetilde{C}_W^{j-\ell,r,\ell\varpi_n+\la} \to \widetilde{A}_P^{j,n\ell+r,\la}.
\end{equation}
by lifting functions on $\G_W \back \SO_0(W_{\R})$ to $\G_M \back M$.  This induces a map on cohomology 
\begin{align}
& \tilde{\iota}_P:  C^{\infty}(\G', j,\la) \otimes 
H^{n(q-\ell)}(X_W,\mathbb{S}_{[\ell\varpi_n+\la]}(W_{\C}))  \notag \\ 
&\hookrightarrow C^{\infty}(\G', j,\la) \otimes
H^{n(q-\ell)}(X_M, H^{n\ell}(\mathfrak{n},\mathbb{S}_{[\la]}(V_{\C}))) \\
&\hookrightarrow C^{\infty}(\G', j,\la) \otimes
H^{nq}(e'(P),\mathbb{S}_{[\la]}(V_{\C})). \notag
\end{align}

We also introduce
\begin{equation}
\overline{C}^{j,r,\la}_V =
C^{\infty}\left( \G' \back G'; T^{j}(U) \otimes
\mathbb{S}_{\la}(\C^n)^{\ast} \right)^{K'} {\otimes}
\mathcal{A}^{r}(\overline{X}, \mathbb{S}_{\la} (\calV_{\C})),
\end{equation}
the complex associated to the differential forms on the compactification $\overline{X}$ with values in $\mathbb{S}_{\la} (\calV_{\C})$, the local system associated to
$\mathbb{S}_{\la} (V_{\C})$. We then have a restriction map
\begin{equation}
\tilde{r}_P: \overline{C}_V^{\bullet} \to
\widetilde{A}_P^{\bullet}
\end{equation}
induced by the inclusion $e'(P) \hookrightarrow \overline{X}$. 

\subsubsection{Theta series}

Using the Schr\"odinger model $\calS(V_{\R}^n)$ of the Weil
representation, we now introduce for $\varphi \in C_V^{j,r,\la}$,
its theta series $\theta(\varphi)$ as follows. For $g' \in G'$, we then define for $z \in D$, the
theta series
\begin{equation}
\theta_{\calL_V}(g',z,\varphi) = \sum_{\x \in \calL_V} \omega(g')
\varphi(\x,z).
\end{equation}
We easily see that the series is $\G$-invariant as $\G$ stabilizes
$\calL_V$. Thus  $\theta_{\calL_V}$ descends to a closed differential $nq$-form on the
locally symmetric space $X = \G \back D$. More precisely, by the
standard theta machinery, we have
\begin{equation}
\theta_{\calL_V}(\varphi) \in \widetilde{C}_V^{j,r,\la}
\end{equation}
for some congruence subgroup $\G' \subseteq \Sp(n,\Z)$. Summarizing,
the theta distribution $\theta_{\calL_V}$ associated to $\calL$ gives
rise to a $G' \times G$ intertwining map of complexes
\begin{equation}
\theta_{\calL_V}: C_V^{\bullet} \longrightarrow
\widetilde{C}_V^{\bullet}.
\end{equation}

\begin{remark}
The main point of \cite{FMII} is that for the Schwartz forms $\varphi_{nq,[\la]}$ one has 
\[
[\theta_{\calL_V}(\varphi_{nq,[\la]})] \in Mod(\G',j,\la) \otimes H^{nq}(X, \mathbb{S}_{\la} (\calV_{\C})),
\]
and the Fourier coefficients are Poincar\'e dual classes of special cycles with nontrivial local coefficients.
\end{remark}

For a similar theta intertwiner for $A_P$, we note that
$A_P$ involves the Weil representation for $W=
E^{\perp}/E$. Recall (see Proposition~\ref{intertwiner} and
Definition~\ref{def:localresWeil}) that we can extend the action of
$\Orth(W_{\R})$ on $S(W_{\R}^n)$ to $P$ such that the Weil
representation intertwining map $r_P^{\calW}$ becomes an
$MN$-intertwiner. In particular, $N$ and $M_P'$ act trivially on
$S(W_{\R}^n)$. We let $\calL_W$ be a linear combination of delta functions of (cosets of) lattices in $W^n$, which is stabilized by $\G_P$, that is, by $\G_W$. Recall that we defined $\G_W$ as the image of $\G_P$ when acting on $E^{\perp}/E$. It contains $\G \cap \SO_0(W_{\R})$ as a finite subgroup of finite index.  Applying the theta distribution associated to $\calL_{W}$ we obtain an intertwining map
\begin{equation}
\theta_{\calL_W}: A_P^{\bullet} \to \widetilde{A}_P^{\bullet}.
\end{equation}
Furthermore, $\theta_{\calL_W}$ commutes with $\iota_P$:
\begin{equation}
 \theta_{\calL_W} \circ \iota_P = \tilde{\iota}_P \circ \theta_{\calL_W}.
\end{equation}
More general, we let 
\begin{equation}
A_P^{\bullet,\calL_W,\G_W} = \{ \varphi \in A_P^{\bullet}; \, \theta_{\calL_W}(\varphi) \;\text{is $\G_W$-invariant}\}. 
\end{equation} 
and obtain a map $\theta_{\calL_W}: A_P^{\bullet,\calL_W,\G_W} \to \widetilde{A}_P^{\bullet}$ as before. 

We will be interested in a particular $\calL_W$, which naturally
arises from $\calL_V$ as follows. Let $\pi_E: E^{\perp} \to E^{\perp}/E$ be the natural projection map. We then set
\begin{equation}
\widehat{\calL}_W := \pi_E(\calL_V \cap E^{\perp}). 
\end{equation}
For this definition, it is crucial to view $W = E^{\perp}/E$ as a
subquotient of $V$ and not as the subspace $E^{\perp} \cap (E')^{\perp}$ of $V$. Namely, $\widehat{\calL}_W$ is in general larger than $W \cap \calL^n$, which can be empty even when 
$\widehat{\calL}_W$ is not.

\begin{remark}
The notation of $\widehat{\calL}_W$ becomes
more transparent if one changes to the adelic setting. Adelically,
$\calL$ corresponds to the characteristic function $\chi_{\calL_V}$ of
the image of $\calL_V$ inside $V(\A_f)$, where $\A_f$ denotes the
finite adeles. Then in this setting, $\widehat{\calL}_W$ corresponds
to the partial Fourier transform of $\chi_{\calL_V}$ with respect to
$E(\A_f)$ when restricted to $W$. From this perspective, the
assignment $\calL \to \widehat{\calL}_W$ is the analogue at the
finite places of the local restriction map $r_{\underline{P}}$ at
the infinite place.
\end{remark}

\subsection{The global restriction } 
\subsubsection{Smooth forms on smooth manifolds with corners}
We begin with a short discussion of the definition of a smooth $\ell$-form on a smooth $n$-manifold with corners $M$.
For more on smooth manifolds with corners we refer the reader to the Appendix of \cite{BS} or \cite{Lee}, pp363-370. First, for any point $x \in M$ the tangent space $T_x(M)$ is a linear space of dimension $n$. A differential $\ell$-form $\omega$
will be a section of $\bigwedge^{\ell}(T^*(M))$.  To say when an $\ell$-form $\omega$ is smooth on $M$  it suffices to define  smooth $\ell$-forms on the local models $S^n_k = \R^k_{\geq 0} \times \R^{n-k}$.  

\begin{definition}
An $\ell$-form $\omega$ on $\R^k_{\geq 0} \times \R^{n-k}$ is smooth if there exists $\epsilon >0$ and a smooth form $\widetilde{\omega}$
on $\R^k_{> -\epsilon} \times \R^{n-k} \supset \R^k_{\geq 0} \times \R^{n-k}$ such that $\widetilde{\omega}$ restricts to
$\omega$.
\end{definition}

For our purposes we need only two classes of smooth forms.  Recall from the appendix of \cite{BS} that a point $x$ in a neighborhood $U$ that maps by a chart $\varphi$  to  the local model $S^n_k$ above with $\varphi(x) = 0$  is said to have index $k$.  The set of points of index greater than or equal to $k$ is denoted $M^{(k)}$.  The subset $M^{(0)}$ is said to be the interior of $M$, the set $M^{(1)}$ is said to be the boundary of $M$.
The first class of smooth $\ell$-forms on $M$ is obtained by extending by zero from $M^{(0)}$ to $M$  the smooth $\ell$-forms on $M^{(0)}$ whose coefficients relative to one and hence any system of coordinates vanish to
infinite order on $M^{(1)}$.   
 The second class of smooth $\ell$-forms on $M$ will consist of the {\it special} forms.  We define an $\ell$-form $\omega$ in a local model $S^n_k$ to be special if there exists an $\ell$-form $\overline{\omega}$ on $\R^{n-k}$ such that $\omega = p_2^* \overline{\omega}$, where $p_2: S^n_k \to \R^{n-k}$ is projection on the second factor.
We now claim that $\omega$ special implies that is smooth.  Indeed if we let $q_2: \R^k \times R^{n-k} \to \R^{n-k}$
then $\widetilde{\omega} := q_2^* \overline{\omega}$ provides the desired extension of $\omega$. 

\begin{remark}
This  definition of special forms for general smooth manifolds with corners in less restrictive than the definition in
\cite{GHMP}, Definition 13.2, p.169 for the case of Borel-Serre compactifications. In this latter definition the form $\overline{\omega}$ is required to have further
properties (e.g. local left $N_P$-invariance) that use the special features of the Borel-Serre compactification. 
\end{remark}

\subsubsection{The restriction formula}
We now prove
\begin{theorem}\label{resformula}
Assume $V$ is different from the $\Q$-split space for signature $(p,p)$. Then (see Remark \ref{weightedremark} below)  the theta series $\theta_{\calL_V}(\varphi_{nq,\ell'})$,
$\theta_{\calL_V}(\varphi_{nq,\la})$, $\theta_{\calL_V}(\varphi_{nq,[\la]})$ extend to smooth forms on the smooth manifold with corners $\overline{X}$. 

Moreover, for a standard rational parabolic $\underline{P}$, the restrictions 
$\tilde{r}_{{P}}$ to the corresponding boundary component $e'(\underline{P})$ of the three series above are given by
\[
\tilde{r}_{{P}} \left(
\theta_{\calL_V}(\varphi_{nq,\bullet})\right) =
\theta_{\widehat{\calL}_W} (r_{{P}}\varphi_{nq,\bullet}).
\]
\end{theorem}

\begin{remark}
The statement of the theorem is not correct for $\Q$-split case for signature $(p,p)$. In that case, one has to replace the Borel-Serre compactification for $\SO(p,p)$ to the slightly bigger Borel-Serre compactification for $\Orth(p,p)$, as we explain in the final section. With this modification the theorem holds again as stated above. 
\end{remark}

Combining Theorem~\ref{resformula} with Theorem~\ref{localrestrictiontheorem}, we obtain
\begin{corollary}\label{resformula22}
\[
\tilde{r}_{{P}}(\theta _{\mathcal{L}_V}(\varphi^V_{nq,\ell'}))= \widetilde{\iota}_P (\theta_{\widehat{\mathcal{L}}_W}(\varphi^W_{n(q-\ell),n\ell+\ell'})),
\qquad 
\tilde{r}_{{P}}(\theta _{\mathcal{L}_V}(\varphi^V_{nq,\lambda]})) = \widetilde{\iota}_P (\theta_{\widehat{\mathcal{L}}_W}(\varphi^W_{n(q-\ell),\ell\varphi_n+\la})),
\]
and
\[
[\tilde{r}_{{P}}(\theta _{\mathcal{L}_V}(\varphi^V_{nq,[\lambda]}))] = [\widetilde{\iota}_P (\theta_{\widehat{\mathcal{L}}_W}(\varphi^W_{n(q-\ell),[\ell\varphi_n+\la]}))].
\]
\end{corollary}

\begin{remark}\label{weightedremark}
More generally, the proof of Theorem~\ref{resformula} also shows that $\theta _{\mathcal{L}_V}(\varphi^V_{nq,[\lambda]})$ is ``essentially'' a special differential form in the sense of weighted cohomology, see \cite{GHMP}. Namely, $\tilde{r}_{{P}}(\theta _{\mathcal{L}_V}(\varphi^V_{nq,[\lambda]}))$ is $N_P$-invariant and while $\theta _{\mathcal{L}_V}(\varphi^V_{nq,[\lambda]})$ restricted to a neighborhood of $e'(P)$ in $\overline{X}$ is not the pullback by the geodesic retraction of its restriction, the difference 
of $\theta _{\mathcal{L}_V}(\varphi^V_{nq,[\lambda]})$ and this pullback has ``square-exponential'' decrease in the coordinates $t_i$ on $A_P$. In fact, one can distill out of our proof an explicit asymptotic expansion for $\theta _{\mathcal{L}_V}(\varphi^V_{nq,[\lambda]})$. This in particular proves that $\theta _{\mathcal{L}_V}(\varphi^V_{nq,[\lambda]})$ {\it extends to a smooth form on the smooth manifold with corners $\overline{X}$}. 

Moreover, the torus $A_P$ acts on the differential forms in \eqref{weightedforms} with weight
\[
\eta_n:= \alpha_1^{-n} \alpha_2^{-2n} \cdots \alpha_r^{-rn}
\]
(written multiplicatively). Hence (up to the exponentially decreasing part) the forms $\theta _{\mathcal{L}_V}(\varphi^V_{nq,[\lambda]})$ represent weighted cohomology classes with weight profile $\eta_n$. This is independent of the coefficients since $A_P$ acts trivially on the coefficient part of the restriction. In particular, for $p>q$ the classes are $L^2$ if $p>2n+1$. 

Finally, the proof shows that $\theta _{\mathcal{L}_V}(\varphi^V_{nq,\ell'})$ is exponentially decreasing in the direction of $e'(P)$ if $n>p-\ell$. In particular, $\theta _{\mathcal{L}_V}(\varphi^V_{nq,\ell'})$ is exponentially decreasing for $n=p$. 
\end{remark}

\begin{proof}[Proof of Theorem~\ref{resformula}]
It suffices to consider $\varphi^V_{nq,\ell'}$. For $g \in G$
and $g' \in G'$, we let
\begin{equation}
\theta^V_{{\alpha,\underline{\beta},I}}(g',g) = \sum_{\x \in L^n +h}
\omega_V(g')\varphi^V_{\Delta_{\alpha,\underline{\beta},I}}(g^{-1}\x)
\otimes  g^{\ast}\omega_{\alpha} \otimes
 g e_{\underline{\beta}}.
\end{equation}
be the theta series associated to one fixed component of
$\varphi^V_{nq,\ell'}$. For the purposes of studying the restriction
to $e'(\underline{P})$, we can assume $g'=1$ (since it intertwines
with the restriction) and also $g= a(\mathbf{t}) \in A$ (since $g$
varies in a Siegel set and by Lemma~\ref{mixedformulas}).  It also suffices to assume 
\[
\calL_V = \left( L^n_E +h_E\right) \oplus \left(L^n_W +h_w\right)  \oplus \left(L^n_{E'} + h_{E'}\right)
\]
with $L_E,L_W,L_{E'}$ lattices in $E,W,E'$ respectively. 

\begin{lemma}\label{mixedmodelthetaformula}
Let $a(\mathbf{t})\in A$. Then
\begin{align*}
 \theta^V_{{\alpha,\underline{\beta},I}}(a(\mathbf{t}))& = \det(L_{E})^{-n}
 \sum_{\x_W \in L^n_{W} +h_{W}} \sum_{\substack{\xi \in
(L_{E}^{\#})^n
\\ u' \in L^n_{E'} + h_{E'}}} e\left( 2\pi i (\xi,h_{E}) \right)  \\
& \quad \times |\mathbf{t}|^n
\widehat{\varphi}^E_{\Delta_{\alpha,\underline{\beta},I}^{\prime\prime}}
(\widetilde{\mathbf{t}}(\xi{^t}+u'),\widetilde{\mathbf{t}}u')
 \varphi^W_{\Delta'_{\alpha,\underline{\beta},I}}(\x_W)
\otimes a(\mathbf{t})^{\ast} \sigma^{\ast} \omega_{\alpha} \otimes
 a(\mathbf{t}) e_{\underline{\beta}}.
\end{align*}

\end{lemma}

\begin{proof}
This follows directly from Lemma~\ref{mixedformulas} and Poisson summation.
\end{proof}

\begin{lemma}
Assume at least one of the $\alpha_{kj}$ and $\beta_k$ is less
or equal than $\ell$. Then
\[
r_{\underline{P}} \theta^V_{{\alpha,\underline{\beta},I}} =0.
\]
\end{lemma}

\begin{proof}
By the hypothesis we have $W \ne 0$ for all parabolics $P$. Then $W \otimes \R u_i \subset \mathfrak{n}$ is a weight space for the action of $A_P$ with weight $t_i$. So in particular, for $a(\bf{t}) \in A_P$, we have that all components $t_i \to \infty$ as we approach $e'(P)$. Hence by Lemma~\ref{mixedmodelthetaformula} we
clearly see that each term in $\theta^V_{{\alpha,\underline{\beta},I}}(a(\mathbf{t}))$ is rapidly
decreasing as $t_i \to \infty$.
for $\underline{P}$ unless both $\xi=u=0$. But by
Lemma~\ref{mixedvarphi}, we have
\begin{equation}
\widehat{\varphi}^E_{\Delta_{\alpha,\underline{\beta},I}^{\prime\prime}}(0,0)
= \widehat{\varphi}^V_{\Delta_{\alpha,\underline{\beta},I}}\left(
\begin{smallmatrix}  0 \\ \x_W \\ 0
\end{smallmatrix} \right) = 0. \qedhere
\end{equation}
\end{proof}

Now for the remainder of the proof of Theorem~\ref{resformula},
assume that
\begin{equation}\label{primecond}
\alpha_{kj}, \beta_k \geq \ell+1.
\end{equation}
Again, each term in Lemma~\ref{mixedmodelthetaformula} is rapidly
decreasing unless $\xi=u=0$. So it suffices to consider
\begin{equation}\label{W-series}
\widehat{a(\mathbf{t})\varphi^V_{\Delta_{\alpha,\underline{\beta},I}}}
\left(
\begin{smallmatrix}  0 \\ \x_W \\ 0
\end{smallmatrix} \right)=
|\mathbf{t}|
\varphi^W_{\Delta'_{\alpha,\underline{\beta},I}}(\x_W) \otimes
a(\mathbf{t})^{\ast} \sigma^{\ast} \omega_{\alpha} \otimes
 a(\mathbf{t}) e_{\underline{\beta}}.
\end{equation}
Now $a(\mathbf{t}) e_{\underline{\beta}} = e_{\beta}$ by
\eqref{primecond}. We have
\begin{equation}
\sigma^{\ast} \omega_{\underline{\alpha_j}} = \frac{(-1)^{\ell}}{2^{\ell/2}} \omega_{\alpha_{j1}
p+1} \wedge \cdots \wedge \omega_{\alpha_{j q-\ell} m-\ell} \wedge
\nu_{\alpha_{j q-\ell+1} \ell } \wedge \cdots \wedge \nu_{\alpha_{j
q} 1 },
\end{equation}
and $A$ acts trivially on the $\omega_{\bullet}$'s, while for the
$\nu_{\bullet}$'s 
we have
$a(\mathbf{t})^{\ast}\nu_{ji} = \frac{db_{ji}}{t_i}$,
where $1 \leq i \leq \ell$ and $\ell+1 \leq j \leq m-\ell$. Here
$b_{ji}$ is the coordinate of $W \otimes E$ for $e_j \otimes u_i$
and $t_i$ is the parameter in $a(t_1,\dots,t_i,\dots,t_{\ell}) \in
A$. We obtain
\begin{align}
|\mathbf{t}| a(\mathbf{t})^{\ast}\sigma^{\ast} \omega_{\alpha} & = \frac{(-1)^{n\ell}}{2^{n\ell/2}}
|\mathbf{t}| \omega_{\alpha_{1\; 1}p+1} \wedge \cdots \wedge
\omega_{\alpha_{q-\ell \; 1} m-\ell} \wedge
\frac{db_{\alpha_{q-\ell+1\; 1} \ell}}{t_{\ell}} \wedge \cdots
\wedge
\frac{db_{\alpha_{q \;1} 1}}{t_1} \notag\\
& \quad \wedge \cdots \notag \\
& \quad \wedge \omega_{\alpha_{1\; n}p+1} \wedge \cdots \wedge
\omega_{\alpha_{q-\ell \; n} m-\ell} \wedge
\frac{db_{\alpha_{q-\ell+1\; n} \ell}}{t_{\ell}} \wedge \cdots
\wedge \frac{db_{\alpha_{q \;n} 1}}{t_1} \notag \\
&=  \frac{(-1)^{n\ell}}{2^{n\ell/2}}\omega_{\alpha_{1\; 1}p+1} \wedge \cdots \wedge
\omega_{\alpha_{q-\ell \; 1} m-\ell} \wedge db_{\alpha_{q-\ell+1\;
1} \ell} \wedge \cdots \wedge
db_{\alpha_{q \;1} 1} \label{weightedforms} \\
& \quad \wedge \cdots \notag\\
& \quad \wedge \omega_{\alpha_{1\; n}p+1} \wedge \cdots \wedge
\omega_{\alpha_{q-\ell \; n} m-\ell} \wedge db_{\alpha_{q-\ell+1\;
n} \ell} \wedge \cdots \wedge db_{\alpha_{q \;n} 1}. \notag
\end{align}
This shows for \eqref{W-series} we have
\begin{equation}
\widehat{a(\mathbf{t})\varphi^V_{\Delta_{\alpha,\underline{\beta},I}}}
\left(
\begin{smallmatrix}  0 \\ \x_W \\ 0
\end{smallmatrix} \right) =
r_{\underline{P}}\varphi^V_{\Delta_{\alpha,\underline{\beta},I}}
(\x_W)
\end{equation}
independent of $\mathbf{t}$. This completes the proof of
Theorem~\ref{resformula}.
\end{proof}

\subsection{Nonvanishing}

We now prove Theorem~\ref{TH:nonvan}.

By the hypotheses we can find a rational parabolic $\underline{P}$ such that $\dim E =\ell=q$, so $W$ is positive definite and $X_W$ is a point. Then by Theorem~\ref{resformula},
\begin{align}
[\tilde{r}_{P} \theta_{\calL_V}(\tau,\varphi_{q,[\la]}^V)] &= \tilde{\iota}_P[\theta_{\widehat{\calL}_W}(\tau,\varphi_{0,[\ell\varpi_n+\la]}^W] \\
& \quad  \in Mod(\G',m/2,\la) \otimes \tilde{\iota}_P \left(H^{0}(X_W,\mathbb{S}_{[\ell\varpi_n+\la]}(W_{\C}))\right) \notag \\
& \quad \quad \simeq  Mod(\G',m/2,\la)  \otimes \tau_{nq,\ell'} \left(\mathbb{S}_{[\ell\varpi_n+\la]}(W_{\C})\right) \notag \\
& \quad \quad \simeq  Mod(\G',m/2,\la)  \otimes \mathbb{S}_{[\ell\varpi_n+\la]}(W_{\C}). \notag
\end{align}
So in this case $\tilde{\iota}_P$ is an embedding. Hence the restriction to $e'(P)$ vanishes if and only if the positive definite theta series $\theta_{\widehat{\calL}_W}(\tau,\varphi_{0,[\ell\varpi_n+\la]}^W)$ vanishes. Furthermore, the restriction of the class $[\theta_{\calL_V}(\tau,\varphi_{q,[\la]}^V)]$ cannot arise from an invariant form on $D$, since in that case one would need to obtain the trivial representation in the coefficients.

To obtain the nonvanishing, we first observe

\begin{lemma}
Given $\varphi_{0,[\ell\varpi_n+\la]}^W$ as above, then there exists a coset of a lattice $\calL_W$ which we can take to be contained $\widehat{\calL}_W$ such that 
\[
\theta_{{\calL}_W}(\tau,\varphi_{0,[\ell\varpi_n+\la]}^W) \ne 0.
 \]
\end{lemma}

\begin{proof}
We give a very simple argument which we learned from E. Freitag and R. Schulze-Pillot. We can assume $V = \Q^m$ with the standard inner product. First find a vector $h \in \tfrac1{N_1}(\Z^m)^n$ with $N_1 \in \Z$ such that $\varphi_{0,[\ell\varpi_n+\la]}^W(h) \ne 0 $. Then pick a lattice $L= N_1N_2 \Z^m$ such that $\|\sum_{x\in L^n} \varphi_{0,[\ell\varpi_n+\la]}^W(x)\| < \|\varphi_{0,[\ell\varpi_n+\la]}^W(h)\|$. Such a $N_2 \in \Z$ exists as $\varphi_{0,[\ell\varpi_n+\la]}^W$ is a Schwartz function. Then the theta series associated to $\varphi_{0,[\ell\varpi_n+\la]}^W$ for $\calL_W = L^n+h$ does not vanish. 
\end{proof}

From this data then, we now can find a $\calL_V'$ contained in $\calL_V$ such that $\widehat{\calL'}_W = \calL_W$ with $\theta_{{\calL}_W}(\tau,\varphi_{0,[\ell\varpi_n+\la]}^W) \ne 0$. Replace $\G$ with $\G \cap\Stab{L'}$. Then $[\tilde{r}_{P} \theta_{\calL_V}(\tau,\varphi_{q,[\la]}^V)] \ne 0$. This proves Theorem~\ref{TH:nonvan}.

One feature of our method to establish non-vanishing is that we retain some control over the cover $X'$, since this reduces to the very concrete question of non-vanishing of positive definite theta series. An easy example for this is the following.

\begin{example}
Consider the integral quadratic form given by
\[
y_1y_1' + \cdots + y_q y_q' + 2x_1^2 +\cdots + 2x_k^2
\]
with $y_i,y_i',x_j \in \Z$. So $L = \Z^m$ with $m=2q+k$. Assume $k\geq q$. Note $L^{\#} \subset \tfrac14 \Z^m$. We let $\G$ be the subgroup in $\Stab(L)$ which stabilizes $L^{\#}/L$. Then
\[
H^q(\G,\Z) \neq 0.
\]
Using our method this follows from the non-vanishing of the theta series \linebreak $\sum_{\x \in \Z^k +(\tfrac14,\dots,\tfrac14)} x_1\cdots x_q e^{4\pi i (\sum x_i^2) \tau}$.

\end{example}

\section{The big Borel-Serre compactification for rational $\SO(p,p)$}\label{bigBS}

In this section, $V$ is always a $\Q$-split rational quadratic space of signature $(p,p)$ with Witt basis $u_1,\dots,u_{p-1},u_p, u'_{p},u'_{p-1}, \dots, u'_1$.

We will show that our main Theorem~\ref{resformula} remains true for the case of rational $\SO(p,p)$ but only if we replace the Borel-Serre compactification associated to the usual Tits building of type $D_p$ of (rational) parabolic subgroups of $\SO(p,p)$ by the ``big Borel-Serre compactification'' of type $B_p$ which will be described below. For this we have to change the underlying root system from type $D_p$ to type $B_p$ by adding reflections (and great subspheres in the Tits building). In terms of groups this is achieved by switching from $\SO(p,p)$ to the full orthogonal group $\Orth(p,p)$ (or equivalently, to $\SO(p+1,p)$).

Of course since both compactifications are compactifications of the same locally symmetric space the two boundaries assigned will be the same as topological spaces 
but  {\it their structures as manifolds with corners will be different}.

The main issue for us is that the parabolic subgroups of $\SO(p,p)$ do {\it not} correspond bijectively to isotropic flags, but rather to oriflammes, see Lemma~\ref{oriflamme}.

By switching to the root system $B_p$, i.e., considering $\Orth(p,p)$ or $\SO(p+1,p)$ we do obtain a bijection between parabolics and isotropic flags. This is the crucial aspect in constructing the big Borel-Serre compactification.

We first define the big Borel-Serre compactification {\it extrinsically} by embedding the locally symmetric space $X_{p,p} = \G_{p,p} \back D_{p,p}$ for $\SO(p,p)$ into a suitably constructed space $X_{p+1,p}= \G_{p+1,p} \back D_{p+1,p}$ for signature $(p+1,p)$ and then considering the closure of $X_{p,p}$ inside the Borel-Serre compactfication $\overline{X}_{p+1,p}$. The {\it intrinsic} big Borel-Serre compactification uses the Tits building for parabolic subgroups for the full orthogonal group $\Orth(p,p)$. 

The extension of $\theta(\varphi_{np,[\lambda]})$ is most easily established by pulling back the usual Borel-Serre compactification and restriction formulas for $(p+1,p)$ using the extrinsic definition. We proceed to give the intrinsic definition and compare the two constructions. It is then most instructive to compare the usual and the big Borel-Serre compactification. Finally, we consider the case of signature $(2,2)$ in more detail.

\subsection{The extrinsic big Borel-Serre compactification}

We set $\tilde{V} = V \perp \Q v$ with $(v,v) =1$. Hence $\tilde{V}$ has signature $(p+1,p)$. We rearrange coordinates so that $v$ becomes the $(p+1)$-st standard basis vector $e_{p+1}$. We write $\ell_{p+1} =\Q e_{p+1}$ for the line spanned by $e_{p+1}$. The natural inclusion $V \hookrightarrow \tilde{V}$ defines the inclusion $j_{p+1}:\Orth(p,p)\to \Orth(p+1,p)$. We will often identify $\Orth(p,p)$ with its image under $j_{p+1}$. The inclusion $j_{p+1}$ induces an inclusion (also denoted $j_{p+1}$) of the symmetric spaces $D_{p,p} \to D_{p+1,p}$.  
We let $\Gamma_{p+1,p}$ denote a congruence subgroup in $\SO(\tilde{V})$ stabilizing $\tilde{\calL} = \calL \oplus \Z v$ chosen so that it is torsion free and 
\[
\Gamma_{p,p} = \Orth(p,p) \cap \Gamma_{p+1,p}.
\]
We may assume, for example if $\Gamma_{p+1,p}$ is neat (the intersection of the subgroup of $\C^*$ generated by the elements of $\Gamma_{p+1,p}$ with the roots of unity is $\{1\}$) that this intersection is contained in $\SO(p,p)$. Let $\sigma \in \SO(\tilde{V})$ be the rational element that is $-1$ on $V$ and $1$ on $\ell_{p+1}$. Then $D_{p,p}$ is the fixed point set of $\sigma$ acting on $D_{p+1,p}$, that is 
\[
j_{p+1} D_{p,p} = D_{p+1,p}^{\sigma}.
\]
The inclusion of symmetric spaces induces a map (again denoted $j_{p+1}$) of locally symmetric spaces $j_{p+1}: X_{p,p} \to X_{p+1,p}$. Assume now that $\Gamma_{p+1,p}$ is torsion free. Then it follows from a well-known argument using $\sigma$  (the ``Jaffe Lemma'', Lemma 2.1 of \cite{Millson}) that $j_{p+1}$ induces an embedding of $X_{p,p}$ into $X_{p+1,p}$. 

\begin{definition}(The extrinsic big Borel-Serre compactification)
Assume $\G_{p,p}$ is torsion-free. The big Borel-Serre compactification $\overline{X}_{p,p}$ is the closure of $X_{p,p}$ in $\overline{X}_{p+1,p}$. We note that the inclusion $j_{p+1}$ induces an embedding $\overline{j}_{p+1}: \overline{X}_{p,p} \to \overline{X}_{p+1,p}$.  
\end{definition}

We will discuss the properties extrinsic big Borel-Serre compactification later in detail. At this point we can already give a quick proof that our theta series extend to the big compactification of $\overline{X}_{p,p}$. 

\begin{theorem} The forms $\theta(\varphi_{np,[\lambda]})$ on $X_{p,p}$ extend to the big Borel-Serre compactification $\overline{X}_{p,p}$.
\end{theorem}

\begin{proof}
Let $\widetilde{\varphi}_{np,[\lambda]}$ be the special $np$-cocycle for $\SO(p+1,p)$ and
 $\varphi_{np,[\lambda]}$ be the one for $\SO(p,p)$.
Note that from the explicit formulas for $\widetilde{\varphi}_{np,[\lambda]}$ and $\varphi_{np,[\lambda]}$ we have  
\begin{equation} \label{naturalityofKMforms}
j_{p+1}^* \widetilde{\varphi}_{np,[\lambda]} = \varphi_{np,[\lambda]}\varphi_0^{\ell_{p+1}}.
\end{equation}
Here $\varphi_0^{\ell_{p+1}}$ is the Gaussian associated to the $1$-dimensional positive definite subspace $\ell_{p+1}$. Since the lattice splits we obtain a corresponding restriction formula \begin{equation} \label{naturalityofthetafunctions}
j_{p+1}^*\theta( \widetilde{\varphi}_{np,[\lambda]}) =\theta( \varphi_{np,[\lambda]}) \theta(\varphi_0^{\ell_{p+1}})
\end{equation}
on the level of theta functions. Note that $\theta(\varphi_0^{\ell_{p+1}})$ is constant on $X_{p,p}$, so the product of the two factors on the right of (\ref{naturalityofthetafunctions}) extends to the big Borel-Serre boundary if and only if
the first factor extends.  
Now we have seen above that $\theta(\widetilde{\varphi}_{np,[\lambda]})$ extends over the Borel-Serre boundary of $X_{p+1,p}$.
The lemma then follows by considering the following commutative diagram
(starting with $\theta(\widetilde{\varphi}_{np})$ in the lower left-hand corner).
\[
\begin{CD}
A^{\bullet}(\overline{X}_{p+1,p}) @ >\overline{j}_{p+1}^*>> A^{\bullet}(\overline{X}_{p,p})\\
@VVV                                 @VVV \\
A^{\bullet}(X_{p+1,p}) @ >j_{p+1}^*>> A^{\bullet}(X_{p,p}).
\end{CD}
\]
\vskip-.6cm \end{proof}

\begin{remark}
We are required so far to assume that the lattice $\Gamma_{p+1,p}$ and hence $\Gamma_{p,p}$ is torsion-free. However, after we have given the intrinsic description of our compactification and hence we know that this {\it intrinsic construction} produces a compactification for the quotient of $D$ by a normal torsion-free subgroup
$\Gamma'$ of $\Gamma \subset \SO(p,p)$, then the extension and the restriction formula will hold for the quotient by the larger lattice $\Gamma$ because it is invariantly defined. We leave the details to the reader.
\end{remark}

\subsection{The intrinsic description of the new compactification}

We now give an intrinsic description of the big Borel-Serre compactification, that is, it does not use the embedding $j_{p+1}$. 

In what follows if $G$ is any reductive group we will use $\mathcal{P}(G)$ to denote the set of  parabolic subgroups of $G$. 

There are four key ingredients of a Borel-Serre compactification, see \cite{BJ}, III.9 (and Section~\ref{BS-C1} above). 

\begin{enumerate}

\item The Tits building $\mathcal{B}(G)$ (or rather its quotient by the arithmetic group $\Gamma \subset G$ under consideration). 

\item For each rational parabolic $P$ of $G$ there is the split torus $A_P$ which is the connected component of the identity of the center of $P/N$. 

\item For each rational parabolic subgroup $P$ there is the associated ``Borel-Serre face'' $e(P):= P/A_P K_P$. Here $K_P= P \cap K$ is as before the subgroup of $P$ that stabilizes the basepoint $z_0$ of the associated symmetric space.  

\item 
The set $\Phi(P,A_P)$ of restrictions of the set of positive roots to $A_P$, which governs the topology around the boundary faces, in particular, convergence to a point in the boundary. The reader will note the definition of convergence will not be changed if the elements of $\Phi(P,A_P)$ are replaced by positive scalar multiples. Furthermore, one obtains the same set of convergent sequences if in the rule \cite{BJ}, p.328, one replaces $\Phi(P,A_P)$ by $\Delta(P,A_P)$ the set restrictions to $A_P$ of the {\it simple} roots in the root system associated to the maximal torus $A_{P_0}$ for a chosen minimal parabolic $P_0$. 
\end{enumerate}

\begin{definition}(The intrinsic big Borel-Serre compactification)
The intrinsic big Borel-Serre compactification $\overline{X}_{p,p}$ is obtained
by applying the ``uniform construction of Borel-Ji'' (\cite{BJ}, \S III.9) to the Tits building $\mathcal{B}(\Orth(p,p))$ for the full orthogonal group together with the root system $B_p$. 
\end{definition}

The term ``intrinsic compactification'' is a bit premature since $\Orth(p,p)$ since one still needs to show that the construction really gives a compact space. At this point it is only a formal procedure. Moreover, it is a priori not clear that we can freely change the root system from $D_p$ to $B_p$. Only once we have established the equivalence to the extrinsic description this will be justified. Note however, that the full orthogonal group $\Orth(p,p)$ gives rise to the same symmetric space as $\SO(p,p)$.

We now describe some of the features of the new construction. 

\subsubsection{The new building $\mathcal{B}(\Orth(p,p))$ and the map of parabolic subgroups}

Recall that we defined the standard totally isotropic subspaces $E_k = \Span(u_1,\dots, u_k)$ in $V$ and the spaces $E_+ = E_p=\Span(u_1,\dots,u_{p-1}, u_p)$ and $E_- = \Span(u_1,\dots u_{p-1}, u_p')$. 

We first note (see eg \cite{AB}, \cite{Garrett})

\begin{lemma}\label{simplicialcomplex}
The (standard) parabolic subgroups of $\Orth(p,p)$ are the stabilizers of the (standard) isotropic flags (in $E_p$), and every isotropic flag determines a parabolic. Thus the associated Tits building $\mathcal{B}(\Orth(p,p))$ is the spherical building associated to the partially ordered set of isotropic flags in $V$ and the parabolic subgroups of $\Orth(p,p)$ are the stabilizers of the faces of the building. 
\end{lemma}

\begin{example}\label{p1flag2}
We illustrate this fundamental difference to the special orthogonal group $\SO(p,p)$. Let $P \subset \Orth(p,p)$ be the stabilizer of the isotropic subspace $E_{p-1}$. Then
\[
P =\left\{
\begin{pmatrix}
g  & c_2 \; c_3& \dots \\  
0   &  h & \dots \\
0   &  0      & g^{\ast}
\end{pmatrix}\right\} 
\]
with $g \in \GL_{p-1}(\R)$, $h \in \Orth(1,1)$, $c_i \in\R^{p-1}$ (column vectors) and $g^{\ast}$ as in \eqref{O-GL-embed}. Note $\Orth(1,1)  =  \SO(1,1) \cup w \SO(1,1)$ and $\SO(1,1) =\left\{ \kzxz{b}{}{}{b^{-1}} \right\}$. Here $w= \kzxz{0}{1}{1}{0}$. Hence $P$ is a maximal parabolic subgroup of $\Orth(p,p)$. 

Now consider $P' = P \cap \SO(p,p)$, the stabilizer of $E_{p-1}$ in $\SO(p,p)$. Now we have
\[
P'  = \left\{
\begin{pmatrix}
g  & c_2 & c_3  & \dots \\
0   &  b  & 0  & \dots \\
0   &  0  & b^{-1} & \dots \\
0   &  0  &   0  & g^{\ast}
\end{pmatrix} \right\}.
\]
Thus $P'$ is strictly contained in the stabilizer of both isotropic $p$-planes $E_+$ and $E_-$. Hence {\it is not a maximal parabolic} and we can associate $P'$ to {\it two} isotropic flags, namely $(E_{k-1},E_+)$ and $(E_{k-1},E_-)$; i.e., the oriflamme $(E_+,E_-)$. 
\end{example}

The situation in general is as follows. 

\begin{definition} \label{bad}
We say an isotropic flag ${\bf F}$ in $V$ is bad if an isotropic subspace of dimension $p-1$ occurs in ${\bf F}$. We say a parabolic in $\Orth(p,p)$ is bad if it stabilizes is a bad flag. Otherwise we call ${\bf F}$ and $P_{\bf F}$ good. 
\end{definition}

We then have

\begin{lemma}\label{PP'map}
Let $P \subset \Orth(p,p)$ be a parabolic subgroup stabilizing the flag ${\bf F}$. Set $P' = P \cap \SO(p,p)$. 

\begin{itemize}

\item[(i)]
Assume $P$ is good. Then $P'$ is the stabilizer of the flag ${\bf F}$ (see also Lemma~\ref{oriflamme} (2)). 

\item[(ii)]
Assume $P$ is bad stabilizing a flag $F_1 \subset \dots F_k \subset F_{p-1} (\subset F_p)$ with $\dim F_{p-1} =p-1$ and $\dim F_p = p$. ($F_p$ might or might not be there). Let $F_{p,1}, F_{p,2}$ as in Lemma~\ref{oriflamme} (3). Then $P'$ is the stabilizer of the oriflamme $(F_1,\dots,F_k, F_{p,1},F_{p,2})$. 
\end{itemize}
\end{lemma}

We now describe how each top dimensional simplex of the Tits building $\mathcal{B}(\SO(p,p))$ of type $D_p$ will be bisected to obtain $\mathcal{B}(\Orth(p,p))$. Each spherical chamber (top dimensional i.e. $p-1$ dimensional simplex) contains a distinguished edge $e$ (the edge joining  the two vertices corresponding to highest dimensional isotropic subspaces. Let $f$ be the $p-3$ face that is opposite to $e$. Hence the chamber is the join $e * f$. Let $b$ be the barycenter of $e$. Then we bisect each spherical chamber by the codimension one  interior simplex $b * f$. We make a choice of one of the two halves of the original spherical fundamental chamber $\Delta_{D_p} =\Delta'$ and call it the fundamental spherical chamber $\Delta_{B_p} =\Delta$ of $\mathcal{B}(\Orth(p,p))$. 
The resulting nonthick building is the building of type $B_p$ on which the big Borel-Serre  compactification will be modeled. Note that if $F$ is a face of
$\mathcal{B}(\Orth(p,p))$ then there will be a unique face $F'$ of $\mathcal{B}(\SO(p,p))$ such that the interior $F^0$ is contained in $F'$. 

Since the parabolic subgroups are exactly the subgroups that fix faces of the buildings, the map $F \mapsto F'$ induces a map $\mathcal{P}(\Orth(p,p)) \to \mathcal{P}(\SO(p,p))$ of parabolic subgroups. In fact, it is exactly the assignment $P \mapsto P' = P \cap \SO(p,p)$ in Lemma~\ref{PP'map}. For good flags the claim is obvious, since in that case by definition $P'$ is the subgroup of $\SO(p,p)$ that fixes the same face $F$. Thus the only difficulty is when the face $F$ corresponds to a bad flag. In this case the face $F$ fixed by the original parabolic $P$ has dimension one less than $\mathcal{F}(F)$. But in this case $F^0$ is contained in the interior of $F'$ and if $g \in \SO(p,p)$ fixes an interior point to the face $F'$ then it fixes all of $F'$. The claim follows.
Note that $F \mapsto F'$ is a bijection on faces of dimension less than or equal to $p-1$ but it is two-to-one on top faces.

\subsubsection{The new split central split torus $A_P$.}\label{newcenter}
We define the subtorus $A_P$ of $A_{P_0}$ to be the center of $L = P \cap P^{\theta_0}$ where $\theta_0$ is the Cartan involution corresponding to our chosen basepoint $z_0$. Note that we cannot define it as the annihilator of an appropriate subset $I$ of the simple roots of $SO(p,p)$. However we can define it as the annihilator of an appropriate subset $I$ of the simple roots of the {\it new root system of type $B_p$}, see below, in particular, Lemma \ref{actiononN}.
These roots are defined intrinsically only up to positive multiples but this is enough to unambiguously define $A_P$. We will denote the new torus  $A_P$.  

\subsubsection{The new face $e(P)$}\label{newface}
Given $A_P$, we define the associated face $e(P)$ of the Borel-Serre enlargement by $e(P) = P/ A_P K_P$. Hence the cells $e(P)$ are assembled using the simplicial complex associated  to the partially ordered set of isotropic flags in $V$. 
The point is that the split torus $A_P$ can be {\it strictly ($1$-dimension) smaller} for certain parabolics in the new compactification (because $P$ and its Levi subgroup $L$ will have extra connected components causing its center to be smaller, see Example~\ref{p1flag2}) and consequently the face $e(P)$ will be strictly larger. In Theorem \ref{intrinsiccompactification} we will record this in detail. 

\subsubsection{The new system of roots of type $B_p$ and the set $\Phi(P,A_P$)}\label{newpositiveroots}
Fourth, there is a subset of the positive roots restricted to $A_P$ to be denoted $\Phi(P,A_P)$ and the correponding system of simple roots restricted to $A_P$ to be denoted $\Delta(P,A_P)$. This is the most complicated change to describe intrinsically. We define the Weyl group $W$ of the maximal torus $A_{P_0}$ as usual as the normalizer in $\Orth(p,p)$. But now the element
\[
w = \left( \begin{smallmatrix} I_{p-1} & & & \\  & 0 &1 & \\ & 1 &0 & \\ & & & I_{p-1} \end{smallmatrix} \right)
\]
is in $W$. Hence the Weyl group for $\Orth(p,p)$ is strictly larger than the one for $\SO(p,p)$. In fact, with this additional reflection (which interchanges $u_p$ and $u_p'$) one obtains the Weyl group for the root system $B_p$. While this does not define directly the new roots it defines the root hyperplanes. The choice of the fundamental chamber in the new Tits building defines a positive Weyl chamber, equivalently the correct orientation of the hyperplanes. (This corresponds to the choice of defining the standard parabolics in $\Orth(p,p)$ to be the stabilizer of flags in $E_+$ or $E_-$). For each root hyperplane we choose a linear functional which vanishes on the hyperplane and is positive on the cone on $\Delta_{B_p}$. This new collection of linear functionals we will call the
(new) positive roots to be denoted $\Phi$. In terms of the Tits building this amounts to the following. We have already added the new walls to the spherical building at infinity and choosen the fundamental spherical chamber $\Delta_{B_p}$. We now extend them inside $A_{P_0}$ to obtain the standard linear action of the Weyl group of type $B_p$ as a reflection group. In more detail, given the split torus $A_{P_0}$ which we identify with its Lie algebra $\mathfrak{a}$, we consider the corresponding apartment $\mathcal{A}$ in
$\mathcal{B}(\Orth(p,p))$ (the boundary of $A_{P_0}$). The building structure on $\mathcal{B}(\Orth(p,p))$ gives us a collection of great spheres in the apartment 
$\mathcal{A}$. If we regard the apartment $\mathcal{A}$ as the sphere at infinity of $A_{p_0}$ (each ray leaving the origin of $A_{P_0}$ corresponds to a unique point of $\mathcal{A}$, then the collection of great spheres corresponds (to the boundaries of) a collection of hyperplanes in $A_{P_0}$. Reflections in these hyperplanes give rise to the standard representation of the Coxeter group of type $B_p$. The chosen spherical chamber
$\Delta_{B_p}$ corresponds to a Weyl chamber in $A_{P_0}$ which we will also denote
$\Delta_{B_p}$.

\begin{definition}
$\Phi(P,A_P)$ is  the set of restrictions to $A_P$ of the roots
in $\Phi$.
\end{definition}

\begin{remark}
 We did not use the Lie algebra $\mathfrak{n}$ of $P$ in this definition.  We will see later that what we are doing is pulling back the usual
$A_P$ and $\Phi(P,A_P)$ from $\SO(p+1,p)$ using the embedding $j_{p+1}$.
\end{remark}

\subsection{The intrinsic and the extrinsic big Borel-Serre compactification coincide}

\begin{theorem}  \label{intrinsiccompactification}
The intrinsic and the extrinsic big Borel-Serre compactification of coincide. In particular, the cells $e'(P)$ are assembled using the simplicial complex associated to the partially ordered set of isotropic flags in $V$. 
\end{theorem}

From this we now easily check that all results from Section~\ref{globalressection} carry over with no change to the big Borel-Serre compactification for the split $(p,p)$-case. In particular,

\begin{theorem}
The restriction theorems, Theorem~\ref{resformula} and Corollary~\ref{resformula22}, hold in the big Borel-Serre compactification of $X_{p,p}$. 
\end{theorem} 

\begin{remark}
In fact, the restriction in the small Borel-Serre compactification to faces associated to good parabolics goes through as before as well with no change. It is the restriction to bad faces which causes problems. 
\end{remark}

To prove Theorem~\ref{intrinsiccompactification} we will first prove the analogue of the theorem for the partial compactifications (Borel-Serre enlargements) of the symmetric spaces $D_{p,p}$ 
and $D_{p+1,p}$. We will denote the corresponding enlargements by $\overline{D}_{p,p}$ (constructed using  $\mathcal{P}(\Orth(p,p))$) and $\overline{D}_{p+1,p}$. Recall that earlier we already saw $D_{p,p} = D_{p+1,p}^{\sigma}$. We claim the corresponding equation also holds for the enlargements. We have
 
\begin{proposition} \label{enlargement}
\begin{itemize}
\item[(i)] 
$\overline{D}_{p,p} = \overline{D}_{p+1,p}^{\sigma}.$
\item[(ii)]  
Let $\widetilde{P}$ be the stabilizer of an isotropic flag $\widetilde{\bf F}$ in $\widetilde{V}$ and suppose $\widetilde{P}$ is normalized by  $\sigma$.  
Then the subspaces of the flag $\widetilde{\bf F}$ are in fact contained in $V$.
We let ${\bf F}$ be the associated isotropic  flag in $V$ and $P$ be the stabilizer of ${\bf F}$ whence  
\[
P = \widetilde{P}^{\sigma}.
\]
\item[(iii)]
Suppose $e(\widetilde{P})^{\sigma}$ is nonempty.  Then $\widetilde{P}$ is normalized by $\sigma$ and 
\[
e(\widetilde{P})^{\sigma} = e(P).
\]
\item[(iv)] 
$\overline{D}_{p+1,p}^{\sigma} = D_{p,p}  \coprod \coprod_{P \in  \mathcal{P}(\Orth(p,p))} e(P)$.
 \end{itemize}
 \end{proposition}

On the building level this means that the map $j_{p+1}$ induces a simplicial embedding of $\mathcal{B}(\Orth(p,p))$ onto $\mathcal{B}(\SO(p+1,p))^{\sigma}$ carrying apartments isomorphically onto apartments. The image is the fixed subbuilding $\mathcal{B}(\SO(p+1,p))^{\sigma}$.

The proposition will be a consequence of the following discussion.

We note that the inclusion $\overline{D}_{p,p} \subset  \overline{D}_{p+1,p}^{\sigma}$
is obvious. The reverse inclusion will follow once we have proved (iv). We immediately see
\[
\overline{D}_{p+1,p}^{\sigma} = D_{p+1,p}^{\sigma} \coprod  \coprod_{P \in  \mathcal{P}(\SO(p+1,p))}  e(\widetilde{P})^{\sigma}.
\]
Clearly (iv) will follow from (iii). (ii) and (iii) will be a consequence of the next three lemmas. In order to prove (iii) we need to first prove (ii).

 \begin{lemma}\label{isotropicsubspace}  
 Suppose $\widetilde{E}$ is an isotropic subspace of $\widetilde{V}$ such that $\sigma(\widetilde{E}) = \widetilde{E}$.  Then 
 $\widetilde{E}  \subset V$.
 \end{lemma}
 
 \begin{proof} We have $\widetilde{E} =  (\widetilde{E} \cap \ell_{p+1}) \oplus (\widetilde{E} \cap V)$. But as $\widetilde{E}$ is isotropic we see $\widetilde{E} \cap \ell_{p+1} =0$.
 \end{proof}

We now show that Lemma \ref{isotropicsubspace} implies (ii). Indeed,
$\widetilde{P}$ is the stabilizer of a unique isotropic flag $\widetilde{\bf F}$.
Now since $\widetilde{P}$ is its own normalizer and we are assuming $\sigma$ normalizes $\widetilde{P}$ we find $\sigma \in \widetilde{P}$ and consequently
$\sigma$ carries each of the subspaces in $\widetilde{\bf F}$ into itself. Hence by Lemma \ref{isotropicsubspace} each of these subspaces is contained in $V$. We 
let ${\bf F}$ denote the associated isotropic flag in $V$ and let $P$ be its stabilizer in $\Orth(p,p)$.  We now prove that $\widetilde{P}^{\sigma} =P$.
First we claim  that $\widetilde{P}^{\sigma}$ contained in $\Orth(p,p)$. Indeed, since  $g \in \widetilde{P}^{\sigma}$ we have  $g^{-1} \sigma g = \sigma$ whence $g$ carries the line through
$e_{p+1}$ into itself whence $g \in \Orth(p,p)$. But also by definition $\widetilde{P}^{\sigma}$ fixes ${\bf F}$ whence we have 
\[
\widetilde{P}^{\sigma} =P.
\]

Thus it remains to prove (iii).  This we do in the next two lemmas.
\begin{lemma}
If $e(\widetilde{P})$ contains a fixed point of $\sigma$ then
$\sigma \in \widetilde{P}$ and hence  $\sigma(e(\widetilde{P})) = e(\widetilde{P})$. In fact, we have 
 \begin{equation}\label{stabilizerofface}
 \sigma(e(\widetilde{P})) = e(\widetilde{P}) \iff   \sigma \in \widetilde{P}.
 \end{equation}
\end{lemma} 

\begin{proof}
 It follows from the basic result  of \cite{BS}, Corollary 7.7 (1) (with $P =Q$), that  
 \[
 \sigma(e(\widetilde{P})) \cap e(\widetilde{P}) \neq \emptyset \iff \sigma \in \widetilde{P}. \qedhere
 \]
 \end{proof}

 \begin{lemma}
 Suppose  $P= \widetilde{P}^{\sigma}$. Then we have
\begin{equation} \label{fixedpointsofcells}
e(P) = e(\widetilde{P})^{\sigma}.
\end{equation}
\end{lemma}

\begin{proof}
We only need to show $e(P) \subset e(\widetilde{P})^{\sigma}$. 
So suppose $x \in e(\widetilde{P})$ is fixed by $\sigma$. Let $y$ be the diagonal matrix with $p+1$ ones followed by $p$ minus ones.  Then conjugation by $y$ induces  Cartan involutions of $\SO(p+1,p)$ and $\Orth(p,p)$.  It is standard that we may construct a Levi decomposition  $\widetilde{P} = \widetilde{M} \cdot \widetilde{N}$ with $\widetilde{M} = \widetilde{P} \cap y \widetilde{P} y^{-1}$ whence $\sigma \in \widetilde{M}$. Note that 
$$e(\widetilde{P})= (\widetilde{M}  \widetilde{N})/ \widetilde{K} \cap \widetilde{M}.$$ 
Choose a lift $x' = \widetilde{m} \widetilde{n}$ of $x$ to $\widetilde{P}$.
Then $x$ is fixed under $\sigma$ implies that $\widetilde{n}$ is fixed under
$\sigma$ which implies $\widetilde{n}$ is in the unipotent radical $N$ of $P$.
Also  $\widetilde{m}$ is fixed modulo $\widetilde{K} \cap \widetilde{M}$.
Thus it remains to show that the group $M =\widetilde{M}^{\sigma}$ acts transitively on the fixed point set of $\sigma$ on its associated symmetric space
$ \widetilde{M} / (\widetilde{K} \cap \widetilde{M})$.  But the fixed point set 
is connected (because the unique geodesic joining any two fixed points must also be fixed).  Hence we may obtain the fixed point set by exponentiating the fixed subspace of $\widetilde{\mathfrak{p}}$ the tangent space to $D_{p+1,p}$ at the point $z_0$ fixed by the above Cartan involution.  But this fixed subspace is  $\mathfrak{p}$, the tangent space to $D_{p,p}$ at $z_0$.
\end{proof}

We have now completed the proof of Proposition \ref{enlargement}.

We also need to show that the convergence criterion applied to the topology of $\overline{D}_{p,p}$ is induced from the topology of $\overline{D}_{p+1,p}$ (and hence using the root system of type $B_p$). This follows from the following Lemma which the reader will verify. 

\begin{lemma} \label{actiononN} 
$\Phi(P,A_P)$ is the set of weights of $A_P$ acting on the nilradical $\widetilde{\mathfrak{n}}$ of the parabolic subalgebra of the corresponding parabolic $\widetilde{P}$ ($\widetilde{P}^{\sigma} = P$) via the inclusion $j_{p+1}:L \to \widetilde{L}$. 
\end{lemma}

Theorem \ref{intrinsiccompactification} will follow from the next Lemma. 
  
\begin{lemma} \label{Jaffee} 
Suppose $\Gamma_{p+1,p}$ is torsion free and there exists $\gamma \in \Gamma_{p+1,p}$ such that
$ \gamma(e(P)) \cap e(P) \neq \emptyset$.  Then $\gamma \in P \cap \Gamma_{p,p}$. In particular, 
the image of $e(P)$ in $\overline{X}_{p+1,p}$ is the quotient of $e(P)$ by $P \cap \Gamma_{p,p}$.
\end{lemma}
\begin{proof}  
Suppose $x \in e(P)$ satisfies that $y = \gamma(x) \in e(P)$.
Then  $\sigma \gamma^{-1} \sigma \gamma (x) = x$ since $\sigma$ fixes $x$ and $y$.  But the action of $\Gamma_{p+1,p}$ on the Borel-Serre enlargement of $D_{p+1,p}$ is fixed-point free since by \cite{BS}, Theorem 9.3, it acts properly and we have assumed it is torsion free .  Hence $\sigma \gamma \sigma = \gamma $
and consequently $\gamma \in \Gamma_{p,p}$. The lemma now follows from
Corollary 7.7 (1) of \cite{BS}.
\end{proof}

This concludes the proof of Theorem \ref{intrinsiccompactification}.

\subsection{Relating the small and the big Borel-Serre compactification of $X_{p,p}$}

We now have two compactifications of $X_{p,p}$, the usual Borel-Serre compactification and the new ``big '' Borel-Serre compactification we have just described. 

For $P$ a parabolic in $\Orth(p,p)$, we will write $P'=P \cap \SO(p,p)$ as before. We will denote the corresponding face in the small Borel-Serre enlargement by $e(P')$.

\begin{proposition}\label{comparisonofitemsgoodcase}
Suppose $P$ is a good parabolic in $\Orth(p,p)$. Then 
\begin{enumerate}
\item  $e(P) = e(P')$.
\item  $A_P = A_{P'}$.
\item  If the last subspace in the flag has dimension strictly less than $p$
(hence strictly less than $p-1$) then 
\[
\Phi^{B_p}(P, A_P) = \Phi^{D_p}(P', A_{P'}).
\]
\end{enumerate}
If the last element in the flag has dimension $p$ then $\Phi^{B_p}(P, A_P)$ and $ \Phi^{D_p}(P', A_{P'})$
will coincide except for the last entry which in the first case will
be the restriction of $t_p$ and in the second case will be the restriction
of $t_p^2$ (the squaring makes no difference in terms of the convergence criterion). 
\end{proposition} 

We will leave the proof of this proposition to the reader.

We now state what happens if $P$ is bad. We may assume that the associated flag is standard, contained in the totally isotropic subspace $E_p=E_+$.

\begin{proposition}\label{comparisonofitemsbadcase}
Suppose $P$ is a bad parabolic in $\Orth(p,p)$. There are two cases. 

\smallskip
 
(i) Suppose first the last subspace in the flag has dimension $p-1$. Then 
\begin{enumerate}
\item $e(P) \cong e(P')  \times \R_+.$
\item  $A_P \times \R_+ = A_{P'}$. Note that there is a projection
map $\pi_p:A_{P'} \to A_P$ which omits the last coordinate $t_p$.
This map is split by the map $i_p:A_P \to A_{P'}$ which puts a one
in the last component. 
\item Then $\Delta^{B_p}(P, A_P)$ is the set of restrictions of the old simple
roots of type $D_p$ to $A_P$ and $\Delta^{D_p}(P', A_{P'})$ is the set of restrictions of the old simple roots of type $D_p$ to the larger torus
$A_{P'}$. This may be restated as follows. We may identify $A_P$ and $A_{P'}$ with quotient tori of $A$ and hence we may identify their character groups with subgroups of the character group of $A$.  Suppose that $A_{P'}$
has dimension $r+1$ whence $A_P$ hence dimension $r$.  Then $|\Delta^{D_p}(P', A_{P'})| =r+1$ and $|\Delta^{B_p}(P, A_P)| =r$.
Then the first $r-1$ elements of the two sets of restricted simple roots ``coincide'' in the sense that as characters of $A$ they are the  pull-backs of the restrictions of the  roots
$t_i/t_{i+1},1 \leq i \leq p-2,$ to $A_P$ and $A_{P'}$ (so some of these may be trivial),  the last element of $\Delta^{B_p}(P,, A_P)$ is $t_{p-1}$
and the last two elements of $\Delta^{D_p}(P', A_{P'})$ are $t_{p-1}/t_p$ and $t_{p-1}t_p$. 
\end{enumerate}

\smallskip 

(ii) Now suppose the last element in the flag has dimension $p$, so the last two elements
are $E_{p-1}$ and $E_p$, then 
\begin{enumerate}
\item $e(P) = e(P')$.
\item $A_P = A_{P'}$.
\item  $\Delta^{D_p}(P', A_{P'})$  and $\Delta^{B_p}(P, A_P)$ have the same cardinality $r$, and  their first $r-1$ elements coincide. The last two nontrivial elements of $\Delta^{B_p}(P, A^{B_p})$ are $t_{p-1}/t_p$ and $t_p$
and the last two nontrivial elements of
$\Delta^{D_p}(P', A_{P'})$ are the restrictions of $t_{p-1}/t_p$ and $t_{p-1}t_p$.  
\end{enumerate}
\end{proposition}  

\begin{proof}  We prove (i) for the special case in which $P$ is the stabilizer of the isotropic subspace $E_{p-1}$, see Example~\ref{p1flag2}. For $P' = P \cap \SO(V_\R)$ we easily see
 \[
 A_{P'} =  \left\{ \left( \begin{smallmatrix}
aI_{p-1}  & 0 & 0 & 0 \\
0   &  b  & 0  & 0 \\
0   &  0  & b^{-1} & 0\\
0   &  0  &   0  & a^{-1} I_{p-1}
\end{smallmatrix} \right); \; a,b \in \R_+ \right\}
\]
and 
\[
\Delta^{D_p}(P',A_{P'}) = \{ a/b,ab \}.
\]
Consequently if $Y_{p-1}$ denotes the symmetric space associated to $\SL(E_{p-1})$ we have a diffeomorphism (ignoring the fiber bundle structure)
 \[
 e(P')) \cong Y_{p-1}  \times (W \otimes E_{p-1}) \times \wwedge{2} E_{p-1}
 \]
with $W = \Span(u_p,u_p')$. But for the Levi of $P$ in the full group $\Orth(p,p)$, we have 
\[
Z(L) = Z(L \cap \SO(p,p)) \cap Z(w)
\]
with $w$ as in Example~\ref{p1flag2} whence we have 
\[
 A_P  = \left\{ \left( \begin{smallmatrix}
aI_{p-1}  & 0 & 0 & 0 \\
0   &  1 & 0  & 0 \\
0   &  0 & 1 & 0\\
0   &  0  &   0  & a^{-1} I_{p-1}
\end{smallmatrix} \right); \; a \in \R_+ \right\}
\]
and 
\[
\Delta^{B_p}(P,A_P) = \{ a \}.
\]
Hence  we have a diffeomorphism
\[
e(P)  \cong Y_{p-1} \times \R_+ \times (W \otimes E_{p-1}) \times \wwedge{2} E_{p-1}.
\]

For (ii) suppose the last subspace has dimension $p$. For convenience we assume $P$ is the stabilizer of the flag $(E_{p-1}, E_+)$. Then 
\[
 P =\left\{
\begin{pmatrix}
g  & c  &  \dots & \dots\\
0   &  b  & \dots &\dots \\
0   &  0  & b^{-1} & \dots \\
0   &  0  & 0  & g^{\ast}
\end{pmatrix} \right\} 
\qquad
\text{and}
\quad L = \left\{
\begin{pmatrix}
g  &  0  &  0 & 0\\
0   &  b  & 0 & 0 \\
0   &  0  & b^{-1} & 0 \\
0   &  0  &  0 & g^{\ast}
\end{pmatrix} \right\} 
\]
with $g \in \GL_{p-1}(\R)$, $c \in \R^{p-1}$, $b \in \R^{\ast}$. Hence
\[
A_P =A_{P'} =\left\{
\begin{pmatrix}
aI_{p-1}  & 0  & 0 & 0\\
0   &  b  & 0 & 0 \\
0   &  0  & b^{-1} & 0 \\
0   &  0  & 0  & a^{-1} I_{p-1}
\end{pmatrix}; \; a,b \in \R_+ \right\}, 
\]
but
\[
\Delta^{B_p}(P,A_P) = \{a/b,b\} \qquad \text{and} \qquad \Delta^{D_p}(P',A_{P'}) = \{a/b,ab\}. \qedhere
\]
\end{proof}

\subsection{Signature $(2,2)$}

We now consider the case of signature $(2,2)$ in detail. In particular, we illustrate in this case the failure of the restriction formula for the small Borel-Serre compactification. 

\subsubsection{Comparison of the two compactifications} \label{thedifferences}

For $\SO(2,2)$ each apartment of the underlying Tits building (the building of parabolic subgroups of $\SO(2,2)$) is a square; the building of type $D_2 = A_1 \times A_1$. In the usual Borel-Serre compactification each of the four vertices is blown up to a circle bundle over a quotient of the upper half plane by a subgroup of finite index in $SL(2,\Z)$, i.e., modular curves.  Each edge is blown up to a $2$-torus, the two circle bundles over the modular curves corresponding to the two vertices of the edge are glued along this torus.

We now describe the big Borel-Serre compactification. In this case the underlying building (the nonthick Tits building associated to the complex of isotropic flags in $\Q^{2,2}$) has apartments which are octagons.  We will regard these octagons as the barycentric subdivisions of the above squares. We blow up the original vertices to the same circle bundles over modular curves as before. We blow up the four new vertices (the barycenters of the original edges) to trivial $2$-torus bundles over $\R_+$ compactified by adding two points $0$ and $\infty$.  We can glue the four new three manifolds to the four old ones because each has  boundary components homeomorphic to the $2$-torus. There is one such glueing for each of the eight edges of the octagon. It is critical to observe that not only do we use a new glueing scheme, the nonthick building of type $B_2=C_2$ associated to the isotropic flag complex but also there are some new cells $e(P)$ that do not occur in the usual Borel-Serre compactification.

In detail, we consider one fixed edge of the apartment of the Tits building for $\SO(2,2)$ corresponding to the basis $\{u_1,u_2, u_2', u_1'\}$. Namely, we let $Q_\pm'$ be the maximal parabolic in $\SO(2,2)$ stabilizer of the isotropic plane $E_\pm$ spanned by $u_1,u_2$ and $u_1,u_2'$ respectively. The intersection $P' = Q'_+ \cap Q'_-$ stabilizes the oriflamme $(E_+,E_-)$. Recall that in this situation the maximal split torus $A$ is given by $\{ a(t_1,t_2) = \diag(t_1,t_2,t_2^{-1},t_1^{-1}); \, t_i >0 \}$. We set $W:= \Span(u_2,u_2')$. Then 

\begin{itemize}

\item[(i)] $e(Q'_+) \simeq \h \times \R $ with trivial bundle structure. The collar neighborhood in $\overline{D}$ is given by $e(Q'_+) \times \{ a(t,t); \, t^2 >T \}$. 

\item[(ii)] $e(P') = N_{P'} \simeq W \simeq \R^2 $.  The collar neighborhood in $\overline{D}$ is given by $e(P') \times \{ a(t_1,t_2); t_1t_2>T, t_1/t_2 >T \}$. 

\item[(iii)] $e(Q'_-) \simeq \h \times \R  $ with trivial bundle structure. The collar neighborhood in $\overline{D}$ is given by $e(Q'_+) \times \{ a(t,t^{-1}); \, t^2 >T \}$. 

\end{itemize}

Furthermore, $e(Q'_+)$ and $e(Q'_-)$ are glued in $e(P')$ with the respective $\R$-fibers glued to the ``$x$-direction" of $\h$. 

  Now we consider the analogous picture for $\Orth(2,2)$. The faces $e(Q_{\pm})$ for the stabilizers $Q_{\pm}$ of the planes $E_\pm$ stay the same (with slightly different neighborhoods). But now there are three parabolics $P,P_+,P_-$ whose restriction to $\SO(2,2)$ is $P'$, and we blow up $e(P')$ by $e(P) \simeq e(P') \times \R_+$ and glue $e(P)$ to $e(Q_\pm)$ along $e(P_\pm)$. The blow-up variable in $\R_+$ in the neighborhood of $e(P')$ is given by $t_1/t_2$. We have

\begin{itemize}

\item[(i)]
$P$ is the stabilizer of the line $E_1 =\R u_1$. Then $e(P) = \{ a(1,t_2) \} \times W$ with collar neighborhood $e(P) \times \{a(t,1); \, t_2>T \}$. 

\item[(ii)]
$P_\pm$ are the stabilizers of the flag $\R u_1 \subset E_\pm$. Then $e(P_\pm) \simeq W$. Collar neighborhoods are given by $e(P_+) \times \{ a(t_1,t_2); \, t_1t_2 >T, t_2 >T \}$ and $e(P_-) \times \{ a(t_1,t_2); \, t_1t^{-1}_2 >T, t_2^{-1} >T \}$ respectively. 

\end{itemize}

Inside $e(P)\simeq \{ a(1,t_2) \} \times W$ one approaches $e(P_\pm)$ by $t_2 \to \infty$ and $t_2 \to 0$ resp..

\subsubsection{Nonexistence and existence of the restriction for the case of $\SO(2,2)$}
In this subsection we will explain why $\theta(\varphi_{2,0})$ does not extend to $\overline{X}_{p,p}$  if $\overline{X}_{p,p}$ is the {\it small} Borel compactification of $\SO(p,p)$. 

Namely, $\theta(\varphi_{2,0})$ does not extend to the $2$-torus $e'(P')$, where $P'$ is the stabilizer of the oriflamme $(E_+,E_-)$. The limit as we approach $e'(P')$ is undefined (it depends on the way we approach the corner). We have just seen that the corner $e'(P')$ is the intersection of the two maximal faces $e'(Q'_\pm)$, trivial circle bundles over quotients of the upper half plane. 

It suffices to study $\theta(\varphi_{2,0})(a(t_1,t_2)) = \sum_{y_1,y_2,y_2',y_1'} \varphi_{2,0}(t^{-1}_1y_1,t_2^{-1}y_2,t_2y_2,t_1y_1')$ as we go to the corner. Here $y_i,y_i'$ are the Witt coordinates of $V$. In this case the $2$-form $\theta(\varphi_{2,0})$ has four components. Three of the components go to zero as $\alpha_1 = t_1t_2$ and $\alpha_2= t_1/t_2$ go to infinity; essentially because $t_1 = \sqrt{\alpha_1\alpha_2} \to \infty$ we can apply the partial Fourier transform and Poisson summation argument from Section~\ref{globalressection} on the sum on $y_1$. We find that the limit
coincides (up to a constant) with 
 $\sum_{(y_2,y_2')}  
\widetilde{H}_2(t_2 ^{-1} y_2 + t_2 y'_2) 
e^{ - \pi (t_2^{-2} y_2^2 + t_2^2 (y'_2)^{2})}
\frac{dt_2}{t_2} \wedge (\tfrac{dw_2}{t_2} + t_2{dw_2'} )$. Here $w_2,w_2'$ are the variables for the $2$-torus $e'(P')$ realized as a quotient of $W = \R u_2 \oplus \R u_2'$. Now the resulting limit is supposed to be a $2$-form on the corner $e'(P')$, that is, a form in the coordinates $w_2,w_2'$ on the torus. However note that {the limit depends on $t_2$} (and also involves the coordinate differential $dt_2$).  Thus it depends on how we approach the boundary and consequently is not well defined. In particular, as claimed,
the form  $\theta(\varphi_{2,0})$ {does not extend to a well-defined $2$-form on the manifold with corners $\overline{X}$.

In the big Borel-Serre compactification the problems go away. For the face $e'(P)$, $t_2$ is the extra variable for $e'(P) = e'(P') \times \R_+$, and we obtain the above form as the limit as $t_1 \to \infty$. At the other faces $e(P_\pm)$ which as sets are again the $2$-torus $e'(P')$ but now are approached by $t_1/t_2, t_2 \to \infty$ resp. $t_1t_2, t_2^{-1} \to \infty$. Then the Poisson summation argument on the sum on $y_1,y_2$ resp $y_1,y_2'$ gives vanishing.

\end{document}